%% file: paper_HILS-RVRP.tex
\documentclass[a4papert,10pt]{article}

\usepackage{graphicx}          

\usepackage[latin1]{inputenc}
\usepackage{lscape}
\usepackage[lined,algoruled,linesnumbered,procnumbered,onelanguage]{algorithm2e}
\usepackage{enumerate}
\usepackage{enumitem}
\usepackage{setspace}
\usepackage{amssymb}
\usepackage{amsmath}
\usepackage{setspace}

\usepackage[hang]{subfigure}

\usepackage{longtable}
\usepackage{pdflscape}

\usepackage{url}

\usepackage{colortbl}
\usepackage{xcolor}
\definecolor{LightCyan}{rgb}{0.88,1,1}

\usepackage{multicol}
\usepackage{multirow}
\usepackage{rotating}
\usepackage{enumitem}

\usepackage[comma]{natbib}

\usepackage[pdftex,hidelinks]{hyperref}

\begin{document}
\sloppy
\onehalfspacing


\title{\textbf{A Hybrid Heuristic for a Broad Class of Vehicle Routing Problems with Heterogeneous Fleet}}

\author{
	\textbf{Puca Huachi Vaz Penna} \\ 
	\vspace{-2pt} 
	Dep. de Computa\c{c}\~{a}o, Universidade Federal de Ouro Preto - Ouro Preto, Brazil \\ 
	puca@iceb.ufop.br 
	\and
	\textbf{Anand Subramanian} \\
	\vspace{-2pt} 
	Dep. de Sistemas de Computa\c{c}\~{a}o, Universidade Federal da Para\'{i}ba, Jo\~{a}o Pessoa, Brazil \\ anand@ci.ufpb.br
	 \and
	\textbf{Luiz Satoru Ochi }\\
	\vspace{-2pt} 
	Universidade Federal Fluminense - Instituto de Computa\c{c}\~{a}o, Niter\'{o}i,  Brazil \\ satoru@ic.uff.br
	\and
	\textbf{Thibaut Vidal }\\
	\vspace{-2pt} 
	Pontifícia Universidade Católica do Rio de Janeiro,  Rio de Janeiro, Brazil \\ vidalt@inf.puc-rio.br
	\and
	\textbf{Christian Prins} \\
	\vspace{-2pt} 
	ICD-LOSI, UMR CNRS 6281, Université de Technologie de Troyes,  Troyes, France \\ christian.prins@utt.fr
}


\date{}

\maketitle

 \paragraph{Abstract.} We consider a family of Rich Vehicle Routing Problems (RVRP) which have the particularity to combine a heterogeneous fleet with other attributes, such as backhauls, multiple depots, split deliveries, site dependency, open routes, duration limits, and time windows.
To efficiently solve these problems, we propose a hybrid metaheuristic which combines an iterated local search with variable neighborhood descent, for solution improvement, and a set partitioning formulation, to exploit the memory of the past search. Moreover, we investigate a class of combined neighborhoods which jointly modify the sequences of visits and perform either heuristic or optimal reassignments of vehicles to routes.
To the best of our knowledge, this is the first unified approach for a large class of heterogeneous fleet RVRPs, capable of solving more than 12 problem variants. The efficiency of the algorithm is evaluated on 643 well-known benchmark instances, and 71.70\% of the best known solutions are either retrieved or improved. Moreover, the proposed metaheuristic, which can be considered as a matheuristic, produces high quality solutions with low standard deviation in comparison with previous methods. Finally, we observe that the use of combined neighborhoods does not lead to significant quality gains. Contrary to intuition, the computational effort seems better spent on more intensive route optimization rather than on more intelligent and frequent fleet re-assignments.




\section{Introduction}

The capacitated Vehicle Routing Problem (VRP) is one of the most studied problems in the field of combinatorial optimization. Since the seminal work of \cite{Dantzig1959}, many additional constraints, objectives and decision subsets, called problem attributes, have been combined with the classical version of the problem.
Such attributes include multiple depots, pickup and delivery, backhauls, heterogeneous fleets, time windows, among others. The reader is referred to \cite{Vidaletal-survey2013} for a recent survey and classification on the most common attributes adopted in the VRP literature.
Several variants that consider each of these attributes individually received a lot of attention over the past few years. However, in practical applications, many attributes tend to appear together, thus increasing the resolution challenges.  
Several recent articles have attempted to cope with this increasing variety of problems. The term \emph{rich}, in particular, has being widely adopted to describe VRP versions composed of multiple attributes.
Given the importance of solving real-world problems, unified solution methods capable of tackling many VRP variants are of high importance. This explains the recent trend in the development of this type of approach \citep{RopkePisinger2006, Subramanianetal2013a, DerigsVogel2014, Vidaletal2014}. 

One important aspect of practical VRP applications is the frequent use of a \emph{heterogeneous} fleet of vehicles \citep{Hoff2010}, with different capacities and operational costs. This type of VRPs was excluded from most unified frameworks available in the literature. For example, the framework of \cite{Vidaletal2014} considered VRPs with heterogeneous fleet, but only in the case where the fleet is unlimited (fleet size and mix VRP -- FSMVRP).

The contributions of this work are as follows.
\begin{itemize}[nosep]

	\item We propose the first unified algorithm designed to solve a broad class of Heterogeneous Fleet RVRPs (HFRVRPs), thus filling the methodological gap of previous works and extending the range of applications to richer and more challenging variants. 
	
	\item Our algorithm is capable of dealing with at least 10 distinct and well-known attributes, such as multiple depots, time windows, (mixed) backhauls, site dependency, and split deliveries. 
The attributes may be considered one at a time (classical variants) or simultaneously (rich variants), leading to a wide gamut of problems. In practice, hundreds of variants could be formed with those attributes. Obviously, for comparison purposes, we have decided to test our algorithm only on those variants 
where there has been substantial research and publicly available instances. 
	
	\item The algorithm generalizes the ILS-RVND-SP matheuristic of \cite{Subramanianetal2013a}, originally designed for vehicle routing problems with homogeneous fleet, and the alternative version of this matheuristic \citep{Subramanian2012a} that was developed for classical Heterogeneous Fleet VRPs (HFVRPs).  ILS-RVND-SP combines an Iterated Local Search (ILS -- \citealt{Lourenco2010}) with Randomized Variable Neighborhood Decent (RVND) and an integer programming-based optimization over a Set Partitioning (SP) formulation. The routes generated by the ILS-RVND heuristic \citep{Pennaetal2013} are used to create a pool of promising routes for the SP. The SP problem is then solved by a Mixed Integer Programming (MIP) solver, which interacts with the ILS-RVND during its execution. 
	
	\item The generalization includes, among other features, the addition of alternative constructive procedures, neighborhood structures and perturbation mechanisms in order to cope with attributes that were not considered in \cite{Subramanianetal2013a}.  We also introduce a novel perturbation scheme that exploits the heterogeneous characteristic of the fleet. Moreover, in contrast to \cite{Subramanian2012a, Subramanianetal2013a}, the proposed generalized algorithm now accepts infeasible solutions, which seem to be crucial for obtaining high quality solutions when dealing variants with time windows \citep{Vidaletal-survey2013,Vidaletal2013}. In addition, move evaluations for this class of variants are now performed via the efficient approach presented in \cite{Vidaletal2014}. Overall, this led to a successful algorithm, which produces solutions of high quality, and sometimes new best solutions, for 12 difficult VRP variants with heterogeneous fleet, 20 sets of benchmark instances and overall 643 test problems.
	
	\item Finally, we followed the reasonable belief that optimality gaps could be related to errors in fleet assignment rather than errors in routing. For this reason, we tested more advanced neighborhoods with combine moves on the sequence of visits with a re-optimization (either heuristic or exact) of the fleet assignment decisions. From our experiments, this approach did not lead to significant improvements. As this (negative) result seems very counter-intuitive at first glance, we believe that it deserves some discussion in the present paper.
\end{itemize}

The remainder of this paper is structured as follows. Section \ref{sec:problems} describes the problems under study and Section \ref{sec:RelatedWorks} reviews the main works related to RVRP and heterogeneous fleet variants. Section \ref{sec_ILSRVND} describes the proposed matheuristic as well as the combined neighborhoods with joint routing and assignment optimization. Section \ref{sec_Results} reports the computational results and establishes a comparison with the current literature. Section~\ref{sec_conclusoes} finally concludes.

\section{Class of Problems Considered} \label{sec:problems}

In what follows, we provide a formal description of the classical HFVRPs as well as of the main attributes considered in this work.

The HFVRP can be defined as follows: let $G = (V, A)$ be a directed graph where $V = \{0, 1, \ldots, n\}$ is a set composed of $n + 1$ vertices, and $A = \{(i, j): i, j \in V, i \neq j\}$ is the set of arcs. Vertex $0$ denotes the depot, where the vehicle fleet is located, while the set $V' = V\setminus \{0\}$ includes the remaining vertices which represent the $n$ customers. Each customer $i \in V'$ has a non-negative demand $q_i$. The fleet is composed by $K$ different types of vehicles, with $M = \{1, \ldots, K \}$. For each $k \in M$, there are $m_k$ available vehicles, each with a capacity~$Q_k$.
Every vehicle type is also associated with a fixed cost denoted by $f_k$. Finally, each arc $(i, j) \in A$ has a length $d_{ij}$ and its traversal cost for vehicle type $k$ is $c_{ij}^k = d_{ij} \times r_k$, where $r_k$ is a cost per distance unit, also called dependent cost or variable cost in the literature.
The objective is to determine a fleet composition as well as a set of routes, $R^k = (i_1, i_2, \ldots, i_{|R|})$, that minimize the sum of fixed and travel costs in such a way that: 
(a) every route $R^k$ starts and ends at the depot ($i_1 = i_{|R|} = 0$ and $\{i_2, \ldots, i_{|R|-1}\subseteq V^{\prime}\}$) and is associated with a vehicle type $k \in M$; 
(b) each customer belongs to exactly one route; 
(c) vehicle capacity is not exceeded;
(d) for each vehicle type $k$, the number of vehicles actually used does not exceed $m_k$.

The HFVRP is $\mathcal{NP}$-hard since it includes the classical VRP as a special case when all vehicles are identical. The problem was introduced by \cite{Golden1984} under the name Fleet Size and Mix (FSM), which is a variant that assumes an unlimited number of vehicles of each type, i.e., $m_k = +\infty, \forall k \in M$. Fifteen years later, \cite{Taillard1999} proposed the Heterogeneous Fixed Fleet VRP (HFFVRP), a variant in which the number of vehicles of each type is limited. 

Since the seminal works of \cite{Golden1984} and \cite{Taillard1999}, several HFVRP extensions considering well-known real-life VRP attributes were presented in the literature. We were able to deal with the following ones:

\begin{itemize}[nosep]
\item \textbf{Asymmetry (A)}: the costs of two opposite arcs may differ, i.e., $c_{ij}$ is not necessarily equal to $c_{ji}$, for $i \in V$ and $j \in V$. Although this is the case in cities with one-way streets, a vast majority of VRP articles consider undirected networks which lead to simpler local search procedures.
\item \textbf{Open (O) routes}: in open VRPs, the vehicle does not need not return to the depot after visiting the last customer, i.e., $c_{i0} = 0$ for all $i \in V$.
\item \textbf{Multiple depots (MD)}: more than one depot is available, but each vehicle must start and end at the same depot in a route. The number of vehicles per depot is usually limited.
\item \textbf{Multiple trips (MT)}: each vehicle may perform a sequence of successive trips called multitrip, often limited by a maximum length or duration.
\item \textbf{Backhauls (B)}: two different types of customers are considered, more precisely, linehaul and backhaul. The first type includes the customers with delivery demands, while the second one includes those with pickup demands. In this case, backhaul customers can only be visited after the last linehaul customer, and a route cannot be only composed of backhaul customers. The vehicle leaves the depot with a load that is equal to the sum of the delivery demands (lineheauls) in the route, and returns to the depot with a load that is equal to the sum of the pickup demands (backhauls).
\item \textbf{Mixed Backhauls (MB)}: similar to the previous case, but there are no constraints on the order in which linehaul  and backhaul customers should be visited. The load of a vehicle may increase or decrease along its route.
\item \textbf{Site dependency (SDep)}: some customers can only be visited by a subset of the existing  vehicles. This usually happens when a customer site limits its access to vehicles with some particular characteristics.
\item \textbf{Split Deliveries (SD)}: customers can be visited more than once and their demands can be split among different vehicles. In this case, it is necessary to decide on the amount of goods to be delivered to each customer by each vehicle.
\item \textbf{Time Windows (TW)}: a time window $[a_i, b_i]$ and a service time $s_i$ are defined for each vertex $i \in V$. The service of a customer should start within its time window and a vehicle is allowed to arrive at customer $i$ before $a_i$, but not after $b_i$. Note that waiting times are allowed.
\item \textbf{Route Duration (RD) constraints}:  there exists a limit $D$ on the duration of a route in terms of distance or time.
\end{itemize}

\section {Related Works}\label{sec:RelatedWorks} 


\textbf{The classical HFVRP} has been widely studied in the literature, as can be observed in the surveys of \cite{BaldacciBook2008}, \cite{Hoff2010}, \cite{Irnichetal2014} and \cite{Kocetal-review2015}. RVRPs have also been an object of active research interest. The reader is referred to \cite{Caceres-Cruzetal2014}, \cite{DerigsVogel2014} and \cite{Lahyanietal2015} for a comprehensive literature review on this topic. In this section we review some of the main works and milestones for the HFRVRPs.

As in most vehicle routing problems, the progress on solution methods went through several phases: from early constructive methods, towards local search-based heuristics, metaheuristics, and hybrid methods.
First, \cite{Golden1984} developed several heuristics for the FSMVRP, while \cite{Taillard1999} presented a column generation heuristic for the HFFVRP.
In the following years, nearly all classical metaheuristic frameworks have been considered:
\begin{itemize}[nosep]
\item \emph{Tabu Search (TS)} in \cite{Gendreau1999,Wassan2002,Lee2008,Brandao2009} and \cite{Brandao2011}; 
\item \emph{Variable Neighborhood Search (VNS)} in \cite{Imran2009}; 
\item \emph{Iterated Local Search (ILS)} in \cite{Pennaetal2013} and a hybrid ILS in \cite{Subramanianetal2013a};
\item \emph{Threshold Accepting (TA)} in \cite{Tarantilis2003,Tarantilis2004};
\item \emph{Record-to-Record (RTR)} in \cite{Li2007};
\item \emph{Adaptive Memory Programming (AMP)} in \cite{RochatTaillard1995,Li2010};
\item \emph{Multi-Start (MS) / Evolutionary Local Search (ELS)} in \cite{Duhameletal2011,Duhameletal2013};
\item \emph{Hybrid Genetic Algorithms (HGA)} in \cite{Ochi1998a,Ochi1998b,Lima2004,Prins2009,Liu2009} and \cite{Vidaletal2014}.
\end{itemize}

Not all these methods were equally successful, and their relative performance can even vary between different benchmark sets. After years of research on metaheuristics for the HFVRP, it was impossible to conclude on a \emph{more suitable} metaheuristic framework. The common viewpoint remains that success cannot be attributed to one method's name, but rather to its specific components, efficient neighborhood structures, a proper balance between intensification and diversification strategies \citep{Blum2003,Vidaletal-survey2013}, as well as a good distributions of the search effort dedicated to the sequencing and customer-to-vehicle assignment decision sets.
Finally, the research on exact methods has also progressed over the last few years. The current state-of-the-art approaches \citep{Choi2007,Pessoaetal2009,Baldaccietal2009, Baldacci2010, Baldacci2010a} can consistently solve instances with 75 customers in a few minutes, as well as some instances of a few 100 customers. However, this performance is still insufficient for many practical applications. \\

\textbf{Rich variants of the HFVRP} can include a large variety of attributes.
From early works on separate problems, which tends to reduce our ability to compare method performances, the literature has evolved towards unified solution frameworks with the potential to be applied and compared on a wide gamut of problems. 
With this aim, some of the most successful algorithms include the Unified Tabu Search (UTS) of \cite{Cordeauetal2001}, the Adaptive Large Neighborhood Search (ALNS) of \cite{PisingerRopke2007}, the Iterated Local Search (ILS) of \cite{Subramanianetal2013a} and the Unified Hybrid Genetic Algorithm (UHGS) of \cite{Vidaletal2014}. 
The success of UTS can be mainly explained by very simple neighborhoods and diversity mechanisms.
To achieve high quality solutions for a variety of problems, ALNS uses a family or destruction and reconstruction procedures, along with an adaptive selection strategy.
ILS exploits again simple neighborhoods along with mathematical programming over a set partitioning formulation.
Finally, UHGS relies on generic route evaluation operators, a giant-tour solution representation (see \citealt{Prinsetal2014} for a survey on VRP algorithms based on this representation), as well as advanced diversity management mechanisms that promote different and good solutions for survival.


Table \ref{tbl:literature} summarizes a large variety of works related to rich HFVRPs. The table is organized by method family (constructive heuristic or local search; neighborhood-centered; population-based; mathematical programming). The first column, \textbf{Authors}, shows the authors names. The second, named \textbf{Approach} list the method used by the authors to solve the problem. Next, columns \textbf{Fleet} and \textbf{Costs} present the main attributes of HFVRP, the fleet size and the associated vehicles costs, respectively. Finally, the last set of columns, \textbf{Additional Attributes}, shows the multi-attributes characteristic of variant studied by the authors, where \textbf{MD}, \textbf{TW}, \textbf{RD}, \textbf{SDep}, \textbf{SD}, \textbf{MT}, \textbf{Op} and \textbf{Others} are, respectively, Multi-Depot, Time Windows, Route Duration, Site Dependency, Split Delivery Multi-Trip, Open and others attributes.

\begin{table}[!htbp]

\def\arraystretch{1.22}
\setlength{\tabcolsep}{0.6mm}
\caption{HFRVRP related works and  attributes}
\scalebox{0.90} {
\begin{tabular}{@{\extracolsep{\fill}}l@{\hspace*{0.3cm}}ll@{\hspace*{0.4cm}}cclccclcccccccccc}
\noalign{\smallskip}\hline\noalign{\smallskip}
 & &  & \multicolumn{2}{l}{Fleet}  &  & \multicolumn{3}{l}{Costs} &  & \multicolumn{8}{l}{Additional Attributes}  \\ 
  \cline{4-5} \cline{7-9} \cline{11-18}\noalign{\smallskip}
 \multicolumn{1}{c}{Authors} &  & Approach(es) & U & L &  & F & V & FV &  & MD & RD & TW & SDep & SD & MT & Op & Others \\ 
\noalign{\smallskip}\hline\noalign{\smallskip}
\cite{LiuShen1999} &  & CH & $\checkmark$ &  &  & $\checkmark$ &  &  &  &  &  & $\checkmark$ &  &  &  & &\\ 
\cite{DellAmico2007} &  & CH+RR & $\checkmark$ &  &  & $\checkmark$ &  &  &  &  &  & $\checkmark$ &  &  &  & & \\ 
\cite{DondoCerda2007} &  & CH &  & $\checkmark$ &  &  &  & $\checkmark$ &  & $\checkmark$ &  & $\checkmark$ &  &  &  & &\\ 
\cite{Caceres-Cruzetal2013} & & CH & & $\checkmark$ &  &  & $\checkmark$ &  &  &  &  $\checkmark$ &  &  &  & $\checkmark$ & $\checkmark$ & A \\ 
\cite{Caceres-Cruz2014} &  & CH &  & $\checkmark$ &  &  & $\checkmark$ &  &  &  & $\checkmark$ &  &  &  & $\checkmark$ & & \\ 
\cite{Prins2002} &  & CH+LS/TS &  & $\checkmark$ &  &  & $\checkmark$ &  &  & \multicolumn{1}{l}{} & $\checkmark$ &  &  &  & $\checkmark$ & & \\ 
\cite{Salhi1997} &  & LS & $\checkmark$ &  &  &  &  & $\checkmark$ &  & $\checkmark$ &  &  &  &  &  & &\\ 

\cite{CordeauLaporte2001} &  & TS &  & $\checkmark$ &  &  &  &  &  &  & $\checkmark$ &  & $\checkmark$ &  &  & &\\ 
\cite{Cordeauetal2004} &  & TS &  & $\checkmark$ &  &  &  &  &  &  & $\checkmark$ & $\checkmark$ & $\checkmark$ &  & & & \\ 
\cite{Paraskevopoulos2008} &  & TS & $\checkmark$ &  &  & $\checkmark$ &  &  &  &  &  & $\checkmark$ &  &  &  & & \\ 
\cite{Lietal2012} &  & TS+AMP &  & $\checkmark$ &  &  &  & $\checkmark$ &  &  &  &  &  &  & & $\checkmark$ & \\ 
\cite{Yousefikhoshbakht2014} &  & TS &  & $\checkmark$ &  &  &  & $\checkmark$ &  &  &  &  &  &  & & $\checkmark$ &  \\ 
\cite{CordeauMaischberger2012} &  & ILS+TS &  & $\checkmark$ &  &  &  &  &  &  & $\checkmark$ & $\checkmark$ & $\checkmark$ &  &  & & \\ 

\cite{Tavakkoli-Moghaddam2007} &  & SA &  & $\checkmark$ &  &  & $\checkmark$ &  &  &  &  &  &  & $\checkmark$ & & &  \\ 
\cite{Braysyetal2008} &  & SA & $\checkmark$ &  &  & $\checkmark$ &  &  &  &  &  & $\checkmark$ &  &  &  & & \\ 
\cite{Braysyetal2009} &  & TA+GLS & $\checkmark$ &  &  & $\checkmark$ &  &  &  &  &  & $\checkmark$ &  &  &  & & \\ 
\cite{PisingerRopke2007} &  & ALNS &  & $\checkmark$ &  &  &  &  &  &  & $\checkmark$ &  & $\checkmark$ &  &  & &\\ 
\cite{Amorimetal2014} &  & ALNS & $\checkmark$ &  &  &  &  & $\checkmark$ &  &  &  & $\checkmark$ & $\checkmark$ &  & & & GI \\ 
\cite{Mancini2015} &  & ALNS & $\checkmark$ & \multicolumn{1}{l}{} &  &  & $\checkmark$ &  &  & $\checkmark$ & $\checkmark$ &  & $\checkmark$ &  & & & P,PI\\

\cite{RepoussisTarantilis2010}  &  & AMP & $\checkmark$ &  &  & $\checkmark$ &  &  &  &  &  & $\checkmark$ &  &  &  & & \\ 
\cite{Tutuncu2010} &  & GRASP+AMP &  & $\checkmark$ &  &  & $\checkmark$ &  &  &  &  &  &  &  & & & B \\ 
\cite{Duhameletal2011} &  & GRASP+ELS &  & $\checkmark$ &  &  &  & $\checkmark$ &  &  &  &  &  &  & & &  \\ 
\cite{Duhameletal2013} &  & GRASP+ ELS &  & $\checkmark$ &  &  &  & $\checkmark$ &  &  &  &  &  &  & & &  \\ 
\cite{Mar-Ortizetal2013} &  & GRASP &  & $\checkmark$ &  &  &  &  &  &  & $\checkmark$ & $\checkmark$ & $\checkmark$ & $\checkmark$ & $\checkmark$ & & \\ 
\cite{GoelGruhun2008} &  & VNS/LNS &  & $\checkmark$ &  &  & $\checkmark$ &  &  &  &  & $\checkmark$ &  &  & & $\checkmark$ & PD,MC \\ 

\cite{Salhietal2014} &  & VNS & $\checkmark$ &  &  &  &  & $\checkmark$ &  & $\checkmark$ &  &  &  &  &  & & \\ 
\cite{Armasetal2015} &  & VNS & \multicolumn{1}{l}{} & $\checkmark$ &  &  & $\checkmark$ &  &  &  & $\checkmark$ & $\checkmark$ & $\checkmark$ &  & & & CPr \\ 
\cite{Armas2015} &  & VNS & \multicolumn{1}{l}{} & $\checkmark$ &  &  & $\checkmark$ &  &  &  & $\checkmark$ & $\checkmark$ & $\checkmark$ &  & & & CPr,DR \\
\cite{Dominguez2016} &  & MS & $\checkmark$ &  &  &  &  & $\checkmark$ &  &  &  &  &  &  &  & & 2L\\ 

\noalign{\smallskip}\hline\noalign{\smallskip}
\cite{CruzReyesetal2007} &  & ACO &  & $\checkmark$ &  &  &  &  &  & $\checkmark$ &  & $\checkmark$ & $\checkmark$ & $\checkmark$ & & & VS \\ 
\cite{Pellegrinietal2007} &  & ACO &  & $\checkmark$ &  &  &  & $\checkmark$ &  &  & $\checkmark$ & $\checkmark$ &  &  & & & P \\ 
\cite{Belmecherietal2013} &  & ACO/PSO &  & $\checkmark$ &  &  & $\checkmark$ &  &  &  &  & $\checkmark$ &  &  & & & MB \\ 
\cite{Belfiore2009} &  & SS &  & $\checkmark$ &  &  &  & $\checkmark$ &  &  &  & $\checkmark$ &  & $\checkmark$ &  & & \\ 
\cite{Vidaletal2014} &  & GA & $\checkmark$ &  &  & $\checkmark$ &  &  &  & $\checkmark$ & $\checkmark$ & $\checkmark$ & $\checkmark$ &  &  & & \\ 
\cite{BerghidaBoukra2015} &  & EA &  & $\checkmark$ &  &  & $\checkmark$ &  &  &  &  & $\checkmark$ &  &  & & & MB \\ 
\cite{Kocetal2015} &  & EA & $\checkmark$ & $\checkmark$ &  &  $\checkmark$ &  &  &  &  &  & $\checkmark$  &  &  &  & & \\ 
\cite{Yao2016} &  & PSO &  & $\checkmark$ &  &  &  & $\checkmark$ &  & $\checkmark$ &  &  &  &  &  & & CD \\ 

\noalign{\smallskip}\hline\noalign{\smallskip}
\cite{Cesellietal2009} &  & CG &  & $\checkmark$ &  &  &  & $\checkmark$ &  & $\checkmark$ & $\checkmark$ & $\checkmark$ & $\checkmark$ & $\checkmark$ & & & WT,EC,PI \\ 
\cite{Goel2010} &  & CG &  & $\checkmark$ &  &  & $\checkmark$ &  &  &  &  & $\checkmark$ &  &  & & $\checkmark$ & PD,MC \\ 
\cite{Bettinellietal2011} &  & BCP/CG & \multicolumn{1}{l}{} & $\checkmark$ &  & $\checkmark$ &  &  &  & $\checkmark$ & $\checkmark$ & $\checkmark$ &  &  &  & & \\ 
\cite{Ozfirat2010} &  & CP & $\checkmark$ & $\checkmark$ &  & $\checkmark$ &  &  &  &  &  &  &  & $\checkmark$ & & &  \\ 
\cite{RieckZimmermann2010} &  & Math. Model &  & $\checkmark$ &  &  & $\checkmark$ &  &  &  &  & $\checkmark$ &  &  & & & SPD,WT,D \\ 
\cite{Salhietal2013} &  & Math. Model/LS & $\checkmark$ &  &  & $\checkmark$ &  &  &  &  &  &  &  &  & & & B \\ 

\noalign{\smallskip}\hline\noalign{\smallskip}
HILS-RVRP (this work)&  & ILS+SP & $\checkmark$ & $\checkmark$ &  & $\checkmark$ & $\checkmark$ & $\checkmark$ &  & $\checkmark$ & $\checkmark$ & $\checkmark$ & $\checkmark$ & $\checkmark$ & $\checkmark$ & $\checkmark$ & B, MB \\ 
\noalign{\smallskip}\hline\noalign{\smallskip}
  \multicolumn{18}{l}{{\scriptsize 
  2L: Two-dimensional Loading;
  A: Asymmetric;
  B: Backhauls;
  CD: Collection Depot;
  Cpr: Customer Priority;
  D: Docking;\vspace*{-0.2cm}
  }} \\
  \multicolumn{18}{l}{{\scriptsize
  DR: Dynamic Request;
  EC: External Courier;
  MB: Mixed Backhauls;
  MC: Multi-Dimensional Capacity;
  P: Periodic;
  PD: Pickup and Delivery; \vspace*{-0.2cm}
  }} \\
  \multicolumn{18}{l}{{\scriptsize
  PI: Products Incompatibility;
  SPD: Simultaneous Pickup and Delivery;
  VS: Vehicle Scheduling;
  WT: Working Time Regulations.
  }}
\label{tbl:literature}
\end{tabular}
}
\end{table}

In general, most research emphasis has been focused, in recent years, on metaheuristic principles and new problems rather than studies on neighborhoods and elementary building blocks of the methods. We aim to contribute to this latter point, in order to progress towards unified methods and simple concepts which can be efficiently applied to a good variety of HFRVRPs. Finally, even after nearly a hundred articles dealing with HFVRPs, some important questions remain open:
\begin{itemize}[nosep]
	\item What is a good neighborhood for HFFVRP?
	\item Is it beneficial or not to revise frequently vehicle-to-route assignments during computations?
\end{itemize}
This is a second gap in the literature on which we aim to contribute.

\section{HILS-RVRP Algorithm}
\label{sec_ILSRVND}

The proposed unified hybrid heuristic, called HILS-RVRP, extends the algorithms of  \cite{Pennaetal2013} and \cite{Subramanian2012a,Subramanianetal2013a}. It combines a multi-start ILS \citep{Lourenco2010} with mathematical programming over a SP formulation.



Algorithm \ref{alg:HILS-RVRP} presents the pseudocode of HILS-RVRP. This matheuristic performs $I_{MS}$ restarts \mbox{(Lines \ref{hhurvfh-line-beginMaxIter}--\ref{hhurvfh-line-endMaxIter})}.
At each restart, an ILS-RVND heuristic is responsible for improving an initial solution (Lines \ref{hhurvfh-line-InitialSolution}--\ref{hhurvfh-line-IILS}), as well as populating a pool of routes ($r_{pool}$) associated with locally optimal solutions. The pool is then used to build a restricted SP model, which is solved by a MIP solver. 
Two strategies are proposed to exploit the SP model. If the instance is considered to be of sufficient size ($n \geq n_{large}$), then the algorithm calls the SP approach after each restart (Line \ref{hhurvfh-line-endResolvePC-ILS}). Otherwise, if the problem is of smaller size, then the SP model is solved only once at the end of the last restart. As in \cite{Subramanianetal2013a}, we have assumed that $n_{large} = 150$.


\begin{algorithm}[!ht]
	\caption{HILS-RVRP($I_{MS}$, $I_{ILS}$,  $T_{max}$, $n_{large}$, $RGap_{max}$)}
	\label{alg:HILS-RVRP}
	\DontPrintSemicolon
	\Begin{
		\SetAlgoVlined
		Initialize fleet \; \nllabel{alg:hhurvfh-line-initfleet}
		$v \leftarrow$ total number of vehicles\; {\nllabel{hhurvfh-line-vehicles}}
		$f(s^{*}) \leftarrow \infty$ \;
		$r_{pool} \leftarrow\ \emptyset$ \;
		\For {$i \leftarrow 1$ \textbf{to} $I_{MS}$}
		{ \nllabel{hhurvfh-line-beginMaxIter}
			$s \leftarrow$ GenerateInitialSolution($v$) \nllabel{hhurvfh-line-InitialSolution}\;
			$[s^{*\prime}, r_{pool}] \leftarrow$ ILS-RVND($s$, $I_{ILS}$, $r_{pool}$) \nllabel{hhurvfh-line-IILS}\;
			
			\If {\textnormal{(}$n \geq n_{large}$ \emph{\textbf{ or }} $i = I_{MS}$\textnormal{)} }
			{\nllabel{hhurvfh-line-beginResolvePC-ILS}
				$[s^{*\prime}, r_{pool}] \leftarrow$ SolveSP($r_{pool}$, $s^{*\prime}$, $I_{ILS}$, $T_{max}$, $RGap_{max}$) \nllabel{hhurvfh-line-endResolvePC-ILS}\;
			}
			
			\If {\textnormal{(}$f(s^{*\prime}) < f(s^{*})$\textnormal{)}}
			{ \nllabel{hhurvfh-line-UpdPerturb}
				$s^{*} \leftarrow$ $s^{*\prime}$\;
			}
		} \nllabel{hhurvfh-line-endMaxIter}
		\Return $s^{*}$ \;
	}
\end{algorithm}


	The pseudocode of the ILS-RVND sub-procedure is depicted in Algorithm \ref{alg:HHURVFH}. This heuristic iteratively performs local search and perturbation operations until a maximum number of consecutive iterations without improvements of the best current solution ($I_{ILS}$) is achieved (lines \ref{hhurvfh-line-beginI_{ILS}}--\ref{hhurvfh-line-endI_{ILS}}). The pool of routes is updated after each call to the local search procedure (Lines \ref{hhurvfh-line-Inicializar_{pool}} and \ref{hhurvfh-line-Adicionar_{pool}}).  In order to avoid $r_{pool}$ to grow arbitrarily large, the algorithm only adds to the pool those routes associated with solutions whose gap with respect to the best current solution is relatively small, as described in \cite{Subramanianetal2013a}. 
	
	
	\begin{algorithm}[!ht]
		\caption{ILS-RVND($s$, $I_{ILS}$, $r_{pool}$)}
		\label{alg:HHURVFH}
		\DontPrintSemicolon
		\Begin{
			\SetAlgoVlined
			$s^{*\prime} \leftarrow$ LocalSearch($s$) \nllabel{hhurvfh-line-RVND-out}\;
			AddTemporaryRoutes($r_{pool}$, $s^{*\prime}$, $f(s^{*})$) \nllabel{hhurvfh-line-Inicializar_{pool}}\;
			$iter_{ILS} \leftarrow$ 0\;
			\While {\textnormal{(}$iter_{ILS} \leq I_{ILS}$\textnormal{)}}
			{ \nllabel{hhurvfh-line-beginI_{ILS}}
				$s^{\prime} \leftarrow$ Perturbation($s^{*\prime}$) \nllabel{hhurvfh-line-Perturbation}\;
				$s^{\prime\prime} \leftarrow$ LocalSearch($s^{\prime}$) \nllabel{hhurvfh-line-RVND}\;
				AddTemporaryRoutes($r_{pool}$, $s^{\prime\prime}$, $f(s^{*})$) \nllabel{hhurvfh-line-Adicionar_{pool}}\;
				\If {\textnormal{(}$f(s^{\prime\prime}) < f(s^{*\prime})$\textnormal{)}}
				{ \nllabel{hhurvfh-line-UpdSolucao}
					$s^{*\prime} \leftarrow$ $s^{\prime\prime}$\;
					$iter_{ILS} \leftarrow$ 0 \nllabel{hhurvfh-line-endUpdSolucao}\;
				}
				$iter_{ILS} \leftarrow$ $iter_{ILS} + 1$\;
			} \nllabel{hhurvfh-line-endI_{ILS}}
			\Return [$s^{*\prime}, r_{pool}$]\;
		}
	\end{algorithm}

Algorithm \ref{alg:HHURVFH.ResolvedorPC} describes the SolveSP procedure. The routes associated with the local optima of each iteration are permanently stored in the pool (Lines \ref{line-hhurvfh_PC-AddPermanenteRoutes} and \ref{hhurvfh-line-endUpdSolucaoMIP1}), while the remaining ones are treated as temporary routes that are removed at each iteration (Line \ref{line-hhurvfh_PC-RemoveTempRoutes}).
A time limit $T_{max}$ is imposed to the solver to avoid prohibitively large CPU time. During the resolution of the SP, the solver may find a new best integer solution. In that case, it calls the ILS-RVND to improve it (Line \ref{line-hhurvfh_PC-mipSolver}). If this improvement leads to a solution whose value is better than the current lower bound of the SP model, then the solver execution is naturally interrupted, otherwise the value of the solution is used as a cut-off bound.
Finally, when solving FSM variants, this model can be further modified in case the gap between the linear relaxation of the the root node and the incumbent solution exceeds a given input value ($RGap_{max}$). As thoroughly detailed in \cite{Subramanian2012a}, the modification consists of adding constraints to forbid the vehicle fleet associated with the best current solution to be changed. This may remove potential improved solutions from the integer linear program, but decreases the root gap and lets the model be more computationally tractable.

\begin{algorithm}[!ht]
	\caption{SolveSP($r_{pool}$, $s^{*}$, $I_{ILS}$, $T_{max}$, $RGap_{max}$)}
	\label{alg:HHURVFH.ResolvedorPC}
	\DontPrintSemicolon
	\Begin{
		\SetAlgoVlined
		AddPermanentRoutes($r_{pool}$, $s^{*}$) \label{line-hhurvfh_PC-AddPermanenteRoutes}\;
		$improvement \leftarrow$ true \;
		\While {\textnormal{(}$improvement$\textnormal{)}}
		{\label{line-hhurvfh_PC-beginLoop}
			$SP_{model} \leftarrow $ CreateSPModel($r_{pool}$, $v$) \label{line-hhurvfh_PC-createSP} \;
			$s^{\prime} \leftarrow$ MIPSolver($SP_{model}, s^{*}, T_{max}$,  $RGap_{max}$, ILS-RVND$(s^{*}, I_{ILS}, r_{pool}))$\label{line-hhurvfh_PC-mipSolver} \;
			\eIf {\textnormal{(}$f(s^{\prime}) < f(s^{*})$\textnormal{)}}
			{ \nllabel{hhurvfh-line-UpdSolucaoMIP1}
				$s^{*} \leftarrow$ $s^{\prime}$\;
				AddPermanentRoutes($r_{pool}$, $s^{*}$) \nllabel{hhurvfh-line-endUpdSolucaoMIP1}\;
			} 
			{
				$improvement \leftarrow$ false\;
			}
		}\label{line-hhurvfh_PC-endLoop}

				RemoveTemporaryRoutes($r_{pool}$) \label{line-hhurvfh_PC-RemoveTempRoutes}\;

		\Return [$s^{*}, r_{pool}$]\;
	}
\end{algorithm}


Our HILS-RVND also significantly differs from the previously cited works, as it has been extended to a much wider family of VRP variants with heterogeneous fleet and integrates several key improvements: (i) the ability of handling infeasible solutions with respect to time constraints; (ii) the incorporation of new neighborhoods and perturbations to manage multiple-depots, backhauls, split and time-windows attributes; and (iii) the addition of pre-processing phases and auxiliary data structures for efficient move evaluations. We finally investigate, within this method, a class of large neighborhoods which aim to jointly change the sequence of visits and optimize the assignment of vehicles to routes. 
The next subsections now provide a detailed description of each component of the method.

\subsection{Constructive Procedure}
\label{ConstructiveProcedure}


Initial solutions are generated via a simple insertion heuristic. Firstly, each route is filled with a random customer, and the remaining ones are inserted either according to (i) a nearest insertion criterion, or (ii) a modified cheapest insertion criterion which promotes the insertion of customers located far from the depot. At each restart, one of these criteria is randomly selected. Infeasible solutions are incorporated and penalized in accordance with the characteristics of the problem as follows:



\paragraph{\textbf{Fixed Fleet.}} 
When it is no longer possible to perform a feasible insertion of an unrouted customer due to fleet capacity, an extra vehicle with large fixed/dependent costs and capacity is added to the partial solution to accommodate the remaining customers. The routes associated with this extra vehicle tend to be emptied during the local search, leading to a feasible solution. Once this happens, this additional vehicle is eliminated from the solution.

\paragraph{\textbf{Unlimited fleet.}}
In this case, it is always possible to generate feasible solutions with respect to the vehicle capacity because there is no limit on the number of vehicles of each type. Once a complete initial solution is generated, a vehicle associated with each type is added to the solution so as to allow for a possible fleet resizing during the local search.

\paragraph{\textbf{Multiple Depots.}} 
The same rationale as in the previous cases is used when multiple-depots are considered, but an extra vehicle is added for each depot.

\paragraph{\textbf{Backhauls.}}
As specified by the problem, infeasible solutions are avoided by not allowing backhaul customers to be inserted before linehaul customers in the route. The insertion of backhaul customers in empty routes is also forbidden.
This is done via a simple modification of the distance matrix,  by setting a large cost to those arcs which connect the depot to one backhaul customer, or a backhaul customer to a linehaul customer.


\paragraph{\textbf{Site-Dependencies.}}

Customers are only allowed to be inserted in compatible routes. The cheapest insertion criterion also takes into account vehicle restrictions of the customers by including an insertion incentive that is inversely proportional to the number of vehicle types that a customer can be visited by. 
Moreover, a similar policy as the one used for the fixed fleet case is adopted when feasible insertions are no longer possible. In this case, the extra vehicle does not consider site-dependency constraints.


\paragraph{\textbf{Split Deliveries.}}
Splits are not allowed during the construction phase.

\paragraph{\textbf{Time Windows.}}
Time-window constraints are ignored during the construction phase.

\subsection{Local Search}
\label{LocalSearch}

The local search is performed by a procedure based on Randomized Variable Neighborhood Descent (RVND). Inter-route neighborhood classes are explored in random order as in a VND \citep{Hansenetal2010}. Intra-route neighborhoods are further applied to those routes that have been modified by one of the inter-route neighborhoods.

\subsubsection{Inter-Route Neighborhood Structures}
\label{sec:Inter-Route}

The method relies on several inter-route neighborhoods, described in the following. As indicated below, a few neighborhood structures are only applied for some specific attributes.

\paragraph{\textbf{General neighborhoods.}}

A set of seven inter-route neighborhood structures were adopted for all variants, namely: (i) \textsc{Shift(1,0)}; (ii) \textsc{Shift(2,0)}; 
(iii) \textsc{Swap(1,1)}; (iv) \textsc{Swap(2,1)}; (v) \textsc{Swap(2,2)}; (vi) \textsc{2-opt*}; and (vii) \textsc{$k$-Shift}. The first two consists of transferring one and two customers, respectively, from one route to another one. Neighborhoods (iii), (iv) and (v) consists of interchanging customers between two routes. For example, \textsc{Swap}(2,1) interchanges two consecutive customers from one route with one customer from another one. Neighborhood (vi) is an inter-route version of the classical \textsc{2-opt} neighborhood. Finally, neighborhood (vii) consists of moving $k$ consecutive customers from one route to the end of another one. 
 
\paragraph{\textbf{Multi-depot neighborhoods.}}
 

In the presence of multiple depots, two additional neighborhoods are considered: \textsc{ShiftDepot} and \textsc{SwapDepot}. The first one transfers an entire route from one depot to another one, whereas the second interchanges routes between two different depots. 
 
\paragraph{\textbf{Split-delivery neighborhoods.}}


When customer demands are allowed to be split, four additional neighborhoods are included in the RVND, namely: (i) \textsc{Swap(1,1)*}; (ii) \textsc{Swap(2,1)*}; (iii) \textsc{RouteAddition}; and (iv) \textsc{$k$-Split}. The first two were proposed by \citet{Boudia2007}, and they generalize the neighborhoods \textsc{Swap(1,1)} and \textsc{Swap(2,1)} by changing the quantities to be delivered to the customers that are moved. The third one was proposed by \citet{Dror1990}. It consists in adding an extra route by combining the segments of two routes that contain a customer that is visited more than once. The last neighborhood removes a customer from the solution and inserts it back using a procedure called \textsc{SplitReinsertion} \citep{Boudia2007}. 
A detailed description of these four neighborhoods can be found in \cite{Silvaetal2015}. \\

All neighborhoods are examined in an exhaustive fashion, that is, all possible moves are evaluated and the best improving move is applied.


\subsubsection{Intra-Route Neighborhood Structures}
\label{Intra-Route}

Five well-known intra-route neighborhood structures were adopted, namely: \textsc{Reinsertion}, \textsc{Or-opt} (with two and three adjacent customers), \textsc{2-opt} and \textsc{Exchange}. Such neighborhoods were also implemented using a RVND procedure, and they are called every time a solution is modified during the intra-route search.


\subsection{Perturbation Mechanisms}
\label{Perturbation}


Five perturbation mechanisms were adopted. They may generate infeasible solutions, but only with respect to time window constraints. Three of them are used in all variants, namely, \textsc{Multiple-Swap(1,1)}, \textsc{Multiple-Shift(1,1)} and \textsc{Split}. The first two consists of applying random \textsc{Swap(1,1)} and \textsc{Shift(1,1)} moves consecutively, respectively. Note that the \textsc{Shift(1,1)} operator moves one customer from route $R_{1}$ to route $R_{2}$ and vice-versa, but not necessarily interchanging their positions as in \textsc{Swap(1,1)}. The third perturbation divides a route into smaller ones as described in \cite{Pennaetal2013}. 

The fourth mechanism is only applied when the deliveries are allowed to be split. Such mechanism, called \textsc{Multiple-k-Split}, was proposed by \cite{Silvaetal2015} and consists in applying the \textsc{k-Split} operator several times at random. The last \textsc{Merge} operator, introduced in this work, selects a route whose capacity is smaller than the largest one and attempts to merge it with another route that is selected using a criterion similar to classical savings heuristics of \cite{Clarke1964}.


\subsection{Efficient Move Evaluations}\label{sec:MoveEval}

One important component of any local search based heuristic is the ability to efficiently perform move evaluations. This evaluation can be done in $\mathcal{O}(1)$ time for most VRPs, e.g., when the objective function is to minimize the total distance (or travel time) and also when penalties due to constraint violations are not incurred, that is, only feasible solutions are considered during the search. In our case, since we accept and penalize infeasible solutions with respect to time windows \citep{Vidaletal2013,Vidal2013a}, a more sophisticated approach must be adopted to compute the cost of moves in constant time.


To evaluate moves in the presence of time-window infeasible solutions, we apply the methodology of \cite{Vidaletal2013}. This technique has been designed for a specific relaxation where one is penalized for a \emph{return in time} (time warp) rather than for a late arrival \citep{Nagata2010}. Any classical move (e.g., \textsc{Swap} or \textsc{Shift}) produces a pair of routes which result from the concatenation of a bounded number of sequences of consecutive visits from the incumbent solution. By preprocessing and updating some meaningful information on sequences of the incumbent solution, it is possible to speed-up the subsequent evaluation of the moves. In our context, any sequence $\sigma$ will be characterized by six values: distance \mbox{$DIST(\sigma)$},
demand \mbox{$Q(\sigma$)},
duration \mbox{$D(\sigma$)},
earliest visit \mbox{$E(\sigma$)}, 
latest visit \mbox{$L(\sigma$)},
and return in time \mbox{$TW(\sigma$)}.
For a sequence $\sigma_0$ containing only one customer~$i$,
\mbox{$DIST(\sigma_0$) = $Q(\sigma_0$) = $TW(\sigma_0$) = $0$},
 \mbox{$D(\sigma_0$) = $s_i$},
 \mbox{$E(\sigma_0$) = $a_i$}, and
 \mbox{$L(\sigma_0$) = $b_i$}. Let $\sigma_1 = (\sigma_{1(i)}, \ldots, \sigma_{1(j)})$ and $\sigma_2=(\sigma_{2(i)}, \ldots, \sigma_{2(j)})$ be two sequences, then the sequence corresponding to the concatenation of $\sigma_1 \oplus \sigma_2$ (where the first customer of $\sigma_2$ is visited immediately after the last customer of $\sigma_1$) can be evaluated as:
\begin{align}
&DIST(\sigma_1 \oplus \sigma_2) = DIST(\sigma_1) + d_{\sigma_1(|\sigma_1|)\sigma_2(1)} + DIST(\sigma_2) \\
&Q(\sigma_1 \oplus \sigma_2) = Q(\sigma_1) + Q(\sigma_2) \\
&D(\sigma_1 \oplus \sigma_2) = D(\sigma_1) + d_{\sigma_1(|\sigma_1|)\sigma_2(1)} + D(\sigma_2) + \Delta_{WT} \label{eqn:DuracaoTW}\\
&E(\sigma_1 \oplus \sigma_2) = \max \{ E(\sigma_2) - \Delta, E(\sigma_1)\}  - \Delta_{WT} \\
&L(\sigma_1 \oplus \sigma_2) = \min \{ L(\sigma_2) - \Delta, L(\sigma_1)) \} + \Delta_{TW} \\
&TW(\sigma_1 \oplus \sigma_2) = TW(\sigma_1) + TW(\sigma_2) + \Delta_{TW}, \\
& \text{with } \hspace*{0.3cm} \Delta = D(\sigma_1) - TW(\sigma_1) + d_{\sigma_1(|\sigma_1|)\sigma_2(1)} \\
&\Delta_{WT} = \max \{ E(\sigma_2) - \Delta - L(\sigma_1), 0\} \\
&\Delta_{TW} = \max \{ E(\sigma_1) + \Delta - L(\sigma_2), 0 \}. \label{eqn:DeltaPJT}
\end{align}

By using Eqs. \eqref{eqn:DuracaoTW}--\eqref{eqn:DeltaPJT} one can obtain in $\mathcal{O}(1)$ time the penalized cost of a new route $\sigma$ issued from a move and assigned to a vehicle $k$, which is formally expressed in Eq.~\eqref{eqn:FOTW}.
\begin{align}
F(\sigma,k) &=   f_k + r_k \times DIST(\sigma) + \omega \times TW(\sigma) \label{eqn:FOTW}
\end{align}
The term $f_k + r_k \times DIST(\sigma)$ corresponds to the total distance 
of the route $\sigma$ when using the vehicle of type~$k$. The term $\omega \times TW(\sigma)$ computes the penalty due to time windows violations existing in route $\sigma$, where $\omega$ is parameter that controls the level of intensity of such penalty. Note that if all customers are served within their time windows, then the total $TW(\sigma)$ violation equals zero.
It is finally worth mentioning that the partial loads on sequences are used during the search, as in \cite{Pennaetal2013}, to filter subsets of moves that are known in advance to be infeasible with respect to the vehicle capacity.

\subsection{Compound Neighborhood Structures }\label{sec:CNS}


In most methods from the literature, the number of feasible moves is limited by the fleet composition, i.e., the current vehicle type associated with each route. In an attempt to perform more systematic fleet optimization along with route improvements, we designed a generalized version of our local search, that we will reference as \emph{combined neighborhood search} (CNS). 
In the proposed approach, the same inter-route neighborhoods (\textsc{Swap}, \textsc{Shift}, \textsc{2-Opt*} etc...) are generalized with a combined optimization of route-to-vehicle assignments during the move evaluations. This means that a generalized \textsc{Swap} move, for example, is evaluated by jointly swapping the visits \textbf{and} determining the best assignment of vehicle types to each route in the newly created solution. Besides this change, the overall local search scheme and exploration strategy remains the same as described in Section~\ref{LocalSearch}.

To optimize fleet-assignment decisions we tested two approaches.
The first uses an exact method based on the Primal-Dual Algorithm to solve the Assignment Problem (AP) and finds the optimal fleet composition according to the neighborhood tested. Although very slow, this procedure should be viewed as a benchmark, to evaluate what can be gained by means of optimal assignment decisions with each move.
The second approach uses a simple reassignment heuristic, considering only the vehicles that are still unemployed.

%

\paragraph{\textbf{Primal-Dual Algorithm}}

For a given solution, the optimal fleet composition and assignment to routes can be found by solving an Assignment Problem (AP) expressed in Equations (\ref{ObjFunc}-\ref{Domain}).
Let $\mathcal{R}$ be the set of routes and $\mathcal{P}$ be the set of available vehicles.
The model relies on binary decision variables $x_{ij}$, which take value $1$ if and only if route~$i$ is associated to vehicle $j$. 
For each route $i$, let $q_i$ be the load and $d_i$ the distance associated to the route.
For each vehicle $k$, let $Q_k$ be the capacity, $F_k$ the fixed cost and $U_k$ the cost per distance unit.
The cost of an assignment of a vehicle $j \in \mathcal{P}$ to a route $i \in \mathcal{R}$ is given by $c_{ij}$, where $c_{ij} = F_j+U_j\times d_i$ if $q_i \leq Q_j$, otherwise $c_{ij} = \infty$.
\begin{align}
\text{Min} \sum_{i \in \mathcal{R}} \sum_{j \in \mathcal{P}} c_{ij}x_{ij} \label{ObjFunc}
\end{align}
\qquad \qquad \qquad \qquad subject to
\begin{align}
\qquad \qquad \qquad \qquad \sum_{j \in \mathcal{P}} x_{ij} = 1 & & \forall i \in \mathcal{R} \label{R1} \\
\qquad \qquad \qquad \qquad \sum_{i \in \mathcal{R}} x_{ij} = 1 & & \forall j \in \mathcal{P} \label{R2}\\
\qquad \qquad \qquad \qquad x_{ij} \in \{0,1\}  & & \forall i \in \mathcal{R}, \forall j \in \mathcal{P}. \label{Domain}
\end{align}

The objective function (\ref{ObjFunc}) minimizes the sum of the costs by choosing the best assignment of routes to vehicles. Constraints (\ref{R1}) state that a single route from the set $\mathcal{R}$ is associated to only one vehicle \mbox{$j \in \mathcal{P}$}. Constraints (\ref{R2}) requires that a single vehicle from the set $\mathcal{P}$ is assigned to only one route $i \in \mathcal{R}$. Constraints~(\ref{Domain}) define the domain of the decision variables.
Note that AP requires $|\mathcal{R}| = |\mathcal{P}|$, if $|\mathcal{R}| < |\mathcal{P}|$ some dummy routes are created and assigned to vehicles with null costs.
This AP is solved in $\mathcal{O}(|\mathcal{R}|^3)$ operations using the primal-dual algorithm (PDA) of  \cite{McGinnis1983}.

\paragraph{\textbf{Simple Fleet Reassignment}}

The previously-described PDA is exact but computationally expensive.
Thus, we developed an alternative heuristic, called Simple Fleet Reassignment (\textsc{SFR}), which takes into account only the vehicles that are still available (i.e., not assigned to a route) and the routes involved in the move. First, a list with all available vehicles $LV$ is created. Next, for routes $r_1$ and $r_2$ associated with the move, sequentially, the method finds in $LV$ the vehicle type that best fits,  i.e., the vehicle that meets the demands of the customers of the route with the least associated cost. If a better vehicle is found for these routes, the vehicle type assigned to routes $r_1$ and $r_2$ is updated and the solution cost is returned.

\section{Computational Experiments}
\label{sec_Results}

HILS-RVRP was coded in C++ (g++ 4.8.2) and executed in an Intel Core i7 Processor 2.93 GHz with 8~GB of RAM running Ubuntu Linux 14.04 (Kernel 3.10 -- 64 bits). A single thread was used for all tests. The SP models were solved using CPLEX 12.5.1. The proposed algorithm was tested on well-known benchmark instances, available in the literature for each variant considered.

The values of most of the parameters are the same as in \cite{Subramanian2012a}, that is, 
\mbox{$I_{MS}=30$}, $T_{max} = 30 \text{ seconds}$ and $RGap_{max}=0.02$. The maximum number of iterations per ILS was set to $I_{ILS}=n + 5 \times v$ ($n$ and $v$ being the number of customers and vehicles, respectively) as in \cite{Pennaetal2013}.
Infeasibility can occur when an extra vehicle is used while building an initial  solution for the fixed fleet variants. In this case, the extra vehicle has the following penalty values: $f_{|M|+1}=10 \times f_{|M|}$, $r_{|M|+1}=100 \times r_{|M|}$ and $Q_{|M|}=\sum_{i=1}^{n} q_{i}$ (w.l.o.g., the set of vehicle types $M$ is in increasing order of costs and the vehicle type $|M|$ is the one with the largest cost). When the fleet is non-correlated (see Subsection \ref{sec:sacn}), we assume that $|M|$ is the vehicle type associated with the largest fixed cost (or with the largest variable cost in case there are no fixed costs in the problem). 

Since \cite{Subramanian2012a} and \cite{Pennaetal2013} did not consider HFVRPs with time windows or infeasible solutions, the value of the time-window penalty $\omega$, in Eq.~\eqref{eqn:FOTW}, could not be inherited from these works. Therefore, we conducted a series of experiments with different values of $\omega$ on 168 challenging instances of the FSMTW (138), SDepVRPTW (10) and HFFVRPMBTW (20) to calibrate the value of $\omega$. We ran our algorithm ten times for each instance and for $\omega \in \{1,10,100,1000\}$. The idea behind this experiment is to determine a value of $\omega$ for which all runs yield feasible solutions. Table \ref{tbl:omega} reports the results of the experiments. The value $\omega$ is the only one to lead to feasible solutions on all test instances, this value has thus been used in our computational experiments.

\begin{table}[htbp]
	\def\arraystretch{1.2}
	\caption{Percentage of feasible runs for each value of $\omega$ on time windows variants}
	\begin{tabular}{cccccccc}
		\hline\noalign{\smallskip}
		$\omega$& & FSMTW (duration) & FSMTW (distance) & HFFVRPMBTW & SDepVRPTW & & Average \\ \cline{1-1} \cline{3-6} \cline{8-8}
		1 & & 21.30 &	17.97 & 66.50	& 40.00	& & 36.44 \\
		10 & & 59.86 &	39.28 &	91.00	& 74.00	& & 66.04 \\
		100 & & 77.25 &	61.30 &	97.50  & 93.00 & & 82.26 \\
		1000 & & 100.00 & 100.00 & 100.00 & 100.00 & & 100.00\\
		\noalign{\smallskip}\hline
	\end{tabular}
	\label{tbl:omega}
\end{table}

The next set of computational experiments aims to compare the proposed algorithm with the current state-of-the-art methods on a wide range of VRP variants with heterogeneous fleet. The benchmark instances used in those tests are described in \mbox{Table~\ref{tbl:instances}}. In this table, \textbf{Authors} and \textbf{Acr.} indicate the name and the acronym of the authors who proposed the instances, respectively. \textbf{\#Inst.} represents the number of instances of the dataset. \textbf{Variant} denotes the name of the VRP variant, while \textbf{\textit{n}} is the number of customers and \textbf{\textit{m}} the number of vehicles types. \textbf{Costs} corresponds to the types of costs considered, that is, fixed, variable, or both fixed and variable. \textbf{Fleet} represents the fleet size, i.e., limited (L) or unlimited (U), while \textbf{Depots} is the number of depots and \textbf{RD} indicates if the instances of the referred benchmark impose route duration constraints (Y=yes and N=no). Letters F, V and FV in problem names indicate whether fixed costs, variable costs or both are tackled.

The algorithm was executed 10 times on each instance with different random seeds, and a summary of the results is presented in Table \ref{tbl:ResultadosFinal1-HILS-RVRP}. Detailed results, for each instance, are provided in Appendix. A comparison was performed with the best known algorithms reported in the literature.

\begin{table}[htbp]
\def\arraystretch{1.2}
\caption{List of benchmark instances}
\setlength{\tabcolsep}{1.1mm}
\scalebox{0.92} {
\begin{tabular}{llclcccccc}
\hline\noalign{\smallskip}
\multicolumn{1}{c}{\textbf{Authors}} & \multicolumn{1}{l}{\textbf{Acr.}} & \textbf{\#Inst.} & \multicolumn{1}{c}{\textbf{Variant}} & \textbf{\textit{n}} & \textbf{\textit{m}} & \textbf{Costs} & \textbf{Fleet} & \textbf{Depots} & \textbf{RD} \\ 
\noalign{\smallskip}\hline\noalign{\smallskip}
\cite{Golden1984}  & G84 & 12 & FSMVRP-F & [20-100] & [3-6] & F & U & 1 & N \\ 
 &  &  & HFFVRPSD &  &  &  &  &  &  \\ 
\noalign{\smallskip}\hline\noalign{\smallskip}
\cite{Taillard1999}  & T99 & 12 & FSMVRP-V & [20-100] & [3-6] & V/FV & L/U & 1 & N \\ 
 &  &  & FSMVRP-FV &  &  &  &  &  &  \\ 
 &  &  & HFFVRP-V &  &  &  &  &  &  \\ 
 &  &  & HFFVRP-FV &  &  &  &  &  &  \\ 
 &  &  & HFFOVRP &  &  &  &  &  &  \\ 
\noalign{\smallskip}\hline\noalign{\smallskip}
\cite{Brandao2011}  & B11 & 5 & FSMVRP-V & [100-199] & [4-6] & V & L/U & 1 & N \\ 
 &  &  & HFFVRP-V &  &  &  &  &  &  \\ 
\noalign{\smallskip}\hline\noalign{\smallskip}
\cite{Li2007} & LGW07 & 5 & HFFVRP-V & [200-360] & [5-6] & V & L & 1 & N \\ \noalign{\smallskip}\hline\noalign{\smallskip}
\cite{Duhameletal2011} & DLP11 & 96 & HFFVRP-FV & [20-256] & [2-8] & FV & L & 1 & N \\ 
\noalign{\smallskip}\hline\noalign{\smallskip}
\cite{Salhi1997} & SS97 & 23 & MDFSMVRP & [50-360] & 5 & FV & U & [2-9] & N/Y \\ 
\noalign{\smallskip}\hline\noalign{\smallskip}
\cite{Tutuncu2010} & T10 & 18 & HFFVRPB & [50-100] & [3-5] & V & L & 1 & N \\ \noalign{\smallskip}\hline\noalign{\smallskip}
\cite{Salhietal2013} & SNM13 & 36 & FSMVRPB & [20-100] & [3-6] & F & U & 1 & N \\ 
\noalign{\smallskip}\hline\noalign{\smallskip}
\cite{CordeauLaporte2001} & CL01 & 35 & SDepVRP & [27-1008] & [2-6] & - & L & 1 & N/Y \\ 
 &  & 20 & SDepVRPTW & [48-1008] &  &  &  &  &  \\ \hline
\cite{LiuShen1999} & {LS99a} & 168 & FSMVRPTW & 100 & [3-6] & F & U & 1 & N \\
 &  &  & {Min. duration} &  &  &  &  &  &  \\ 
 \noalign{\smallskip}\hline\noalign{\smallskip}
\cite{LiuShen1999} & {LS99b} & 168 & FSMVRPTW & 100 & [3-6] & F & U & 1 & N \\
 &  &  & {Min. distance} &  &  &  &  &  &  \\ 
  \noalign{\smallskip}\hline\noalign{\smallskip}
\cite{Belmecherietal2013} & BPYA13 & 56 & HFFVRPMBTW & 100 & 5 & V & L & 1 & N \\ 
\noalign{\smallskip}\hline
\end{tabular}
}
\label{tbl:instances}
\end{table}

\subsection{Comparison with state-of-the-art methods for HFRVRPs}

Table \ref{tbl:ResultadosFinal1-HILS-RVRP} describes the results obtained by HILS-RVRP. In this table, \textbf{Variant} denotes the HFVRP variant name, \textbf{Bench.} denotes the benchmark set name, \textbf{\textit{n}} is the number of customers, \textbf{Authors} represents the authors of state-of-the-art methods reported in the literature, \textbf{Best Gap} indicates the average gap between the best solution found by HILS-RVRP and the best known solution, \textbf{Avg. Gap} corresponds to the gap between the average solution found by HILS-RVRP and the best known solution,  \textbf{\#Sol. found}  represents the number of best known solutions found, {\textbf{Avg. T} indicates the average CPU time in seconds for} each \textbf{CPU} model, scaled for our 2.93 GHz PC using the performance factors listed in \cite{Dongarra2010}. The best algorithms, according to the \textbf{Best Gap}, are highlighted in boldface.

The algorithm was tested on 13 HFRVRP variants using 12 well-known benchmark datasets with up to 1008 customers. The table is divided into 3 groups of results. The first group describes the results on a classical version of the HFFVRP. The second group presents the results for HFRVRP variants that involves heterogeneous vehicles and another attribute, while the third group displays the results for more complex HFRVRP variants, adding time windows attributes to previously described HFRVRP versions.

 \begin{table}[!hbp]
\def\arraystretch{1.15}
\caption{Summary of the results}
\setlength{\tabcolsep}{.9mm}
\scalebox{0.81} {
\begin{tabular}{lccllrrrlcl}
\hline\noalign{\smallskip}
\multicolumn{1}{c}{} & \multicolumn{1}{c}{} & \multicolumn{1}{c}{} & & \multicolumn{7}{l}{State-of-the-art methods}  \\
\cline{5-11}\noalign{\smallskip}
\multicolumn{1}{c}{} & \multicolumn{1}{c}{} & \multicolumn{1}{c}{} & \multicolumn{1}{c}{} & \multicolumn{1}{c}{} & \multicolumn{1}{c}{\textbf{Best}} & \multicolumn{1}{c}{\textbf{Avg.}} & \multicolumn{1}{c}{{\textbf{Avg. T}}} & \multicolumn{1}{c}{} & \multicolumn{1}{c}{\textbf{\textbf{\#Sol.}}} & \multicolumn{1}{c}{\textbf{CPU}} \\
\multicolumn{1}{c}{\textbf{Variant}} & \multicolumn{1}{c}{\textbf{Bench.}} & \multicolumn{1}{c}{\textit{n}} & \multicolumn{1}{c}{} & \multicolumn{1}{c}{\textbf{Authors}} & \multicolumn{1}{c}{\textbf{Gap}} & \multicolumn{1}{c}{\textbf{Gap}} & \multicolumn{1}{c}{\textbf{(s)}} & \multicolumn{1}{c}{} & \multicolumn{1}{c}{\textbf{found}} & \multicolumn{1}{c}{\textbf{(GHz)}} \\
\noalign{\smallskip}\hline\noalign{\smallskip}
 &  &  &  & \cite{Brandao2011} & 0.35 & \multicolumn{1}{c}{--} & \multicolumn{1}{c}{--} &  & 5/5 & P4 2.6 \\
 HFFVRP-V & B11 & [100 - 199] &  & \cite{Subramanian2012a} & 0.00 & 0.39 & 53.59 &  & 5/5 & I7 2.93 \\
 &  &  &  & \bf{HILS-RVRP} & \textbf{-0.05} & \textbf{0.20} & \textbf{42.25} &  & \textbf{5/5} & I7 2.93 \\
\noalign{\smallskip}\hline\noalign{\smallskip}
 &  &  &  & \cite{Brandao2011} & 0.09 & \multicolumn{1}{c}{--} & 1246.28 &  & 2/5 & P4 2.6 \\
HFFVRP-V  & LGW07  & [200 - 360] &  & \cite{Subramanian2012a} & 0.92 & 2.15 & 302.77 &  & 1/5 & I7 2.93 \\
 &  &  &  & \textbf{HILS-RVRP} & \textbf{-0.02} & \textbf{1.42} & \textbf{549.41} &  & \textbf{4/5} & I7 2.93 \\
\noalign{\smallskip}\hline\noalign{\smallskip}
  &  &  &  & \cite{Duhameletal2011} & 0.86 & \multicolumn{1}{c}{--} & 468.57 &  & 7/96 & P 2.2 \\
 HFFVRP-FV & DLP11 & [20 - 256] &  & \bf{HILS-RVRP} & \textbf{-0.14} & \textbf{0.26} & \textbf{203.84} &  & \textbf{85/96} & I7 2.93 \\
\noalign{\smallskip}\hline\hline\noalign{\smallskip}
  &  &  &  & \small{\cite{Yousefikhoshbakht2014}} & 0.00 & \multicolumn{1}{c}{--} & 116.72 &  & 1/8 & P4 3.0 \\ 
 HFFOVRP & T99 & [50 - 100] &  & \bf{HILS-RVRP} & \textbf{-6.24} & \textbf{-6.15} & \textbf{6.31} &  & \textbf{7/8} & I7 2.93 \\
\noalign{\smallskip}\hline\noalign{\smallskip}
  &  &  &  & {\cite{Salhietal2014}} & {1.52} & \multicolumn{1}{c}{--} & {247.30} & & {12/23} & P4-M \\
 MDFSMVRP & SS97 & [50 - 360] &  & {{\textbf{\cite{Vidaletal-MD-2014}}}} & {\textbf{0.05}} & {\textbf{0.07}} & {\textbf{715.07}} & & {\textbf{21/23}} & Xe 2.93 \\
  &  &  &  & {HILS-RVRP} & {0.06} & {0.76} & {91.43} &  & {17/23} & I7 2.93 \\
\noalign{\smallskip}\hline\noalign{\smallskip}
  &  &  &  & \cite{Tutuncu2010} & 0.00 & \multicolumn{1}{c}{--} & \multicolumn{1}{c}{--} &  & 16/16 & P4 2.66 \\
 HFFVRPB & T10 & [50 - 100]  &  & \bf{HILS-RVRP} & \textbf{-10.41} & \textbf{-9.76} & \textbf{3.26} &  & \textbf{11/11} & I7 2.93 \\
\noalign{\smallskip}\hline\noalign{\smallskip}
  &  &  &  & \cite{Salhietal2013} & 1.24 & \multicolumn{1}{c}{--} & 668.66 &  & 21/36 & P4 3.0 \\
 FSMVRPB & SNM13 & [50 - 100]  &  & \bf{HILS-RVRP} & \textbf{-1.42} & \textbf{-0.77} & \textbf{5.06} &  & \textbf{33/36} & I7 2.93 \\
\noalign{\smallskip}\hline\noalign{\smallskip}
  &  &  &  & \cite{PisingerRopke2007} & 0.24 & {0.84} & 186.82 &  & 19/35 & P4 3 \\
SDepVRP & CL01 & [27 - 1008] &  & \small{\textbf{\cite{CordeauMaischberger2012}}} & \textbf{0.04} & \multicolumn{1}{c}{\textbf{--}} & \multicolumn{1}{c}{\textbf{--}} &  & \textbf{33/35} & X7 2.93 \\ 
  &  &  &  & HILS-RVRP & 0.31 & 1.44 & 451.51 &  & 22/35 & I7 2.93 \\
\noalign{\smallskip}\hline\noalign{\smallskip}
  &  &  &  & \cite{Ozfirat2010} & 0.00 & \multicolumn{1}{c}{--} & 259.42 &  & 12/12 & P4 3 \\
 HFFVRPSD & G84 &  [20 - 100] &  & \textbf{HILS-RVRP}$^*$ & \textbf{-1.80} & \textbf{-1.53}  & \textbf{278.44} &  & \textbf{12/12} & I7 2.93 \\
\noalign{\smallskip}\hline\hline\noalign{\smallskip}
 
 {FSMVRPTW} & {LS99a} & [100] &  & {\bf{\cite{Vidaletal2014}}} & \textbf{0.24} & {\textbf{0.32}} & \textbf{304.90} &  & \textbf{120/168} & Opt 2.2 \\
  & {Min. duration} &  &  & \cite{Kocetal2015} & 0.36 & \multicolumn{1}{c}{--} & \multicolumn{1}{c}{326.51} &  & 29/168 & Xe 2.6 \\ 
  &  &  &  & {HILS-RVRP} & {0.26} & {0.70} & {122.02} &  & {91/168} & I7 2.93 \\
\noalign{\smallskip}\hline\noalign{\smallskip}
 
   {FSMVRPTW} & {LS99b} & [100]  &  & {\cite{Vidaletal2014}} & {0.10} & {0.22} & {282.88} &  & {124/168} & Opt 2.2 \\
  & {Min. distance} &  &  & \bf{\cite{Kocetal2015}} & \textbf{0.07} & \multicolumn{1}{c}{--} & \multicolumn{1}{c}{\textbf{266.42}} &  & \textbf{118/168} & Xe 2.6 \\ 
  &  &  &  & {HILS-RVRP} & {0.10} & {0.39} & {142.19} &  & {110/168} & I7 2.93 \\
\noalign{\smallskip}\hline\noalign{\smallskip}
 
  &  &  &  & \cite{Belmecherietal2013} & 6.96 & \multicolumn{1}{c}{--} & \multicolumn{1}{c}{--} &  & 15/56 & \multicolumn{1}{c}{--} \\
   HFFVRPMBTW & BPYA13 & [100] &  & \cite{BerghidaBoukra2015} & \textbf{2.86} & -- & -- &  & \textbf{23/56} & i7 2.20 \\
  &  &  &  & \bf{HILS-RVRP} & \textbf{-20.02} & \textbf{-19.37} & \textbf{103.46} &  & \textbf{55/56} & I7 2.93 \\
\noalign{\smallskip}\hline\noalign{\smallskip}
 &  &  &  & \small{\cite{CordeauMaischberger2012}} & 0.56 & \multicolumn{1}{c}{--} & \multicolumn{1}{c}{--} &  & 6/20 & X7 2.93 \\
  SDepVRPTW & CL01 & [48 - 288] &  & {\bf{\cite{Vidaletal2013}}} & \textbf{0.10} & {\textbf{0.36}} & \multicolumn{1}{c}{\textbf{328.95}} &  & \textbf{11/20} & Xe 2.93 \\
  &  &  &  & {HILS-RVRP} & {0.34} & {0.97} & {1424.57} &  & \textbf{8/20} & I7 2.93 \\
\noalign{\smallskip}\hline\noalign{\smallskip}
\multicolumn{11}{l}{$^*$: HILS-RVRP without the SP module.}
\end{tabular}
} 
\label{tbl:ResultadosFinal1-HILS-RVRP}
\end{table}

For classical HFVRP variants (group 1), a total of 106 benchmark instances were considered, and they were divided into three different datasets (LGW07, B11, DLP11) containing 5, 5 and 96 instances respectively. HILS-RVRP was capable of obtaining 68 new solutions and equaling 26 of them, meaning that the proposed algorithm found or improved the BKS in 88.7\% of the instances. The gap between the best solution found by HILS-RVRP and the BKS varied from -0.14\% to -0.02\%.
Note that LGW07 contains larger instances than B11, involving up to 360 customers, and  HILS-RVRP was capable of obtaining three new best solutions for this dataset. The DLP11 instances, that are based on real distances from French cities, the HILS-RVND found 61 new solutions.


The second group includes six variants, each with a different set of instances, leading to a total 125 test-problems involving up to 1008 customers. From Table  \ref{tbl:ResultadosFinal1-HILS-RVRP}, we observe that the proposed algorithm always found an average gap lower or equal than 0.76\%, except for the SDepVRP variant, for which HILS-RVRP obtained a value of 1.44\%. One possible reason for not achieving improved solutions for this latter problem is that we did not implement any particular neighborhood structure or perturbation mechanism for this case. 
Furthermore, several improved solutions were found for problems HHFFOVRP, HFFVRPB and FSMVRPB. Regarding the CPU times, HILS-RVRP appears to be in many cases faster than other algorithms from the literature, such as those of \citet{Salhietal2013}  and \citet{Yousefikhoshbakht2014}.

Finally, the third group considers four variants with time-window constraints. HILS-RVRP has been tested on three benchmark datasets, also considering two cases for the FSMVRPTW: either minimizing the sum of the durations, or the total travel distance. Overall, we considered a total of 392 instances involving up to 288 customers. The proposed heuristic found average gaps always smaller than 1\% and many improved solutions were obtained, even for problem FSMVRPTW, which is the most well-studied variant after the classical ones. Variant HFFVRPMBTW is the one that combines the most attributes, and an average improvement of 20\% was achieved, thus suggesting that HILS-RVRP is robust enough to cope with problems that consider several attributes at once. Moreover, the proposed algorithm appears to be competitive in terms of CPU time, except for problem SDepVRPTW, where HILS-RVRP was slower than the state-of-the-art method of \citet{Vidaletal2014}.

\subsection{Sensitivity Analysis -- Combined Neighborhoods} \label{sec:sacn}

\begin{figure}[h]
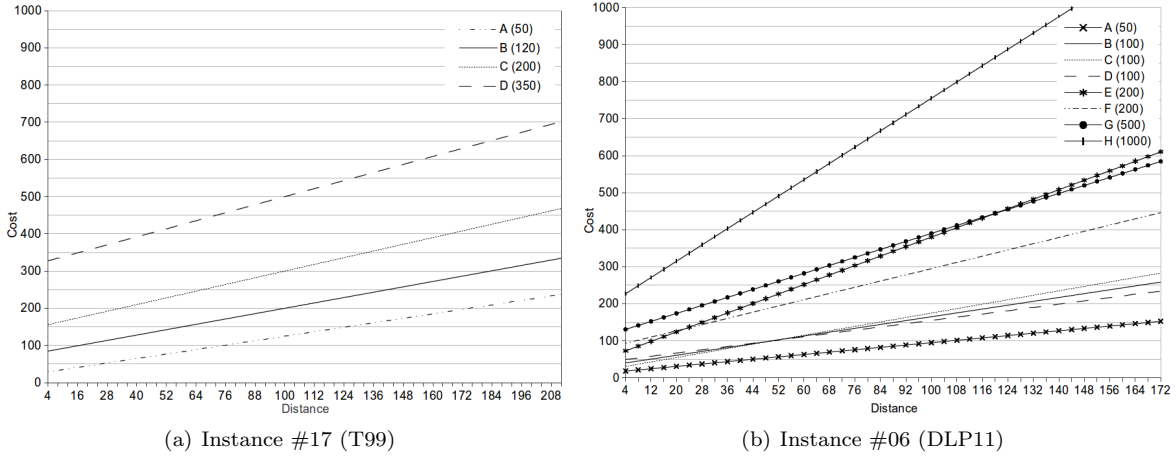

\label{fig:Instances}
\centering
\subfigure[imgTaillard17][Instance \#17 (T99)]{\label{imgTaillard17} \includegraphics[width=7.5cm]{Taillard_17}}
\subfigure[imgDLP06][Instance \#06 (DLP11)]{\label{imgDLP06} \includegraphics[width=7.8cm]{DLP_06}}
\caption{Correlated and uncorrelated Instances}
\label{figInstances}
\end{figure}

The impact of the CNS was investigated on two instance sets of the classical HFVRP. The first set, including T99 and B11 benchmark problems, involves correlated vehicle costs and capacities, i.e., if vehicle types are considered in ascending order of capacities, the fixed costs and variable costs also increase (Fig. \ref{imgTaillard17}). This situation is the rule when one considers vehicle with the same age and technology. The second set, named DLP11, is based on road distances between major cities in different districts of France. The fleet composition is uncorrelated in most problems, see Fig. \ref{imgDLP06}, for instance the larger fixed cost of a hybrid vehicle is compensated by a smaller operating cost. Figure \ref{figInstances} shows the route cost per distance for each vehicle type for one instance of each set. In Figure \ref{imgDLP06}, a route with a customer demand of $100$ can be associated to vehicle type B, C or D. If the distance is smaller than $20$ it is better use vehicle C. Otherwise, if the route distance is greater than $50$, then vehicle D leads to smaller costs. This threshold effect does not happen when fleet capacity and costs are correlated.

We implemented three versions of the algorithm to study the impact of the CNS:
\begin{itemize}[nosep]
\item \textbf{HILS-RVRP}: The HILS-RVRP without the CNS;
\item \textbf{HILS-RVRP-SFR}: The HILS-RVRP with CNS based on the Simple Fleet Reassignment procedure;
\item \textbf{HILS-RVRP-PDA}: The HILS-RVRP with CNS based on the Primal-Dual procedure.
\end{itemize}
These three method variants have been tested with and without the \textsc{Merge} ($P^{(4)}$) neighborhood, leading to overall six versions of the proposed algorithm.
Each version was executed 10 times for each instance and the number of multi-start iterations of the HILS-RVRP parameter was set to $100$.
The results are presented in Tables \ref{Results_Correlated} and \ref{Results_NonCorrelated}. In these two tables, \textbf{Gap} denotes the gap between the average solution, on 10 runs, found by each version of the algorithm and the best known solution of the literature. \textbf{Time} corresponds to the average time, in seconds, of these runs. The best average gap for each version is highlighted in boldface.

Table \ref{Results_Correlated} displays the results on the 8  HFFVRP-D and HFFVRP-FD benchmark instances T99, and the HFFVRP-D instances B11. 
All versions achieved a very similar performance in terms of average gap, but the version HILS-RVRP-PDA without the \textsc{Merge} perturbation slightly outperformed the others. In this set of instances, the use of the \textsc{Merge} procedure also generated slightly worse solutions. This can be related to the structure of the instances, since this perturbation drives the search towards fewer routes, associated to vehicles with larger capacities, and consequently higher costs. Still, these differences of performance remain very limited.
In terms of the average CPU time, using the \textsc{Merge} leads to a time reduction of 3.42\%, while the CNS displays larger CPU times due to the resolution of vehicle-to-route assignment problems.

\begin{table}[htbp]
\caption{Comparative results for all versions of the algorithm on correlated instances }
\def\arraystretch{1.2}
\setlength{\tabcolsep}{1.05mm}
\scalebox{0.83}
{
\begin{tabular}{lclrrrrrrrrlrrrrrrrr}
\noalign{\smallskip}\hline\noalign{\smallskip}
 &  &  & \multicolumn{8}{c}{$P^{(1)}+P^{(2)}+P^{(3)}$} &  & \multicolumn{8}{c}{$P^{(1)}+P^{(2)}+P^{(3)}+P^{(4)}$}  \\
 \cline{4-11} \cline{13-20}
Problem type & \textbf{\textit{n}} &  & \multicolumn{2}{c}{MS-ILS} & & \multicolumn{2}{c}{MS-ILS+SFR} & & \multicolumn{2}{c}{MS-ILS+PDA} & \multicolumn{1}{c}{} & \multicolumn{2}{c}{MS-ILS} & & \multicolumn{2}{c}{MS-ILS+SFR} & & \multicolumn{2}{c}{MS-ILS+PDA} \\
\cline{4-5} \cline{7-8} \cline{10-11} \cline{13-14} \cline{16-17} \cline{19-20}\noalign{\smallskip}
&  &  & \multicolumn{1}{c}{Gap} & \multicolumn{1}{c}{Time} & &  \multicolumn{1}{c}{Gap} & \multicolumn{1}{c}{Time} & & \multicolumn{1}{c}{Gap} & \multicolumn{1}{c}{Time} &  & \multicolumn{1}{c}{Gap} & \multicolumn{1}{c}{Time} & & \multicolumn{1}{c}{Gap} & \multicolumn{1}{c}{Time} & & \multicolumn{1}{c}{Gap} & \multicolumn{1}{c}{Time} \\
&  &  & \multicolumn{1}{c}{(\%)} & \multicolumn{1}{c}{(s)} & & \multicolumn{1}{c}{(\%)} & \multicolumn{1}{c}{(s)} & & \multicolumn{1}{c}{(\%)} & \multicolumn{1}{c}{(s)} &  & \multicolumn{1}{c}{(\%)} & \multicolumn{1}{c}{(s)} & & \multicolumn{1}{c}{(\%)} & \multicolumn{1}{c}{(s)} & & \multicolumn{1}{c}{(\%)} & \multicolumn{1}{c}{(s)} \\
\noalign{\smallskip}\hline\noalign{\smallskip}
HFFVRP-FD$^1$ & 50 -- 100 & & 0.37 & 30.54 & & 0.35 & 32.02 & & \textbf{0.34} & 1639.91 &  & \textbf{0.34} & 27.36 & & 0.36 & 39.51 & & 0.39 & 1448.92 \\
HFFVRP-D$^1$  & 50 -- 100 & & 0.30 & 30.19 & &  0.29 & 30.96 & & \textbf{0.28} & 1900.57 & & 0.30 & 26.37 & & 0.29 & 38.06 & & 0.31 & 1621.06 \\
HFFVRP-D$^2$  & 100 -- 199 & & \textbf{0.21}& 271.40 & & \textbf{0.21} & 384.50 & & \textbf{0.21} & 50779.58 & & \textbf{0.21} & 267.04 & & \textbf{0.21} & 389.86 & & 0.32 & 51331.29 \\
\noalign{\smallskip}\hline\noalign{\smallskip}
\multicolumn{2}{l}{Average} &  & 0.29 & 110.71 & & 0.28 & 149.16 & & \textbf{0.27} & 18106.69 &  & 0.28 & 106.92 & & 0.29 & 155.81 & & 0.34 & 18133.76 \\
\noalign{\smallskip}\hline\noalign{\smallskip}
\multicolumn{6}{l}{\scriptsize $^1$: T99; $^2$: B11.}
\end{tabular}
}
\label{Results_Correlated}
\end{table}

\begin{table}[htb]
\def\arraystretch{1.2}
\caption{Comparative results for all versions of the algorithm on uncorrelated Instances \citep{Duhameletal2011}}
\setlength{\tabcolsep}{1.05mm}
\scalebox{0.83}
{
\begin{tabular}{lclrrrrrrrrlrrrrrrrr}
\hline\noalign{\smallskip}
 &  &  & \multicolumn{8}{c}{$P^{(1)}+P^{(2)}+P^{(3)}$} &  & \multicolumn{8}{c}{$P^{(1)}+P^{(2)}+P^{(3)}+P^{(4)}$}  \\
 \cline{4-11} \cline{13-20}
Problem type & \textbf{\textit{n}} &  & \multicolumn{2}{c}{MS-ILS} & & \multicolumn{2}{c}{MS-ILS+SFR} & & \multicolumn{2}{c}{MS-ILS+PDA} & \multicolumn{1}{c}{} & \multicolumn{2}{c}{MS-ILS} & & \multicolumn{2}{c}{MS-ILS+SFR} & & \multicolumn{2}{c}{MS-ILS+PDA} \\
\cline{4-5} \cline{7-8} \cline{10-11} \cline{13-14} \cline{16-17} \cline{19-20} \noalign{\smallskip}
&  &  & \multicolumn{1}{c}{Gap} & \multicolumn{1}{c}{Time} & &  \multicolumn{1}{c}{Gap} & \multicolumn{1}{c}{Time} & & \multicolumn{1}{c}{Gap} & \multicolumn{1}{c}{Time} &  & \multicolumn{1}{c}{Gap} & \multicolumn{1}{c}{Time} & & \multicolumn{1}{c}{Gap} & \multicolumn{1}{c}{Time} & & \multicolumn{1}{c}{Gap} & \multicolumn{1}{c}{Time} \\
&  &  & \multicolumn{1}{c}{(\%)} & \multicolumn{1}{c}{(s)} & & \multicolumn{1}{c}{(\%)} & \multicolumn{1}{c}{(s)} & & \multicolumn{1}{c}{(\%)} & \multicolumn{1}{c}{(s)} &  & \multicolumn{1}{c}{(\%)} & \multicolumn{1}{c}{(s)} & & \multicolumn{1}{c}{(\%)} & \multicolumn{1}{c}{(s)} & & \multicolumn{1}{c}{(\%)} & \multicolumn{1}{c}{(s)} \\
\noalign{\smallskip}\hline\noalign{\smallskip}
HFFVRP-FD & 20 -- 95    &  & 0.18 & 30.52 & & 0.58 & 40.95 & & 0.57 & 23579.44 &  & 0.19 & 29.22 & & \textbf{0.17} & 42.37 & & \multicolumn{1}{r}{0.18} & 23064.81\\
HFFVRP-FD & 102 -- 147 &  & 0.58 & 164.04 & & 0.98 & 255.80 & & \multicolumn{1}{c}{--} & \multicolumn{1}{c}{--} &  & 0.55 & 169.19 & & \textbf{0.53} & 252.49 & & \multicolumn{1}{c}{--} & \multicolumn{1}{c}{--}  \\
HFFVRP-FD & 152 -- 196 &  & 0.65 & 542.13 & & 0.84 & 820.84 & & \multicolumn{1}{c}{--} & \multicolumn{1}{c}{--} &  & 0.67 & 538.01 & & \textbf{0.64} & 772.70 & &\multicolumn{1}{c}{--} & \multicolumn{1}{c}{--} \\
HFFVRP-FD & 203 -- 256 &  & 0.25 & 1352.96 & & 0.95 & 2197.56 & & \multicolumn{1}{c}{--} & \multicolumn{1}{c}{--} &  & 0.21 & 1244.18 & & \textbf{0.19} & 1935.74 & & \multicolumn{1}{c}{--} & \multicolumn{1}{c}{--} \\
\noalign{\smallskip}\hline\noalign{\smallskip}
Average &  &  & 0.42 & 522.41 & & 0.84 & 828.79 & & \multicolumn{1}{c}{--} & \multicolumn{1}{c}{--} & & 0.41 & 495.15& & \textbf{0.38} & 750.83 & & \multicolumn{1}{c}{--} & \multicolumn{1}{c}{--}  \\
\noalign{\smallskip}\hline
\end{tabular}
}
\label{Results_NonCorrelated}
\end{table}

Table \ref{Results_NonCorrelated} describes the results on the more realistic HFFVRP-FD benchmark instances of \cite{Duhameletal2011}. This set contains 96 instances, ranging from 20 to 256 customers, and with 3 to 8 types of vehicles. This set of instances is divided into four subsets, a ``small'' subset containing 15 instances with less than 100 customers, 38 instances with 100 to 150 customers, 31 instances with 150 to 200 customers, and finally 12 instances and with more than 200 customers.
Due the CPU time requirements of the CNS with \textsc{PDA}, it was only possible to test it on the ``small'' instance subset. On this benchmark, HILS-RVRP-SFR with the \textsc{Merge} perturbation slightly outperforms the other versions in terms of average solution gap, but at the expense of a higher CPU time.
Several new best known solutions have also been generated during these tests for DLP11 benchmark instances. These solutions are presented in the Appendix, and a comparison is established with the best-known solutions obtained by different versions of the GRASPxELS of \cite{Duhameletal2011, Duhameletal2013}.

The conclusion of this experiment is quite counterintuitive. In many previous works on the HFVRP, we commonly assumed that the inherent difficulty of HFVRPs comes from more complex vehicle assignment decisions. However, our attempts to optimize more systematically this decision set, in combination with changes of sequences, did not lead to large solution improvements as they drove the search towards shorter routes. In this context, the CPU effort may be better distributed on a more intensive search on classical neighborhoods, rather than on an extensive search for alternative fleet assignment choices. For this reason, we kept the simpler version of the algorithm, without combined neighborhoods, for the complete tests on all HFVRP variants, and we hope that future research can help to point out efficient strategies for joint assignment and sequencing optimization.



\section{Conclusions and Perspectives} \label{sec_conclusoes}

In this work we proposed an ILS+SP based heuristic algorithm capable of addressing a broad class of VRPs with heterogeneous fleet. We considered variants with asymmetric costs, 
simultaneous pickup and delivery operations, backhauls, multiple depots, site dependency, split deliveries, and time windows. Computational experiments were carried out in hundreds of instances involving up to 1008 customers from more than 10 benchmark datasets. We compared the results found by HILS-RVRP with those obtained by state-of-the-art heuristics in 12 different variants that include at least one of the aforementioned attributes. The proposed heuristic was capable of equaling or improving 71.70\% of the best known solutions and to achieve, in almost all cases, an average gap between the average solutions and the best known solutions below 1\%. The results suggest that HILS-RVRP is a robust and flexible heuristic algorithm that, despite its generality, is still quite competitive with problem-specific algorithm both in terms of solution quality and CPU time.

We finally investigated a larger neighborhood which considers joint changes in the sequences and a systematic optimization of route-to-vehicle assignment decisions. Even if this neighborhood opened the way to a larger class of possible solution refinements, we also observed that it tended to drive the search towards shorter routes and missed high-quality solutions built with longer routes. For future works, the interaction between assignment and sequencing decisions still deserve attention, and other neighborhood and guidance techniques, which aim to better integrate the two decision classes while circumventing the previously-mentioned issues, should be more thoroughly investigated.

\section*{Acknowledgments}
The authors would like to thank the anonymous reviewers for their valuable comments. This research was partially supported by the Conselho Nacional de Desenvolvimento Científico e Tecnológico (CNPq), under grants: 446683/2014-0 (first author); 305223/2015-1, 428549/2016-0 (second author);  308498/2015-1 (fourth author); 400722/2013-5 (third and fifth authors).

\appendix
\include{anexo1}

\bibliographystyle{spbasic}      

\bibliography{referencias}           


\end{document}

%% file: anexo1.tex
\section{Appendix -- Detailed Results} \label{sec:appendix}

{In Tables \ref{tbl:Resultados.Brandao.HFFVRPFV}--\ref{tbl:Resultados.HH.SDepVRPTW}, \textbf{Inst.} denotes the instance identifier, $\textbf{n}$ is the number of customers, \textbf{BKS} represents the value of the best-known solution reported in the literature, \textbf{Sol.}  indicates the value of the solution reported for a single run, \textbf{Best Sol.} and \textbf{Avg. Sol.} correspond to the value of the best and average solution, respectively, \textbf{T} denotes the CPU time in seconds for a single run, \textbf{Avg. T} represents the average CPU time in seconds, \textbf{Gap} denotes 
either the gap between the values of the best solution found by a given algorithm and the BKS, or the mean of the gaps between the values of the best solutions and the BKS.
The last line (\textbf{Average}) presents the average gap between the values of the solutions and the BKS values and the average CPU time in seconds (a dash ``--'', indicates that the time is not available).
Finally, the best results are highlighted in boldface and new improved solutions found by HILS-RVRP are italicized.}

\subsection{HFFVRP-V} \label{sec:Result.HH.HFFVRP}

Detailed results obtained for the HFFVRP-V instances of:  (i) \cite[B11]{Brandao2011}, compared with the TSA of \cite[B11]{Brandao2011} and the ILS-RVND-SP of \cite[SPUO12]{Subramanian2012a} (Table \ref{tbl:Resultados.Brandao.HFFVRPFV}); and (ii)  \cite[LGW07]{Li2007}, compared with \cite[B11]{Brandao2011} and \cite[SPUO12]{Subramanian2012a} (Table \ref{tbl:Resultados.Li.HFFVRPFV}).


\begin{table}[htbp]
\caption{Results for the HFFVRP-V \cite[B11]{Brandao2011}.}
\setlength{\tabcolsep}{1.5mm}
\scalebox{0.82} {
\begin{tabular}{lrrrclrrlrrrrr}
\hline\noalign{\smallskip}
 &  &  & \multicolumn{2}{l}{TSA} &  & \multicolumn{2}{l}{ILS-RVND-SP} & & \multicolumn{3}{l}{ HILS-RVRP } & & \\
 &  &  & \multicolumn{ 2}{l}{B11} &  &  SPUO12 & & & & &  \\
\cline{4-5} \cline{7-8} \cline{10-14}\noalign{\smallskip}
Inst.  & \multicolumn{1}{c}{\textit{n}} & \multicolumn{1}{c}{BKS} & \multicolumn{1}{c}{{Sol.}} & \multicolumn{1}{c}{ T (s)} &  & \multicolumn{1}{c}{Best Sol.} & \multicolumn{1}{c}{{Avg. T$^{a}$ (s)}} &  & \multicolumn{1}{c}{Best Sol.} & \multicolumn{1}{c}{Gap (\%)} & \multicolumn{1}{c}{{Avg. Sol.}} & \multicolumn{1}{c}{{Avg. T$^{a}$ (s)}} & \multicolumn{1}{c}{{Avg. Gap (\%)}} \\
\noalign{\smallskip}\hline\noalign{\smallskip}
N1 & 150 & 2235.87 & 2243.76 & -- &   & \textbf{2235.87} & 51.50 & & \textbf{\textit{2234.13}} & -0.43 & 2241.91 & 39.10 & -0.08 \\
N2 & 199 & 2864.83 & 2874.13 & -- &   & \textbf{2864.83} & 102.77 & & \textbf{\textit{2859.82}} & -0.50 & 2881.54 & 102.25 & 0.26 \\
N3 & 120 & 2378.99 & 2386.90 & -- &   & \textbf{2378.99} & 51.71 & & \textbf{\textit{2378.99}} & -0.33 & 2381.97 & 21.89 & -0.21 \\
N4 & 100 & 1839.22 & 1839.22 & -- &   & \textbf{1839.22} & 9.64 & & \textbf{1839.22} & 0.00 & \textbf{1839.22} & 10.24 & 0.00 \\
N5 & 134 & 2047.81 & 2062.48 & -- &   & \textbf{2047.81} & 52.33 & & \textbf{\textit{2047.81}} & -0.71 & \textbf{\textit{2047.81}} & 37.78 & -0.71 \\
\noalign{\smallskip}\hline\noalign{\smallskip}
 \multicolumn{3}{l}{Average} &  0.35 & -- &  & 0.00 & 53.59 & & & -0.05 &  & 42.25 & 0.20 \\
 \hline
\end{tabular}
} 
\label{tbl:Resultados.Brandao.HFFVRPFV}
\end{table}

%
\begin{table}[hbtp]
\caption{Results for the HFFVRP-V \cite[LGW07]{Li2007}.}
\setlength{\tabcolsep}{1.1mm}
\scalebox{0.8} {
\begin{tabular}{lrrrrrrrrrrrrr}
\hline\noalign{\smallskip}
 &  &  & \multicolumn{2}{l}{TSA} &  & \multicolumn{2}{l}{ILS-RVND-SP} &  & \multicolumn{4}{l}{ HILS-RVRP } \\
 &  &  & \multicolumn{ 2}{l}{B11} &  & \multicolumn{ 2}{l}{SPUO12}   &  &  &  &  &  \\
\cline{4-5} \cline{7-8} \cline{10-14}\noalign{\smallskip}
Inst.  & \multicolumn{1}{c}{\textit{n}} & \multicolumn{1}{c}{BKS} & \multicolumn{1}{c}{{Sol.}} & \multicolumn{1}{c}{T (s)} &  & \multicolumn{1}{c}{Best Sol.} & \multicolumn{1}{c}{{Avg. T$^{b}$ (s)}} &  & \multicolumn{1}{c}{Best Sol.} & \multicolumn{1}{c}{Gap (\%)} & \multicolumn{1}{c}{{Avg. Sol.$^b$}} & \multicolumn{1}{c}{{Avg. T$^{b}$ (s)}} & \multicolumn{1}{c}{{Avg. Gap$^b$ (\%)}} \\
\noalign{\smallskip}\hline\noalign{\smallskip}
H1 & 200 & 12050.08\tiny{ \ } & \textbf{12050.08} & 1395 &   & \textbf{12050.08} & 72.10 &   & \textbf{12050.08} & 0.00 & 12061.45 & 91.38 & 0.09 \\
H2 & 240 & 10208.32$^1$ & 10226.17 & 3650 &   & 10329.15 & 176.43 &   & \textbf{\textit{10207.48}} & -0.01 & 10349.17 & 282.37 & 1.38 \\
H3 & 280 & 16223.39$^1$ & 16230.21 & 2822 &   & 16282.41 & 259.61 &   & \textbf{\textit{16193.48}} & -0.18 & 16317.62 & 364.01 & 0.58 \\
H4 & 320 & 17458.65\tiny{ \ } & \textbf{17458.65} & 8734 &  & 17743.68 & 384.52 &   & 17490.00 & 0.18 & 17842.46 & 741.72 & 2.20 \\
H5 & 360 & 23166.56$^1$ & 23220.72 & 13321 &   & 23493.87 & 621.17 &   & \textbf{\textit{23151.55}} & -0.06 & 23822.40 & 1267.57 & 2.83 \\
\noalign{\smallskip}\hline\noalign{\smallskip}
\multicolumn{3}{l}{Average} &  0.09 & 1246.28 &  & 0.92 & 302.77 & & & -0.02 &  & 549.41 & 1.42 \\
\noalign{\smallskip}\hline\noalign{\smallskip}
\multicolumn{14}{l}{{\scriptsize $^1$: Obtained by TSA with a different parameter tuning.}} \\

\end{tabular}
} 
\label{tbl:Resultados.Li.HFFVRPFV}
\end{table}
%

\newpage
\subsection{HFFVRP-FV}

Detailed results obtained for the HFFVRP instances of \cite[DLP11]{Duhameletal2011}, compared with the sequential version of the GRASPxELS of \cite[DLP13]{ Duhameletal2013}. Column BKS in Tables \ref{tbl:Resultados.HH.HFFVRP_dlp_p}-\ref{tbl:Resultados.HH.HFFVRP_dlp_gg} indicates the best results considering all versions of the GRASPxELS (sequential or parallel). The instances were divided into four sets:
\begin{itemize}[nosep]
\item DLP11--set 1: [20 -- 95] customers (Table~\ref{tbl:Resultados.HH.HFFVRP_dlp_p});
\item DLP11--set 2: [102 -- 147] customers  (Table~\ref{tbl:Resultados.HH.HFFVRP_dlp_m});
\item DLP11--set 3: [152 -- 196] customers  (Table~\ref{tbl:Resultados.HH.HFFVRP_dlp_g});
\item DLP11--set 4: [203 -- 256] customers  (Table~\ref{tbl:Resultados.HH.HFFVRP_dlp_gg}).
\end{itemize}
%

\begin{table}[htbp]
\caption{Results for HFFVRP-FV \small{\cite[DLP11--set 1]{Duhameletal2011}}}
\setlength{\tabcolsep}{1.5mm}
\scalebox{0.90} {
\begin{tabular}{lcrlrrlrrrrr}
\hline\noalign{\smallskip}
 & \multicolumn{1}{l}{} & \multicolumn{1}{l}{} &  & \multicolumn{2}{l}{GRASPxELS} & & \multicolumn{5}{l}{HILS-RVRP}  \\
 & \multicolumn{1}{l}{} & \multicolumn{1}{l}{} &  & \multicolumn{2}{l}{DLP13} & & \multicolumn{5}{l}{}  \\
\cline{5-6} \cline{8-12} \noalign{\smallskip}
\multicolumn{1}{l}{Inst.} & \textit{n} & \multicolumn{1}{l}{BKS} &  & \multicolumn{1}{c}{{Sol.}} & \multicolumn{1}{c}{T (s)} & \multicolumn{1}{c}{} & \multicolumn{1}{c}{Best Sol.} & \multicolumn{1}{c}{Gap (\%)} & \multicolumn{1}{c}{{Avg. Sol.}} & \multicolumn{1}{c}{{Avg. T(s)}} & \multicolumn{1}{c}{{Avg. Gap (\%)}} \\
\noalign{\smallskip}\hline\noalign{\smallskip}
HVRP\_01\_DLP & 92 & 9210.14 &  & \textbf{9210.14} & 52.29 &  & \textbf{9210.14} & 0.00 & 9213.71 & 44.49 & 0.04 \\
HVRP\_08\_DLP & 84 & 4591.75 &  & 4598.49 & 304.85 &  & \textbf{4591.75} & 0.00 & 4595.65 & 12.00 & 0.08 \\
HVRP\_10\_DLP & 69 & 2107.55 &  & \textbf{2107.55} & 24.83 &  & \textbf{2107.55} & 0.00 & 2107.83 & 12.06 & 0.01 \\
HVRP\_11\_DLP & 95 & 3367.41 &  & 3370.47 & 264.61 &  & \textbf{3367.41} & 0.00 & 3373.77 & 22.10 & 0.19 \\
HVRP\_36\_DLP & 85 & 5684.61 &  & 5759.34 & 104.39 &  & \textbf{5684.61} & 0.00 & 5700.14 & 29.75 & 0.27 \\
HVRP\_39\_DLP & 77 & 2923.72 &  & 2934.55 & 182.11 &  & \textbf{\textit{2921.40}} & -0.08 & 2932.75 & 19.51 & 0.31 \\
HVRP\_43\_DLP & 86 & 8737.02 &  & 8764.75 & 219.91 &  & \textbf{8737.02} & 0.00 & 8746.38 & 75.33 & 0.11 \\
HVRP\_52\_DLP & 59 & 4027.27 &  & 4029.42 & 39.97 &  & \textbf{4027.27} & 0.00 & 4030.44 & 16.08 & 0.08 \\
HVRP\_55\_DLP & 56 & 10244.34 &  & 10247.86 & 190.76 &  & \textbf{10244.34} & 0.00 & 10250.98 & 13.23 & 0.06 \\
HVRP\_70\_DLP & 78 & 6685.24 &  & 6689.61 & 120.60 &  & \textbf{\textit{6684.56}} & -0.01 & 6694.43 & 15.67 & 0.14 \\
HVRP\_75\_DLP & 20 & 452.85 &  & \textbf{452.85} & 0.02 &  & \textbf{452.85} & 0.00 & \textbf{452.85} & 1.15 & 0.00 \\
HVRP\_82\_DLP & 79 & 4766.74 &  & 4774.26 & 144.51 &  & \textbf{4766.74} & 0.00 & 4771.33 & 36.22 & 0.10 \\
HVRP\_92\_DLP & 35 & 564.39 &  & \textbf{564.39} & 20.63 &  & \textbf{564.39} & 0.00 & 564.65 & 4.52 & 0.05 \\
HVRP\_93\_DLP & 39 & 1036.99 &  & \textbf{1036.99} & 27.39 &  & \textbf{1036.99} & 0.00 & 1038.34 & 6.98 & 0.13 \\
HVRP\_94\_DLP & 46 & 1378.25 &  & 1378.66 & 15.68 &  & \textbf{1378.25} & 0.00 & \textbf{1378.25} & 31.37 & 0.00 \\
\noalign{\smallskip}\hline\noalign{\smallskip}
 \multicolumn{2}{l}{Average} & \multicolumn{1}{l}{} &  & 0.17 & 114.17 &  & \multicolumn{1}{l}{} & -0.01 & \multicolumn{1}{l}{} & 19.17 & 0.10 \\
 \hline
\end{tabular}
} 
\label{tbl:Resultados.HH.HFFVRP_dlp_p}
\end{table}

\begin{table}[htbp]
\caption{Results for HFFVRP-FV \small{\cite[DLP11--set 2]{Duhameletal2011}}}
\setlength{\tabcolsep}{1.2mm}
\scalebox{0.9} {
\begin{tabular}{lcrlrrlrrrrr}
\hline\noalign{\smallskip}
 & \multicolumn{1}{l}{} & \multicolumn{1}{l}{} &  & \multicolumn{2}{l}{GRASPxELS} & & \multicolumn{5}{l}{HILS-RVRP}  \\
  & \multicolumn{1}{l}{} & \multicolumn{1}{l}{} &  & \multicolumn{2}{l}{DLP13} & & & & & &  \\
\cline{5-6} \cline{8-12} \noalign{\smallskip}
\multicolumn{1}{l}{Inst.} & \textit{n} & \multicolumn{1}{l}{BKS} &  & \multicolumn{1}{c}{{Sol.}} & \multicolumn{1}{c}{T (s)} & \multicolumn{1}{c}{} & \multicolumn{1}{c}{Best Sol.} & \multicolumn{1}{c}{Gap (\%)} & \multicolumn{1}{c}{{Avg. Sol.}} & \multicolumn{1}{c}{{Avg. T(s)}} & \multicolumn{1}{c}{{Avg. Gap (\%)}} \\
\noalign{\smallskip}\hline\noalign{\smallskip}
\small{HVRP\_03\_DLP} & 124 & 10738.28 &  & 11320.58 & 512.1 &  & \textbf{\textit{10709.66}} & -0.27 & \textbf{\textit{10727.52}} & 91.32 & -0.10 \\
\small{HVRP\_05\_DLP} & 116 & 10903.63 &  & 10963.62 & 488.63 &  & \textbf{\textit{10876.48}} & -0.25 & \textbf{\textit{10897.93}} & 26.59 & -0.05 \\
\small{HVRP\_06\_DLP} & 121 & 11692.85 &  & 11792.94 & 367.91 &  & \textbf{\textit{11688.64}} & -0.04 & 11734.52 & 43.64 & 0.36 \\
\small{HVRP\_07\_DLP} & 108 & 8095.88 &  & 8130.50 & 306.09 &  & \textbf{\textit{8089.46}} & -0.08 & 8144.80 & 27.47 & 0.60 \\
\small{HVRP\_12\_DLP} & 112 & 3543.99 &  & \textbf{3543.99} & 71.46 &  & \textbf{3543.99} & 0.00 & 3547.92 & 57.60 & 0.11 \\
\small{HVRP\_13\_DLP} & 119 & 6696.43 &  & 6713.14 & 303.37 &  & \textbf{6696.43} & 0.00 & 6703.23 & 50.97 & 0.10 \\
\small{HVRP\_16\_DLP} & 129 & 4156.97 &  & 4161.61 & 180.91 &  & \textbf{4156.97} & 0.00 & 4164.03 & 68.70 & 0.17 \\
\small{HVRP\_17\_DLP} & 105 & 5362.83 &  & 5370.05 & 172.82 &  & \textbf{5362.83} & 0.00 & 5367.76 & 42.09 & 0.09 \\
\small{HVRP\_2A\_DLP} & 113 & 7793.16 &  & 7885.93 & 298.92 &  & \textbf{7793.16} & 0.00 & 7796.54 & 43.38 & 0.04 \\
\small{HVRP\_2B\_DLP} & 107 & 8464.69 &  & 8537.31 & 303.14 &  & \textbf{\textit{8462.56}} & -0.03 & 8499.95 & 54.40 & 0.42 \\
\small{HVRP\_21\_DLP} & 126 & 5141.49 &  & 5154.38 & 330.23 &  & \textbf{\textit{5139.84}} & -0.03 & 5166.11 & 47.20 & 0.48 \\
\small{HVRP\_25\_DLP} & 143 & 7206.64 &  & 7228.54 & 518.28 &  & 7209.29 & 0.04 & 7230.50 & 123.79 & 0.33 \\
\small{HVRP\_26\_DLP} & 126 & 6446.31 &  & 6481.93 & 350.71 &  & \textbf{\textit{6433.21}} & -0.20 & 6461.05 & 149.06 & 0.23 \\
\small{HVRP\_28\_DLP} & 141 & 5531.06 &  & 5542.76 & 343.06 &  & \textbf{\textit{5530.55}} & -0.01 & 5542.80 & 101.82 & 0.21 \\
\small{HVRP\_30\_DLP} & 112 & 6313.39 &  & 6321.69 & 201.39 &  & 6315.70 & 0.04 & 6342.42 & 82.91 & 0.46 \\
\small{HVRP\_31\_DLP} & 131 & 4091.52 &  & 4103.88 & 308.39 &  & \textbf{4091.52} & 0.00 & 4112.64 & 102.25 & 0.52 \\
\small{HVRP\_34\_DLP} & 136 & {5758.09} &  & 5800.12 & 405.62 &  & \textbf{\textit{5747.25}} & -0.19 & 5785.59 & 65.57 & 0.48 \\
\small{HVRP\_40\_DLP} & 132 & 11123.56 &  & 11172.98 & 614.92 &  & \textbf{\textit{11118.57}} & -0.04 & 11171.17 & 90.62 & 0.43 \\
\small{HVRP\_41\_DLP} & 135 & 7616.17 &  & 7679.32 & {325.80} &  & \textbf{\textit{7597.27}} & -0.25 & 7672.27 & 68.18 & 0.74 \\
\small{HVRP\_47\_DLP} & 111 & 16206.14 &  & 16222.94 & 333.85 &  & \textbf{\textit{16156.12}} & -0.31 & 16247.77 & 41.18 & 0.26 \\
\small{HVRP\_48\_DLP} & 111 & 21318.04 &  & 21413.92 & {371.30} &  & \textbf{\textit{21309.94}} & -0.04 & 21391.58 & 45.75 & 0.34 \\
\small{HVRP\_51\_DLP} & 129 & 7721.47 &  & 7780.88 & {315.60} &  & \textbf{7721.47} & 0.00 & 7787.85 & 58.47 & 0.86 \\
\small{HVRP\_53\_DLP} & 115 & 6434.83 &  & 6470.49 & 418.17 &  & \textbf{6434.83} & 0.00 & 6454.77 & 36.09 & 0.31 \\
\small{HVRP\_60\_DLP} & 137 & 17037.39 &  & 17067.85 & 444.32 &  & \textbf{\textit{17036.59}} & 0.00 & 17055.35 & 73.38 & 0.11 \\
\small{HVRP\_61\_DLP} & 111 & 7295.67 &  & 7300.10 & 108.21 &  & \textbf{\textit{7292.03}} & -0.05 & 7302.40 & 37.38 & 0.09 \\
\small{HVRP\_66\_DLP} & 150 & 12830.82 &  & 13319.73 & 442.89 &  & \textbf{\textit{12783.94}} & -0.37 & 12922.52 & 113.74 & 0.71 \\
\small{HVRP\_68\_DLP} & 125 & 8976.53 &  & 9135.23 & 269.63 &  & \textbf{\textit{8970.63}} & -0.07 & 9123.03 & 67.86 & 1.63 \\
\small{HVRP\_73\_DLP} & 137 & 10195.33 &  & 10243.66 & 598.34 &  & \textbf{10195.33} & 0.00 & 10195.36 & 73.57 & 0.00 \\
\small{HVRP\_74\_DLP} & 125 & 11586.87 &  & 11732.54 & 246.66 &  & \textbf{\textit{11586.58}} & 0.00 & 11591.23 & 82.46 & 0.04 \\
\small{HVRP\_79\_DLP} & 147 & 7259.54 &  & 7314.89 & 473.69 &  & \textbf{7259.54} & 0.00 & 7289.26 & 122.29 & 0.41 \\
\small{HVRP\_81\_DLP} & 106 & 10700.47 &  & 10715.28 & 83.71 &  & \textbf{\textit{10686.31}} & -0.13 & \textbf{\textit{10700.27}} & 58.03 & 0.00 \\
\small{HVRP\_83\_DLP} & 124 & 10019.15 &  & 10019.83 & 332.47 &  & 10020.07 & 0.01 & 10048.17 & 72.48 & 0.29 \\
\small{HVRP\_84\_DLP} & 105 & 7227.88 &  & 7269.55 & 206.41 &  & \textbf{7227.88} & 0.00 & 7237.93 & 54.37 & 0.14 \\
\small{HVRP\_85\_DLP} & 146 & 8779.76 &  & 8874.31 & 382.98 &  & \textbf{\textit{8773.08}} & -0.08 & 8818.55 & 91.21 & 0.44 \\
\small{HVRP\_87\_DLP} & 108 & 3753.87 &  & 3753.87 & 104.11 &  & \textbf{3753.87} & 0.00 & 3756.97 & 31.41 & 0.08 \\
\small{HVRP\_88\_DLP} & 127 & 12402.85 &  & 12443.41 & 632.22 &  & \textbf{\textit{12388.23}} & -0.12 & 12405.80 & 46.52 & 0.02 \\
\small{HVRP\_89\_DLP} & 134 & 7106.84 &  & 7135.36 & 245.63 &  & \textbf{\textit{7086.36}} & -0.29 & \textbf{\textit{7102.98}} & 76.51 & -0.05 \\
\small{HVRP\_90\_DLP} & 102 & 2346.13 &  & 2360.83 & 15.36 &  & 2346.43 & 0.01 & 2356.31 & 47.78 & 0.43 \\
\noalign{\smallskip}\hline\noalign{\smallskip}
 \multicolumn{2}{l}{Average} & \multicolumn{1}{l}{} &  & 0.71 & 327.09 &  & \multicolumn{1}{l}{} & -0.07 & \multicolumn{1}{l}{} & 67.59 & 0.31 \\
\hline
\end{tabular}
} 
\label{tbl:Resultados.HH.HFFVRP_dlp_m}
\end{table}

\begin{table}[htbp]
\caption{Results for HFFVRP-FV \small{\cite[DLP11--set 3]{Duhameletal2011}}}
\setlength{\tabcolsep}{1.2mm}
\scalebox{0.9} {
\begin{tabular}{lcrlrrlrrrrr}

\hline\noalign{\smallskip}
 & \multicolumn{1}{l}{} & \multicolumn{1}{l}{} &  & \multicolumn{2}{l}{GRASPxELS} & & \multicolumn{5}{l}{HILS-RVRP}  \\
 & \multicolumn{1}{l}{} & \multicolumn{1}{l}{} &  & \multicolumn{2}{l}{DLP13} & & & & & &  \\
\cline{5-6} \cline{8-12}\noalign{\smallskip}
\multicolumn{1}{l}{Inst.} & \textit{n} & \multicolumn{1}{l}{BKS} &  & \multicolumn{1}{c}{{Sol.}} & \multicolumn{1}{c}{T (s)} & \multicolumn{1}{c}{} & \multicolumn{1}{c}{Best Sol.} & \multicolumn{1}{c}{Gap (\%)} & \multicolumn{1}{c}{{Avg. Sol.}} & \multicolumn{1}{c}{{Avg. T(s)}} & \multicolumn{1}{c}{{Avg. Gap (\%)}} \\
\noalign{\smallskip}\hline\noalign{\smallskip}
\small{HVRP\_02\_DLP} & 181 & 11678.44 &  & \textbf{11678.44} & 689.81 &  & \textbf{\textit{11675.26}} & -0.03 & 11695.78 & 187.92 & 0.15 \\
\small{HVRP\_04\_DLP} & 183 & 10808.31 &  & 11030.42 & 667.11 &  & \textbf{\textit{10748.17}} & -0.56 & \textbf{\textit{10775.93}} & 171.66 & -0.30 \\
\small{HVRP\_09\_DLP} & 167 & 7619.19 &  & 7654.45 & 319.39 &  & \textbf{\textit{7603.38}} & -0.21 & 7630.55 & 232.55 & 0.15 \\
\small{HVRP\_14\_DLP} & 176 & 5644.92 &  & 5676.98 & 361.72 &  & 5657.62 & 0.22 & 5697.17 & 357.30 & 0.93 \\
\small{HVRP\_15\_DLP} & 188 & {8236.40} &  & 8367.71 & 905.21 &  & \textbf{\textit{8220.64}} & -0.19 & 8285.79 & 158.34 & 0.60 \\
\small{HVRP\_24\_DLP} & 163 & 9101.47 &  & 9186.30 & {443.10} &  & 9119.92 & 0.20 & 9189.22 & 163.77 & 0.96 \\
\small{HVRP\_29\_DLP} & 164 & 9143.69 &  & 9176.51 & 122.02 &  & \textbf{\textit{9142.86}} & -0.01 & 9149.12 & 232.18 & 0.06 \\
\small{HVRP\_33\_DLP} & 189 & 9421.01 &  & 9563.18 & 606.39 &  & \textbf{\textit{9410.99}} & -0.11 & 9471.26 & 344.96 & 0.53 \\
\small{HVRP\_35\_DLP} & 168 & 9574.71 &  & 9817.94 & 811.07 &  & \textbf{\textit{9555.92}} & -0.20 & 9585.91 & 144.06 & 0.12 \\
\small{HVRP\_37\_DLP} & 161 & 6858.23 &  & 6963.61 & 571.37 &  & \textbf{\textit{6850.77}} & -0.11 & 6875.28 & 245.33 & 0.25 \\
\small{HVRP\_42\_DLP} & 178 & 10902.84 &  & 11118.66 & 966.84 &  & \textbf{\textit{10817.90}} & -0.78 & 10995.75 & 246.34 & 0.85 \\
\small{HVRP\_44\_DLP} & 172 & 12197.46 &  & 12351.49 & 744.39 &  & \textbf{\textit{12191.48}} & -0.05 & 12314.24 & 159.25 & 0.96 \\
\small{HVRP\_45\_DLP} & 170 & 10484.23 &  & 10546.69 & 415.02 &  & \textbf{\textit{10476.25}} & -0.08 & 10614.48 & 147.07 & 1.24 \\
\small{HVRP\_50\_DLP} & 187 & 12374.04 &  & 12538.63 & 365.46 &  & \textbf{\textit{12370.94}} & -0.03 & 12430.18 & 374.43 & 0.45 \\
\small{HVRP\_54\_DLP} & 172 & 10393.23 &  & 10426.98 & 565.12 &  & \textbf{\textit{10351.97}} & -0.40 & 10435.58 & 203.19 & 0.41 \\
\small{HVRP\_56\_DLP} & 153 & 31090.71 &  & 31292.64 & 339.08 &  & \textbf{\textit{31030.19}} & -0.19 & 31144.98 & 260.09 & 0.17 \\
\small{HVRP\_57\_DLP} & 163 & 44818.18 &  & 45112.39 & 471.94 &  & \textbf{\textit{44781.64}} & -0.08 & 44899.36 & 250.64 & 0.18 \\
\small{HVRP\_59\_DLP} & 193 & 14282.59 &  & 14367.47 & 476.61 &  & 14304.46 & 0.15 & 14357.81 & 312.23 & 0.53 \\
\small{HVRP\_63\_DLP} & 174 & 19951.76 &  & 20513.10 & {253.10} &  & 20022.94 & 0.36 & 20281.49 & 213.10 & 1.89 \\
\small{HVRP\_64\_DLP} & 161 & 17157.37 &  & 17157.37 & 70.38 &  & \textbf{\textit{17135.16}} & -0.13 & 17157.79 & 106.02 & 0.00 \\
\small{HVRP\_67\_DLP} & 172 & 10937.67 &  & 11090.66 & 506.65 &  & \textbf{\textit{10884.91}} & -0.48 & 10945.00 & 275.36 & 0.07 \\
\small{HVRP\_69\_DLP} & 152 & 9162.78 &  & 9241.75 & 205.32 &  & \textbf{\textit{9147.54}} & -0.17 & 9190.46 & 117.84 & 0.30 \\
\small{HVRP\_71\_DLP} & 186 & 9870.22 &  & 9936.35 & 389.13 &  & \textbf{\textit{9834.40}} & -0.36 & 9915.73 & 108.20 & 0.46 \\
\small{HVRP\_72\_DLP} & 186 & 5905.58 &  & 5948.99 & 458.28 &  & \textbf{\textit{5903.81}} & -0.03 & 5949.29 & 225.19 & 0.74 \\
\small{HVRP\_76\_DLP} & 152 & 12018.26 &  & 12086.57 & 426.51 &  & \textbf{\textit{11994.40}} & -0.20 & 12040.78 & 138.37 & 0.19 \\
\small{HVRP\_77\_DLP} & 190 & 6930.44 &  & 7004.97 & 278.69 &  & \textbf{\textit{6916.01}} & -0.21 & 6974.86 & 271.78 & 0.64 \\
\small{HVRP\_78\_DLP} & 190 & 7035.01 &  & 7066.17 & {439.70} &  & 7053.62 & 0.26 & 7122.27 & 524.00 & 1.24 \\
\small{HVRP\_80\_DLP} & 171 & 6816.89 &  & 6864.75 & 410.38 &  & 6819.71 & 0.04 & 6843.88 & 211.74 & 0.40 \\
\small{HVRP\_86\_DLP} & 153 & 9030.68 &  & 9085.66 & 440.02 &  & \textbf{\textit{9027.84}} & -0.03 & 9048.94 & 252.05 & 0.20 \\
\small{HVRP\_91\_DLP} & 196 & 6377.48 &  & 6419.23 & 672.65 &  & \textbf{\textit{6374.27}} & -0.05 & 6403.29 & 423.12 & 0.40 \\
\small{HVRP\_95\_DLP} & 183 & {6181.60} &  & 6237.61 & 206.09 &  & \textbf{\textit{6175.62}} & -0.10 & 6232.75 & {554.50} & 0.83 \\
\noalign{\smallskip}\hline\noalign{\smallskip}
 \multicolumn{2}{l}{Average} & \multicolumn{1}{l}{} & \multicolumn{1}{l}{} & 0.98 & 470.92 & \multicolumn{1}{l}{} & \multicolumn{1}{l}{} & -0.11 & \multicolumn{1}{l}{} & 245.57 & 0.50 \\
 \hline
\end{tabular}
} 
\label{tbl:Resultados.HH.HFFVRP_dlp_g}
\end{table}

\begin{table}[htbp]
\caption{Results for HFFVRP-FV \small{\cite[DLP11--set 4]{Duhameletal2011}}}
\setlength{\tabcolsep}{1.2mm}
\scalebox{0.9} {
\begin{tabular}{lcrlrrlrrrrr}
\hline\noalign{\smallskip}
 & \multicolumn{1}{l}{} & \multicolumn{1}{l}{} &  & \multicolumn{2}{l}{GRASPxELS} & & \multicolumn{5}{l}{HILS-RVRP}  \\
  & \multicolumn{1}{l}{} & \multicolumn{1}{l}{} &  & \multicolumn{2}{l}{DLP13} & & & & & &  \\
\cline{5-6} \cline{8-12}\noalign{\smallskip}
\multicolumn{1}{l}{Inst.} & \textit{n} & \multicolumn{1}{l}{BKS} &  & \multicolumn{1}{c}{{Sol.}} & \multicolumn{1}{c}{T (s)} & \multicolumn{1}{c}{} & \multicolumn{1}{c}{Best Sol.} & \multicolumn{1}{c}{Gap (\%)} & \multicolumn{1}{c}{{Avg. Sol.}} & \multicolumn{1}{c}{{Avg. T(s)}} & \multicolumn{1}{c}{{Avg. Gap (\%)}} \\
\noalign{\smallskip}\hline\noalign{\smallskip}
\small{HVRP\_18\_DLP} & 256 & 9702.75 &  & 9797.61 & 1216.10 &  & \textbf{\textit{9652.74}} & -0.52 & \textbf{\textit{9688.96}} & 885.35 & -0.14 \\
\small{HVRP\_19\_DLP} & 224 & 11702.77 &  & 11805.34 & 1009.87 &  & \textbf{\textit{11686.39}} & -0.14 & 11745.69 & 274.81 & 0.37 \\
\small{HVRP\_22\_DLP} & 239 & 13068.03 &  & 13162.90 & 835.87 &  & 13091.16 & 0.18 & 13134.19 & 765.90 & 0.51 \\
\small{HVRP\_23\_DLP} & 203 & 7750.27 &  & 7809.20 & 802.30 &  & \textbf{\textit{7741.01}} & -0.12 & 7782.68 & 383.48 & 0.44 \\
\small{HVRP\_27\_DLP} & 220 & 8469.19 &  & 8520.74 & 995.85 &  & \textbf{\textit{8422.92}} & -0.55 & \textbf{\textit{8442.97}} & 372.15 & -0.32 \\
\small{HVRP\_32\_DLP} & 244 & 9417.62 &  & 9537.48 & 1131.44 &  & \textbf{\textit{9382.60}} & -0.37 & 9436.70 & 511.89 & 0.20 \\
\small{HVRP\_38\_DLP} & 205 & 11242.95 &  & 11439.58 & 421.50 &  & \textbf{\textit{11194.68}} & -0.43 & 11254.27 & 531.13 & 0.10 \\
\small{HVRP\_46\_DLP} & 250 & 24674.26 &  & 24805.27 & 1475.05 &  & \textbf{\textit{24566.23}} & -0.44 & 24698.60 & 495.57 & 0.10 \\
\small{HVRP\_49\_DLP} & 246 & 16377.69 &  & 16417.30 & 990.34 &  & \textbf{\textit{16181.17}} & -1.20 & \textbf{\textit{16322.51}} & 650.74 & -0.34 \\
\small{HVRP\_58\_DLP} & 220 & 23397.76 &  & 23530.10 & 1028.25 &  & \textbf{\textit{23370.42}} & -0.12 & 23641.18 & 294.93 & 1.04 \\
\small{HVRP\_62\_DLP} & 225 & 23149.61 &  & 23434.56 & 828.76 &  & \textbf{\textit{23010.35}} & -0.60 & \textbf{\textit{23097.54}} & 342.94 & -0.22 \\
\small{HVRP\_65\_DLP} & 223 & {13053.80} &  & 13077.63 & 635.64 &  & \textbf{\textit{13043.54}} & -0.08 & 13063.89 & 288.39 & 0.08 \\
\noalign{\smallskip}\hline\noalign{\smallskip}
 \multicolumn{2}{l}{Average} & \multicolumn{1}{l}{} & \multicolumn{1}{l}{} & 0.81 & 947.58 & \multicolumn{1}{l}{} & \multicolumn{1}{l}{} & -0.37 & \multicolumn{1}{l}{} & 483.11 & 0.15 \\
 \hline
\end{tabular}
} 
\label{tbl:Resultados.HH.HFFVRP_dlp_gg}
\end{table}

 \newpage
 \subsection{HFFOVRP} \label{sec:Result.HH.HFFOVRP}
Detailed results obtained for the instances of \cite[T99]{Taillard1999} as considered in \cite[YDR14]{Yousefikhoshbakht2014}. The results obtained by HILS-RVRP were compared with those found by the BRMTS heuristic of the referred authors. Moreover, although \citet{Yousefikhoshbakht2014} mentioned that they used fixed and variant vehicle costs, it appears, according to our testing, that they only used variable costs.
Table~\ref{tbl:Resultados.HH.HFFOVRP-V} presents the results involving only variable costs, while Table~\ref{tbl:Resultados.HH.HFFOVRP-FV} reports the results involving both fixed and variable costs.
\begin{table}[htbp]
\caption{Results for the HFFOVRP-V}
\setlength{\tabcolsep}{1.5mm}
\scalebox{0.95} {
\begin{tabular}{rrrlrrlrrrrr}
\hline\noalign{\smallskip}
 &  &  &  & \multicolumn{2}{l}{BRMTS}   & &\multicolumn{5}{l}{HILS-RVRP}  \\
 &  &  &  & \multicolumn{2}{l}{YDR14}   & &\multicolumn{5}{l}{}  \\
\cline{5-6} \cline{8-12}\noalign{\smallskip}
\multicolumn{1}{l}{Inst.} & \multicolumn{1}{c}{\textit{n}} & \multicolumn{1}{l}{BKS} &  & \multicolumn{1}{c}{Best Sol.} & \multicolumn{1}{c}{{Avg. T (s)}} & \multicolumn{1}{c}{} & \multicolumn{1}{c}{Best Sol.} & \multicolumn{1}{c}{Gap (\%)} & \multicolumn{1}{c}{{Avg. Sol.}} & \multicolumn{1}{c}{{Avg. T(s)}} & \multicolumn{1}{c}{{Avg. Gap (\%)}} \\ 
\noalign{\smallskip}\hline\noalign{\smallskip}
13 & 50 & 981.32 &  & 981.32 & 45.63 &  & \textbf{\textit{914.12}} & -6.85 & \textbf{\textit{914.12}} & 2.56 & -6.85 \\ 
14 & 50 & 448.25 &  & 448.25 & 38.72 &  & \textbf{\textit{436.32}} & -2.66 & \textbf{\textit{436.32}} & 1.82 & -2.66 \\ 
15 & 50 & 703.69 &  & 703.69 & 50.93 &  & \textbf{\textit{681.46}} & -3.16 & \textbf{\textit{681.71}} & 2.01 & -3.12 \\ 
16 & 50 & 788.12 &  & 788.12 & 60.34 &  & \textbf{\textit{769.46}} & -2.37 & \textbf{\textit{770.37}} & 2.26 & -2.25 \\ 
17 & 75 & 815.05 &  & 815.05 & 102.61 &  & \textbf{\textit{762.64}} & -6.43 & \textbf{\textit{763.60}} & 5.63 & -6.31 \\ 
18 & 75 & 1596.45 &  & 1596.45 & 159.54 &  & \textbf{\textit{1297.92}} & -18.70 & \textbf{\textit{1299.87}} & 6.73 & -18.58 \\
19 & 100 & 956.62 &  & 956.62 & 208.84 &  & \textbf{\textit{851.94}} & -10.94 & \textbf{\textit{853.89}} & 14.48 & -10.74 \\ 
20 & 100 & 1031.94 &  & \textbf{1031.94} & 267.14 &  & 1044.55 & 1.22 & 1045.69 & 15.02 & 1.33 \\ 
\noalign{\smallskip}\hline\noalign{\smallskip}
 \multicolumn{2}{l}{Average}  & & & 0.00 & 116.72 &  &  & -6.24 &  & 6.31 & -6.15 \\
 \hline
\end{tabular}
} 
\label{tbl:Resultados.HH.HFFOVRP-V}
\end{table}

\begin{table}[htbp]
\caption{Results for the HFFOVRP-FV}
\setlength{\tabcolsep}{1.5mm}
\scalebox{0.95} {
\begin{tabular}{rrrlrrrrr}
\hline\noalign{\smallskip}
 &  &  &  & \multicolumn{5}{l}{HILS-RVRP}  \\
\cline{5-9}\noalign{\smallskip}
 \multicolumn{1}{l}{Inst.} & \multicolumn{1}{c}{\textit{n}} & \multicolumn{1}{l}{BKS} &  &  \multicolumn{1}{c}{Best Sol.} & \multicolumn{1}{c}{Gap (\%)} & \multicolumn{1}{c}{{Avg. Sol.}} & \multicolumn{1}{c}{{Avg. T(s)}} & \multicolumn{1}{c}{{Avg. Gap (\%)}} \\
\noalign{\smallskip}\hline\noalign{\smallskip}
13 & 50 & 2588.65 &  & \textbf{\textit{2588.65}} & 0.00 & 2589.09 & 4.26 & 0.02 \\
14 & 50 & 9961.81 &  & \textbf{\textit{9961.81}} & 0.00 & 9968.21 & 3.74 & 0.06 \\
15 & 50 & 2731.46 &  & \textbf{\textit{2731.46}} & 0.00 & 2731.78 & 30.39 & 0.01 \\
16 & 50 & 2929.78 &  & \textbf{\textit{2929.78}} & 0.00 & 2962.68 & 29.06 & 1.12 \\
17 & 75 & 1792.20 &  & \textbf{\textit{1792.20}} & 0.00 & 1796.87 & 14.39 & 0.26 \\
18 & 75 & 3228.14 &  & \textbf{\textit{3228.14}} & 0.00 & 3236.72 & 31.99 & 0.27 \\
19 & 100 & 10179.70 &  & \textbf{\textit{10179.70}} & 0.00 & 10188.56 & 44.94 & 0.09 \\
20 & 100 & 4344.55 &  & \textbf{\textit{4344.55}} & 0.00 & 4351.39 & 42.39 & 0.16 \\
\noalign{\smallskip}\hline\noalign{\smallskip}
 \multicolumn{2}{l}{Average}  &  &  &  & 0.00 &  & 25.15 & 0.25 \\
 \hline
\end{tabular}
} 
\label{tbl:Resultados.HH.HFFOVRP-FV}
\end{table}

\subsection{Results for the {{MDFSMVRP}}} \label{sec:Result.HH.MDFSM}

Detailed results obtained for the {{MDFSMVRP}} instances  of \cite[SS97]{Salhi1997}, compared with those found by the VNS2 of \cite[SIW14]{Salhietal2014} and the HGSADC
of {\cite[VCGP14]{Vidaletal-MD-2014}} (Table~\ref{Tbl:Resultados.HH.FSMMD}). Column $RD$ in Table \ref{Tbl:Resultados.HH.FSMMD} indicates the instances with route duration.

\begin{table}[htbp]
\caption{Results for the {{MDFSMVRP}}}
\setlength{\tabcolsep}{1.1mm}
\scalebox{0.78} {
\begin{tabular}{llllcrrrrrrrlrrrrr}
\hline\noalign{\smallskip}
 &  &  &  &  &  &  & \multicolumn{ 2}{l}{VNS2} &  & \multicolumn{ 2}{l}{HGSADC} & & \multicolumn{ 5}{l}{HILS-RVRP} \\
 &  &  &  &  &  &  & \multicolumn{ 2}{l}{SIW14} &  &  \multicolumn{ 2}{l}{VCGP14} &  &  &  &  \\
\cline{8-9} \cline{11-12} \cline{14-18}\noalign{\smallskip}
\multicolumn{1}{c}{Inst.} & \multicolumn{1}{c}{\textit{n}} & \multicolumn{1}{c}{\textit{t}} & \multicolumn{1}{c}{\textit{m}} & \textit{RD} & \multicolumn{1}{c}{BKS} &  & \multicolumn{1}{c}{{Sol.}} & \multicolumn{1}{c}{T (s)} &  & \multicolumn{1}{c}{Best Sol.} & \multicolumn{1}{c}{{Avg. T (s)}} & \multicolumn{1}{c}{} & \multicolumn{1}{c}{Best Sol.} & \multicolumn{1}{c}{Gap (\%)} & \multicolumn{1}{c}{{Avg. Sol.}} & \multicolumn{1}{c}{{Avg. T(s)}} & \multicolumn{1}{c}{{Avg. Gap (\%)}} \\ 
\noalign{\smallskip}\hline\noalign{\smallskip}
p01 & 50  & 5 & 4  & \multicolumn{1}{c}{--} & 1477.73 &  & 1499.30 & 12.0 &  & \textbf{1477.73} & 112.2 &  & \textbf{1477.73} & 0.00 & 1494.76 & 4.31 & 1.15 \\ 
p02 & 50  & 5 & 4  & \multicolumn{1}{c}{--} & 957.73 &  & 984.50 & 12.0 &  & \textbf{957.73} & 146.4 &  & \textbf{957.73} & 0.00 & 966.81 & 4.12 & 0.95 \\ 
p03 & 75  & 5 & 5  & \multicolumn{1}{c}{--} & 1569.67 &  & 1588.40 & 30.0 &  & \textbf{1569.67} & 202.2 &  & \textbf{1569.67} & 0.00 & 1585.29 & 9.82 & 0.99 \\ 
p04 & 100  & 5 & 2  & \multicolumn{1}{c}{--} & 2292.64 &  & 2313.70 & 84.0 &  & \textbf{2292.64} & 278.4 &  & \textbf{2292.64} & 0.00 & 2313.42 & 25.41 & 0.91 \\ 
p05 & 100  & 5 & 2  & \multicolumn{1}{c}{--} & 1453.64 &  & 1466.90 & 108.0 &  & \textbf{1453.64} & 445.2 &  & \textbf{1453.64} & 0.00 & 1471.05 & 16.65 & 1.20 \\ 
p06 & 100  & 5 & 3  & \multicolumn{1}{c}{--} & 2208.66 &  & 2246.90 & 60.0 &  & \textbf{2208.66} & 320.4 &  & \textbf{2208.66} & 0.00 & 2229.31 & 25.29 & 0.93 \\ 
p07 & 100  & 5 & 4  & \multicolumn{1}{c}{--} & 2198.91 &  & 2256.40 & 66.0 &  & \textbf{2198.91} & 311.4 &  & \textbf{2198.91} & 0.00 & 2227.96 & 26.61 & 1.32 \\ 
p08 & 249  & 5 & 2  & 310 & 6441.36 &  & 6696.00 & 822.0 &  & \textbf{6441.36} & 1200.0 &  & 6448.93 & 0.12 & 6525.14 & 165.05 & 1.30 \\ 
p09 & 249  & 5 & 3  & 311 & 5998.70 &  & 6068.80 & 438.0 &  & \textbf{5998.70} & 1200.0 &  & 6011.85 & 0.22 & 6074.65 & 171.76 & 1.27 \\ 
p10 & 249  & 5 & 4  & 312 & 5807.53 &  & 6043.00 & 372.0 &  & \textbf{5807.53} & 1200.0 &  & 5826.61 & 0.33 & 5868.29 & 173.56 & 1.05 \\ 
p11 & 249  & 5 & 5  & 313 & 5770.42 &  & 5882.80 & 384.0 &  & \textbf{5770.42} & 1184.4 &  & 5779.64 & 0.16 & 5843.27 & 181.71 & 1.26 \\ 
p12 & 80  & 5 & 2  & \multicolumn{1}{c}{--} & 2072.18 &  & 2076.20 & 72.0 &  & \textbf{2072.18} & 216.0 &  & \textbf{2072.18} & 0.00 & 2075.22 & 6.96 & 0.15 \\ 
p13 & 80  & 5 & 2  & 200 & 2096.39 &  & 2096.40 & 66.0 &  & \textbf{2096.39} & 213.6 &  & \textbf{2096.39} & 0.00 & 2096.63 & 3.48 & 0.01 \\ 
p14 & 80  & 5 & 2  & 180 & 2139.30 &  & \textbf{2139.30} & 66.0 &  & 2160.12 & 250.8 &  & 2160.12 & 0.97 & 2169.66 & 2.89 & 1.42 \\ 
p15 & 160  & 5 & 4  & \multicolumn{1}{c}{--} & 3973.47 &  & 4024.90 & 162.0 &  & \textbf{3973.47} & 576.6 &  & \textbf{3973.47} & 0.00 & 3991.43 & 51.51 & 0.45 \\ 
p16 & 160  & 5 & 4  & 200 & 4119.76 &  & 4148.70 & 162.0 &  & \textbf{4119.76} & 595.8 &  & \textbf{4119.76} & 0.00 & 4128.48 & 23.31 & 0.21 \\ 
p17 & 160  & 5 & 4  & 180 & 4309.09 &  & 4338.20 & 198.0 &  & \textbf{4309.09} & 837.0 &  & \textbf{4309.09} & 0.00 & 4327.05 & 20.12 & 0.42 \\ 
p18 & 240  & 5 & 6  & \multicolumn{1}{c}{--} & 5887.43 &  & 5970.50 & 288.0 &  & \textbf{5887.43} & 1188.0 &  & \textbf{5887.43} & 0.00 & 5917.18 & 166.19 & 0.51 \\ 
p19 & 240  & 5 & 6  & 200 & 6130.36 &  & 6196.30 & 306.0 &  & \textbf{6130.36} & 1168.2 &  & \textbf{6130.36} & 0.00 & 6164.98 & 66.08 & 0.56 \\ 
p20 & 240  & 5 & 6  & 180 & 6469.21 &  & 6567.10 & 318.0 &  & \textbf{6469.21} & 1200.0 &  & \textbf{\textit{6458.07}} & -0.17 & 6486.35 & 56.60 & 0.26 \\ 
p21 & 360  & 5 & 9  & \multicolumn{1}{c}{--} & 8709.26 &  & 8883.10 & 654.0 &  & \textbf{8709.26} & 1200.0 &  & \textbf{8709.26} & 0.00 & 8738.45 & 513.06 & 0.34 \\ 
p22 & 360  & 5 & 9  & 200 & 9151.64 &  & 9294.80 & 468.0 &  & \textbf{9151.64} & 1200.0 &  & 9151.91 & 0.00 & 9209.14 & 210.11 & 0.63 \\ 
p23 & 360  & 5 & 9  & 180 & 9706.60 &  & 9887.70 & 540.0 &  & 9714.41 & 1200.0 &  & \textbf{\textit{9678.75}} & -0.29 & 9732.10 & 178.21 & 0.26 \\ 
\noalign{\smallskip}\hline\noalign{\smallskip}
 \multicolumn{3}{l}{Average}  &  &  &  &  & 1.52 & 247.30 &  & 0.05 & 715.07 &  &  & 0.06 &  & 91.43 & 0.76 \\
 \hline
\end{tabular}
} 
\label{Tbl:Resultados.HH.FSMMD}
\end{table}

\subsection{HFFVRPB} \label{sec:Result.HH.HFFVRPB}

Detailed results obtained for the HFFVRPB instances of \cite[T10]{Tutuncu2010}, compared with those found by the GRAMPS and ADVISER heuristics from the referred authors (Table~\ref{tbl:Resultados.HH.HFFVRPB} describes the results found). Note that we did not report the results for some instances because, according to the values suggested by the authors, they are infeasible.


\begin{table}[htbp]
\caption{Results for the HFFVRPB}
\setlength{\tabcolsep}{1.1mm}
\scalebox{0.82} {
\begin{tabular}{rrrrlrrrclrclrrrrr}
\hline\noalign{\smallskip} 
 &  &  &  &  &  &  & \multicolumn{2}{l}{GRAMPS} & & \multicolumn{2}{l}{ADVISER} &  & \multicolumn{5}{l}{HILS-RVRP} \\
\cline{8-9} \cline{11-12} \cline{14-18} \noalign{\smallskip}
\multicolumn{1}{c}{{Inst.}} & \multicolumn{1}{c}{{$n$}} & \multicolumn{1}{c}{{LH}} & \multicolumn{1}{c}{{BH}} & \multicolumn{1}{c}{\textbf{}} & \multicolumn{1}{c}{{BKS}} &  & \multicolumn{1}{c}{Best Sol.} & \multicolumn{1}{c}{T (s)} &  & \multicolumn{1}{c}{Best Sol.} & \multicolumn{1}{c}{T (s)} & \multicolumn{1}{c}{} & \multicolumn{1}{c}{Best Sol.} & \multicolumn{1}{c}{Gap (\%)} & \multicolumn{1}{c}{{Avg. Sol.}} & \multicolumn{1}{c}{{Avg. T(s)}} & \multicolumn{1}{c}{{Avg. Gap (\%)}} \\
\noalign{\smallskip}\hline\noalign{\smallskip}
1 & 50 & 25 & 25 &  & 1056.44 &  & 1111.67  & -- &  & \textbf{1056.44}  & -- &  & \textbf{\textit{874.60}} & -17.21 & \textbf{\textit{874.76}} & 0.89 & -17.20 \\
2 & 50 & 34 & 16 &  & 982.86 &  & 1067.28  & -- &  & \textbf{982.86}  & -- &  & \textbf{\textit{911.20}} & -7.29 & \textbf{\textit{913.30}} & 0.81 & -7.08 \\
3 & 50 & 40 & 10 &  & 998.22 &  & 1124.14  & -- &  & \textbf{998.22}  & -- &  & \multicolumn{1}{c}{--} & \multicolumn{1}{c}{--} & \multicolumn{1}{c}{--} & \multicolumn{1}{c}{--} & \multicolumn{1}{c}{--} \\
4 & 50 & 25 & 25 &  & 1070.06 &  & 1094.08  & -- &  & \textbf{1070.06}  & -- &  & \textbf{\textit{1050.60}} & -1.82 & \textbf{\textit{1051.11}} & 0.93 & -1.77 \\
5 & 50 & 34 & 16 &  & 1127.97 &  & 1135.21  & -- &  & \textbf{1127.97}  & -- &  & \textbf{\textit{1051.30}} & -6.80 & \textbf{\textit{1052.00}} & 0.83 & -6.74 \\
6 & 50 & 40 & 10 &  & 1183.36 &  & 1200.58  & -- &  & \textbf{1183.36}  & -- &  & \multicolumn{1}{c}{--} & \multicolumn{1}{c}{--} & \multicolumn{1}{c}{--} & \multicolumn{1}{c}{--} & \multicolumn{1}{c}{--} \\
7 & 75 & 37 & 38 &  & 1190.63 &  & \textbf{1190.63}  & -- &  & \textbf{1190.63}  & -- &  & \textbf{\textit{1073.90}} & -9.80 & \textbf{\textit{1096.80}} & 3.13 & -7.88 \\
8 & 75 & 50 & 25 &  & 1182.66 &  & 1211.28  & -- &  & \textbf{1182.66}  & -- &  & \multicolumn{1}{c}{--} & \multicolumn{1}{c}{--} & \multicolumn{1}{c}{--} & \multicolumn{1}{c}{--} & \multicolumn{1}{c}{--} \\
9 & 75 & 60 & 15 &  & 1203.09 &  & 1222.66  & -- &  & \textbf{1203.09}  & -- &  & \textbf{\textit{1003.20}} & -16.61 & \textbf{\textit{1013.97}} & 2.35 & -15.72 \\
10 & 75 & 37 & 38 &  & 1781.50 &  & 1845.75  & -- &  & \textbf{1781.50}  & -- &  & \textbf{\textit{1553.00}} & -12.83 & \textbf{\textit{1557.28}} & 2.74 & -12.59 \\
11 & 75 & 50 & 25 &  & 1941.74 &  & 2035.39  & -- &  & \textbf{1941.74}  & -- &  & \textbf{\textit{1659.80}} & -14.52 & \textbf{\textit{1667.85}} & 3.19 & -14.11 \\
12 & 75 & 60 & 15 &  & 1917.54 &  & 1945.35  & -- &  & \textbf{1917.54}  & -- &  & \multicolumn{1}{c}{--} & \multicolumn{1}{c}{--} & \multicolumn{1}{c}{--} & \multicolumn{1}{c}{--} & \multicolumn{1}{c}{--} \\
13 & 100 & 50 & 50 &  & 1227.81 &  & 1228.24  & -- &  & \textbf{1227.81}  & -- &  & \textbf{\textit{1181.70}} & -3.76 & \textbf{\textit{1195.00}} & 7.98 & -2.67 \\
14 & 100 & 67 & 33 &  & 1109.02 &  & 1136.87  & -- &  & \textbf{1109.02}  & -- &  & \multicolumn{1}{c}{--} & \multicolumn{1}{c}{--} & \multicolumn{1}{c}{--} & \multicolumn{1}{c}{--} & \multicolumn{1}{c}{--} \\
15 & 100 & 80 & 20 &  & 1216.65 &  & 1228.56  & -- &  & \textbf{1216.65}  & -- &  & \textbf{\textit{1114.90}} & -8.36 & \textbf{\textit{1135.00}} & 6.27 & -6.71 \\
16 & 100 & 50 & 50 &  & 1555.35 &  & 1629.47  & -- &  & \textbf{1555.35} & -- &  & \textbf{\textit{1314.50}} & -15.49 & \textbf{\textit{1323.97}} & 6.71 & -14.88 \\
\noalign{\smallskip}\hline\noalign{\smallskip}
\multicolumn{2}{l}{Average} &  &  &  &  &  & \multicolumn{1}{r}{3.31} &  &  & \multicolumn{1}{r}{0.00} &  &  &  & -10.41 &  & 3.26 & -9.76 \\
\hline
\end{tabular}
} 
\label{tbl:Resultados.HH.HFFVRPB}
\end{table}

\subsection{FSMB} \label{sec:Result.HH.FSMB}

Detailed results obtained for the FSMB instances of \cite[SWH13]{Salhietal2013}, compared with those found by \textit{Framework-2} from the same authors (Table~\ref{tbl:Resultados.HH.FSMB}).
\begin{table}[htbp]
\caption{Results for the FSMB}
\setlength{\tabcolsep}{1.5mm}
\scalebox{0.85} {
\begin{tabular}{lrrrrlrrlrrrrr}
\hline\noalign{\smallskip}
 &  &  &  &  &  & \multicolumn{ 2}{l}{Framework-2} &  & \multicolumn{ 5}{l}{HILS-RVRP} \\
 &  &  &  &  &  & \multicolumn{ 2}{l}{SWH13} &  &  &  &  &  &  \\
\cline{7-8} \cline{10-14}\noalign{\smallskip}
\multicolumn{1}{c}{\textbf{Inst.}} & \multicolumn{1}{c}{\textbf{$n$}} & \multicolumn{1}{c}{\textbf{$LH$}} & \multicolumn{1}{c}{\textbf{$BH$}} & \multicolumn{1}{c}{\textbf{BKS}} & \textbf{} & \multicolumn{1}{c}{{Sol.}} & \multicolumn{1}{c}{T (s)} & \multicolumn{1}{c}{} & \multicolumn{1}{c}{Best Sol.} & \multicolumn{1}{c}{Gap (\%)} & \multicolumn{1}{c}{{Avg. Sol.$^b$}} & \multicolumn{1}{c}{{Avg. T$^b$(s)}} & \multicolumn{1}{c}{{Avg. Gap$^b$ (\%)}} \\
\noalign{\smallskip}\hline\noalign{\smallskip}
HWS1 & 20 & 10 & 10 & 720.57 & $^1$ & 726.48 & 0.61 &  & \textbf{720.57} & 0.00 & \textbf{720.57} & 0.12 & 0.00 \\
HWS2 & 20 & 13 & 7 & 818.12 & $^1$ & \textbf{818.12} & 1.09 &  & \textbf{818.12} & 0.00 & \textbf{818.12} & 0.14 & 0.00 \\
HWS3 & 20 & 16 & 4 & 848.23 & $^1$ & 848.59 & 1.64 &  & 848.32 & 0.01 & 848.32 & 0.10 & 0.01 \\
HWS4 & 20 & 10 & 10 & 4342.48 & $^1$ & 4350.65 & 0.91 &  & \textbf{4342.48} & 0.00 & \textbf{4342.48} & 0.12 & 0.00 \\
HWS5 & 20 & 13 & 7 & 5357.98 & $^1$ & 5366.39 & 2.75 &  & \textbf{5357.98} & 0.00 & 5361.47 & 0.12 & 0.07 \\
HWS6 & 20 & 16 & 4 & 5421.65 & $^1$ & 5875.23 & 3.44 &  & \textbf{5421.65} & 0.00 & 5560.13 & 0.14 & 2.55 \\
HWS7 & 20 & 10 & 10 & 729.50 & $^1$ & 767.93 & 0.58 &  & \textbf{729.50} & 0.00 & \textbf{729.50} & 0.14 & 0.00 \\
HWS8 & 20 & 13 & 7 & 838.11 & $^1$ & 872.97 & 1.39 &  & \textbf{838.11} & 0.00 & 838.20 & 0.16 & 0.01 \\
HWS9 & 20 & 16 & 4 & 890.76 & $^1$ & 903.18 & 2.09 &  & \textbf{890.76} & 0.00 & \textbf{890.76} & 0.13 & 0.00 \\
HWS10 & 20 & 10 & 10 & 4349.12 & $^1$ & 4365.44 & 0.88 &  & \textbf{4349.12} & 0.00 & \textbf{4349.12} & 0.12 & 0.00 \\
HWS11 & 20 & 13 & 7 & 5363.58 & $^1$ & 5414.50 & 2.72 &  & \textbf{5363.58} & 0.00 & 5380.38 & 0.15 & 0.31 \\
HWS12 & 20 & 16 & 4 & 5497.98 & $^1$ & 5928.78 & 4.94 &  & \textbf{5497.98} & 0.00 & 5751.82 & 0.15 & 4.62 \\
HWS13 & 50 & 25 & 25 & 1625.70 &  & \textbf{1625.70} & 17.88 &  & \textbf{\textit{1590.47}} & -2.17 & \textbf{\textit{1593.28}} & 1.67 & -1.99 \\
HWS14 & 50 & 33 & 17 & 1811.63 &  & \textbf{1811.63} & 26.19 &  & \textbf{\textit{1771.53}} & -2.21 & \textbf{\textit{1778.19}} & 1.37 & -1.85 \\
HWS15 & 50 & 40 & 10 & 2018.93 &  & \textbf{2018.93} & 38.42 &  & \textbf{\textit{1999.05}} & -0.98 & \textbf{\textit{2004.69}} & 1.21 & -0.71 \\
HWS16 & 50 & 25 & 25 & 5561.67 &  & \textbf{5561.67} & 330.34 &  & \textbf{\textit{5551.19}} & -0.19 & \textbf{\textit{5551.34}} & 1.53 & -0.19 \\
HWS17 & 50 & 33 & 17 & 6570.39 & & \textbf{6570.39} & 996.55 &  & \textbf{\textit{6547.93}} & -0.34 & \textbf{\textit{6547.93}} & 2.29 & -0.34 \\
HWS18 & 50 & 40 & 10 & 7599.08 &  & \textbf{7599.08} & 1120.50 &  & \textbf{\textit{7120.52}} & -6.30 & \textbf{\textit{7523.17}} & 2.86 & -1.00 \\
HWS19 & 50 & 25 & 25 & 1704.41 &  & \textbf{1704.41} & 39.81 &  & \textbf{\textit{1616.21}} & -5.17 & \textbf{\textit{1627.49}} & 1.46 & -4.51 \\
HWS20 & 50 & 33 & 17 & 2037.23 &  & \textbf{2037.23} & 84.95 &  & \textbf{\textit{2015.67}} & -1.06 & \textbf{\textit{2018.79}} & 2.32 & -0.91 \\
HWS21 & 50 & 40 & 10 & 2340.09 &  & \textbf{2340.09} & 103.52 &  & \textbf{\textit{2295.57}} & -1.90 & \textbf{\textit{2304.64}} & 2.67 & -1.51 \\
HWS22 & 50 & 25 & 25 & 1774.71 &  & \textbf{1774.71} & 18.41 &  & \textbf{\textit{1717.60}} & -3.22 & \textbf{\textit{1722.60}} & 1.79 & -2.94 \\
HWS23 & 50 & 33 & 17 & 2166.52 &  & \textbf{2166.52} & 64.77 &  & \textbf{\textit{2096.10}} & -3.25 & \textbf{\textit{2127.37}} & 1.93 & -1.81 \\
HWS24 & 50 & 40 & 10 & 2430.88 &  & \textbf{2430.88} & 49.72 &  & \textbf{\textit{2401.04}} & -1.23 & \textbf{\textit{2407.88}} & 1.39 & -0.95 \\
HWS25 & 75 & 37 & 38 & 1332.02 &  & \textbf{1332.02} & 1006.28 &  & \textbf{\textit{1285.86}} & -3.47 & \textbf{\textit{1292.21}} & 5.66 & -2.99 \\
HWS26 & 75 & 50 & 25 & 1421.04 &  & \textbf{1421.04} & 1779.88 &  & \textbf{\textit{1399.36}} & -1.53 & \textbf{\textit{1401.82}} & 5.19 & -1.35 \\
HWS27 & 75 & 60 & 15 & 1534.65 &  & \textbf{1534.65} & 1996.59 &  & \textbf{\textit{1513.10}} & -1.40 & \textbf{\textit{1524.69}} & 4.60 & -0.65 \\
HWS28 & 75 & 37 & 38 & 1617.85 &  & \textbf{1617.85} & 1351.92 &  & \textbf{\textit{1572.38}} & -2.81 & \textbf{\textit{1574.08}} & 5.24 & -2.71 \\
HWS29 & 75 & 50 & 25 & 1799.76 &  & \textbf{1799.76} & 1513.30 &  & \textbf{\textit{1760.95}} & -2.16 & \textbf{\textit{1761.04}} & 4.45 & -2.15 \\
HWS30 & 75 & 60 & 15 & 1990.46 &  & \textbf{1990.46} & 2662.15 &  & \textbf{\textit{1950.99}} & -1.98 & \textbf{\textit{1951.30}} & 3.80 & -1.97 \\
HWS31 & 100 & 50 & 50 & 4943.29 &  & 5201.81 & 4257.41 &  & 4963.08 & 0.40 & 4966.57 & 12.66 & 0.47 \\
HWS32 & 100 & 66 & 34 & 6035.96 &  & \multicolumn{1}{c}{-} & \multicolumn{1}{c}{-} &  & \textbf{\textit{5993.30}} & -0.71 & \textbf{\textit{5993.91}} & 36.35 & -0.70 \\
HWS33 & 100 & 80 & 20 & 7601.09 &  & \multicolumn{1}{c}{-} & \multicolumn{1}{c}{-} &  & \textbf{\textit{7097.81}} & -6.62 & \textbf{\textit{7330.43}} & 32.22 & -3.56 \\
HWS34 & 100 & 50 & 50 & 2465.41 &  & 2646.52 & 2871.74 &  & 2494.95 & 1.20 & 2522.53 & 17.35 & 2.32 \\
HWS35 & 100 & 66 & 34 & 2971.98 &  & \textbf{2971.98} & 651.80 &  & \textbf{\textit{2927.20}} & -1.51 & \textbf{\textit{2931.90}} & 15.01 & -1.35 \\
HWS36 & 100 & 80 & 20 & 3533.90 &  & \textbf{3533.90} & 1729.30 &  & \textbf{\textit{3450.73}} & -2.35 & \textbf{\textit{3458.24}} & 15.40 & -2.14 \\
\noalign{\smallskip}\hline\noalign{\smallskip}
\multicolumn{2}{l}{Average}  &  &  &  &  & 1.24 & 668.66 &  &  & -1.42 &  & 5.06 & -0.77 \\
\noalign{\smallskip}\hline\noalign{\smallskip}
\multicolumn{14}{l}{{\scriptsize $^1$: Optimality proven.}} \\
\end{tabular}
} 
\label{tbl:Resultados.HH.FSMB}
\end{table}

\newpage
\subsection{SDepVRP} \label{sec:Result.HH.SDepVRP}

Detailed results obtained for the SDepVRP instances of \cite[CL01]{CordeauLaporte2001}, compared with those found by the ALNS 50k of \cite[PR07]{PisingerRopke2007} and the ITS of \cite[CM12]{CordeauMaischberger2012} (Table~\ref{tbl:Resultados.HH.SDepVRPold}, old instances without route duration and Table~\ref{tbl:Resultados.HH.SDepVRPnew}, new instances with route duration).

\begin{table}[htbp]
\caption{Results for the SDepVRP (old set)}
\setlength{\tabcolsep}{1.4mm}
\scalebox{0.8} {
\begin{tabular}{lcrrlrrlrclrrrrr}
\hline\noalign{\smallskip}
 &  &  &  &  & \multicolumn{2}{l}{ALNS 50k} & & \multicolumn{2}{l}{ITS} & & \multicolumn{5}{l}{HILS-RVRP} \\
 &  &  &  &  & \multicolumn{1}{l}{PR07} &  &  & \multicolumn{1}{l}{CM12} &  &  &  &  &  &  &  \\
\cline{6-7} \cline{9-10} \cline{12-16}\noalign{\smallskip}
\multicolumn{1}{c}{Inst.} & \textit{n} & \multicolumn{1}{c}{\textit{t}} & \multicolumn{1}{c}{BKS} &  & \multicolumn{1}{c}{Best Sol.} & \multicolumn{1}{c}{{Avg. T (s)}} &  & \multicolumn{1}{c}{Best Sol.} & \multicolumn{1}{c}{T (s)} & \multicolumn{1}{c}{} & \multicolumn{1}{c}{Best Sol.} & \multicolumn{1}{c}{Gap (\%)} & \multicolumn{1}{c}{{Avg. Sol.}} & \multicolumn{1}{c}{{Avg. T(s)}} & \multicolumn{1}{c}{{Avg. Gap (\%)}} \\
\noalign{\smallskip}\hline\noalign{\smallskip}
p01 & 55  & 3 & 640.32 &  & \textbf{640.32} & 20 &  & \textbf{640.32} & -- &  & \textbf{640.32} & 0.00 & 640.53 & 1.00 & 0.03 \\
p02 & 52  & 2 & 598.10 &  & \textbf{598.10} & 19 &  & \textbf{598.10} & -- &  & \textbf{598.10} & 0.00 & \textbf{598.10} & 0.93 & 0.00 \\
p03 & 80  & 3 & 954.32 &  & 957.04 & 40 &  & \textbf{954.32} & -- &  & \textbf{954.32} & 0.00 & 956.37 & 4.26 & 0.22 \\
p04 & 76  & 2 & 854.43 &  & \textbf{854.43} & 36 &  & \textbf{854.43} & -- &  & \textbf{854.43} & 0.00 & 855.30 & 3.03 & 0.10 \\
p05 & 103  & 3 & 1003.57 &  & \textbf{1003.57} & 68 &  & \textbf{1003.57} & -- &  & \textbf{1003.57} & 0.00 & 1008.67 & 11.02 & 0.51 \\
p06 & 104  & 2 & 1028.52 &  & \textbf{1028.52} & 69 &  & \textbf{1028.52} & -- &  & \textbf{1028.52} & 0.00 & 1036.96 & 7.85 & 0.82 \\
p07 & 27  & 3 & 391.30 &  & \textbf{391.30} & 8 &  & \textbf{391.30} & -- &  & \textbf{391.30} & 0.00 & \textbf{391.30} & 0.11 & 0.00 \\
p08 & 54  & 3 & 664.46 &  & \textbf{664.46} & 24 &  & \textbf{664.46} & -- &  & \textbf{664.46} & 0.00 & \textbf{664.46} & 0.69 & 0.00 \\
p09 & 81  & 3 & 948.23 &  & \textbf{948.23} & 47 &  & \textbf{948.23} & -- &  & \textbf{948.23} & 0.00 & \textbf{948.23} & 4.19 & 0.00 \\
p10 & 108  & 3 & 1218.75 &  & \textbf{1218.75} & 76 &  & \textbf{1218.75} & -- &  & \textbf{1218.75} & 0.00 & 1231.41 & 11.21 & 1.04 \\
p11 & 135  & 3 & 1448.17 &  & 1463.33 & 116 &  & \textbf{1448.17} & -- &  & \textbf{1448.17} & 0.00 & 1482.58 & 49.62 & 2.55 \\
p12 & 162  & 3 & 1665.55 &  & 1678.40 & 157 &  & \textbf{1665.55} & -- &  & 1682.94 & 1.04 & 1708.05 & 64.58 & 2.55 \\
p13 & 54  & 3 & 1194.18 &  & \textbf{1194.18} & 24 &  & \textbf{1194.18} & -- &  & \textbf{1194.18} & 0.00 & 1196.12 & 0.93 & 0.16 \\
p14 & 108  & 3 & 1959.96 &  & 1960.62 & 72 &  & \textbf{1959.96} & -- &  & \textbf{1959.96} & 0.00 & 1960.90 & 7.59 & 0.05 \\
p15 & 162  & 3 & 2685.09 &  & \textbf{2685.09} & 152 &  & \textbf{2685.09} & -- &  & \textbf{2685.09} & 0.00 & 2701.85 & 33.21 & 0.62 \\
p16 & 216  & 3 & 3393.55 &  & 3396.36 & 213 &  & \textbf{3393.55} & -- &  & 3393.86 & 0.01 & 3431.81 & 83.85 & 1.13 \\
p17 & 270  & 3 & 4066.15 &  & 4085.61 & 291 &  & \textbf{4066.15} & -- &  & 4078.19 & 0.30 & 4147.14 & 260.77 & 1.99 \\
p18 & 324  & 3 & 4751.27 &  & 4755.50 & 346 &  & \textbf{4751.27} & -- &  & 4768.23 & 0.36 & 4910.04 & 505.80 & 3.34 \\
p19 & 104  & 3 & 843.15 &  & 846.07 & 85 &  & \textbf{843.15} & -- &  & \textbf{843.15} & 0.00 & 848.59 & 8.32 & 0.64 \\
p20 & 156  & 3 & 1030.78 &  & \textbf{1030.78} & 168 &  & \textbf{1030.78} & -- &  & \textbf{1030.78} & 0.00 & 1044.11 & 24.30 & 1.29 \\
p21 & 209  & 3 & 1263.71 &  & 1271.75 & 217 &  & \textbf{1263.71} & -- &  & \textbf{\textit{1260.01}} & -0.29 & 1278.05 & 61.59 & 1.13 \\
p22 & 122  & 3 & 1008.71 &  & \textbf{1008.71} & 130 &  & \textbf{1008.71} & -- &  & \textbf{1008.71} & 0.00 & 1009.09 & 14.02 & 0.04 \\
p23 & 102  & 3 & 803.29 &  & \textbf{803.29} & 73 &  & \textbf{803.29} & -- &  & \textbf{803.29} & 0.00 & 805.30 & 5.62 & 0.25 \\
\noalign{\smallskip}\hline\noalign{\smallskip}
\multicolumn{2}{l}{Average}  &  &  &  &  0.16 & 106.57 &  & 0.00 &  &  &  & 0.06 &  & 89.22 & 0.80 \\
\hline
\end{tabular}
} 
\label{tbl:Resultados.HH.SDepVRPold}
\end{table}

\begin{table}[htbp]
\caption{Results for the SDepVRP (new instances)}
\setlength{\tabcolsep}{1.28mm}
\scalebox{0.8} {
\begin{tabular}{lcrrlrrlrclrrrrr}
\hline\noalign{\smallskip}
 &  &  &  &  & \multicolumn{2}{l}{ALNS 50k} & & \multicolumn{2}{l}{ITS} & & \multicolumn{5}{l}{HILS-RVRP} \\
 &  &  &  &  & \multicolumn{1}{l}{PR07} &  &  & \multicolumn{1}{l}{CM12} &  &  &  &  &  &  &  \\
\cline{6-7} \cline{9-10} \cline{12-16}\noalign{\smallskip}
\multicolumn{1}{c}{Inst.} & \textit{n} & \multicolumn{1}{c}{\textit{t}} & \multicolumn{1}{c}{BKS} &  & \multicolumn{1}{c}{Best Sol.} & \multicolumn{1}{c}{{Avg. T (s)}} &  & \multicolumn{1}{c}{Best Sol.} & \multicolumn{1}{c}{T (s)} & \multicolumn{1}{c}{} & \multicolumn{1}{c}{Best Sol.} & \multicolumn{1}{c}{Gap (\%)} & \multicolumn{1}{c}{{Avg. Sol.}} & \multicolumn{1}{c}{{Avg. T(s)}} & \multicolumn{1}{c}{{Avg. Gap (\%)}} \\
\noalign{\smallskip}\hline\noalign{\smallskip}
pr01 & 48  & 4 & 1380.77 &  & \textbf{1380.77} & 19 &  & \textbf{1380.77} & -- &  & \textbf{1380.77} & 0.00 & 1408.22 & 0.62 & 1.99 \\
pr02 & 96  & 4 & 2303.89 &  & 2311.54 & 63 &  & \textbf{2303.89} & -- &  & \textbf{2303.89} & 0.00 & 2335.95 & 5.32 & 1.39 \\
pr03 & 144  & 4 & 2575.36 &  & 2602.13 & 140 &  & \textbf{2575.36} & -- &  & 2578.93 & 0.14 & 2598.95 & 16.16 & 0.92 \\
pr04 & 192  & 4 & 3449.84 &  & 3474.01 & 191 &  & \textbf{3449.84} & -- &  & 3454.85 & 0.15 & 3524.97 & 44.90 & 2.18 \\
pr05 & 240  & 4 & 4377.35 &  & 4416.38 & 251 &  & \textbf{4377.35} & -- &  & 4379.41 & 0.05 & 4481.48 & 77.01 & 2.38 \\
pr06 & 288  & 4 & 4422.02 &  & 4444.52 & 314 &  & \textbf{4422.02} & -- &  & 4451.34 & 0.66 & 4531.44 & 141.38 & 2.47 \\
pr07 & 72  & 6 & 1889.82 &  & \textbf{1889.82} & 39 &  & \textbf{1889.82} & -- &  & \textbf{1889.82} & 0.00 & 1935.09 & 2.53 & 2.40 \\
pr08 & 144  & 6 & 2971.01 &  & 2977.50 & 135 &  & \textbf{2971.01} & -- &  & 2973.26 & 0.08 & 3057.03 & 17.43 & 2.90 \\
pr09 & 216  & 6 & 3536.20 &  & \textbf{3536.20} & 226 &  & \textbf{3536.20} & -- &  & 3559.69 & 0.66 & 3624.11 & 64.06 & 2.49 \\
pr10 & 288  & 6 & 4639.62 &  & 4648.76 & 322 &  & \textbf{4639.62} & -- &  & 4663.38 & 0.51 & 4743.69 & 160.19 & 2.24 \\
pr11 & 1008  & 4 & 12719.65 &  & \textbf{12719.65} & 847 &  & 12845.60 & -- &  & 13091.91 & 2.93 & 13305.68 & 6775.12 & 4.61 \\
pr12 & 720  & 6 & 9388.07 &  & \textbf{9388.07} & 658 &  & 9392.84 & -- &  & 9612.09 & 2.39 & 9665.67 & 2872.05 & 2.96 \\
\noalign{\smallskip}\hline\noalign{\smallskip}
 \multicolumn{2}{l}{Average}  &  &  &  & {0.32}  & {267.08}  &  & 0.09 &  &  &  & 0.63 &  & 848.16 & 2.39 \\
 \hline
\end{tabular}
} 
\label{tbl:Resultados.HH.SDepVRPnew}
\end{table}

\newpage
\subsection{HFFVRPSD} \label{sec:Result.HFFVRPSD}

Detailed results obtained for the FSM instances of \cite[G84]{Golden1984} and adapted for the \mbox{HFFVRPSD} by \cite[OO10]{Ozfirat2010}, compared with those found by the CP heuristic of the same authors (Table~\ref{tbl:Resultados.HFVRPSD}). As HFFVRPSD allows visiting the customers more than once, the SP procedure was not considered for this variant.

\begin{table}[htbp]
\caption{Results for the HFFVRPSD} 
\setlength{\tabcolsep}{1.5mm}
\scalebox{0.95} {
\begin{tabular}{rrrlrrlrlrrr}
\hline\noalign{\smallskip}
 &  &  &  & \multicolumn{ 2}{l}{CP} &  & \multicolumn{5}{l}{HILS-RVRP}   \\
 &  &  &  & \multicolumn{2}{l}{OO10} &  &  &  &  &  &  \\
\cline{5-6} \cline{8-12}\noalign{\smallskip}
\multicolumn{1}{c}{{Inst.}} & \multicolumn{1}{c}{{\textit{n}}} & \multicolumn{1}{c}{{BKS}} & \multicolumn{1}{c}{} & \multicolumn{1}{c}{{Sol.}} & \multicolumn{1}{c}{{T (s)}} &  & \multicolumn{1}{c}{{Best Sol.}} & \multicolumn{1}{c}{{Gap}} & \multicolumn{1}{c}{{Avg. Sol.}} & \multicolumn{1}{c}{{Avg. T}} & \multicolumn{1}{c}{{Avg. Gap }} \\
\noalign{\smallskip}\hline\noalign{\smallskip}
3 & 20 & 970.53 &  & \textbf{970.53} & 104 &  & \textbf{\textit{943.12}} & -2.82 & \textbf{\textit{943.12}} & 11.03 & -2.82 \\
4 & 20 & 6421.88 &  & \textbf{6421.88} & 45 &  & \textbf{\textit{6399.20}} & -0.35 & \textbf{\textit{6400.01}} & 5.20 & -0.34 \\
5 & 20 & 998.74 &  & \textbf{998.74} & 84 &  & \textbf{\textit{970.77}} & -2.80 & \textbf{\textit{970.77}} & 10.11 & -2.80 \\
6 & 20 & 6514.09 &  & \textbf{6514.09} & 137 &  & \textbf{\textit{6497.44}} & -0.26 & 6582.88 & 4.54 & 1.06 \\
13 & 50 & 2440.78 &  & \textbf{2440.78} & 146 &  & \textbf{\textit{2374.03}} & -2.73 & \textbf{\textit{2376.20}} & 198.26 & -2.65 \\
14 & 50 & 9138.25 &  & \textbf{9138.25} & 412 &  & \textbf{\textit{9114.07}} & -0.26 & 9203.83 & 90.81 & 0.72 \\
15 & 50 & 2616.11 &  & \textbf{2616.11} & 298 &  & \textbf{\textit{2568.96}} & -1.80 & \textbf{\textit{2569.86}} & 106.94 & -1.77 \\
16 & 50 & 2719.89 &  & \textbf{2719.89} & 302 &  & \textbf{\textit{2690.18}} & -1.09 & \textbf{\textit{2691.51}} & 98.02 & -1.04 \\
17 & 75 & 1783.33 &  & \textbf{1783.33} & 317 &  & \textbf{\textit{1725.40}} & -3.25 & \textbf{\textit{1726.62}} & 385.42 & -3.18 \\
18 & 75 & 2394.16 &  & \textbf{2394.16} & 486 &  & \textbf{\textit{2333.26}} & -2.54 & \textbf{\textit{2336.85}} & 781.51 & -2.39 \\
19 & 100 & 8722.49 &  & \textbf{8722.49} & 506 &  & \textbf{\textit{8654.31}} & -0.78 & \textbf{\textit{8655.85}} & 845.61 & -0.76 \\
20 & 100 & 4130.49 &  & \textbf{4130.49} & 396 &  & \textbf{\textit{4011.30}} & -2.89 & \textbf{\textit{4033.08}} & 803.84 & -2.36 \\
\noalign{\smallskip}\hline\noalign{\smallskip}
\multicolumn{2}{l}{Average} & \multicolumn{1}{l}{} &  & 0.00 & 269.42 &  & \multicolumn{1}{l}{} & -1.80 & \multicolumn{1}{l}{} & 278.44 & -1.53 \\
\hline
\end{tabular}
}
\label{tbl:Resultados.HFVRPSD}
\end{table}

\newpage
\subsection{FSMVRPTW -- {{minimizing}} route duration} \label{sec:Result.HH.FSMTW-duration}


Detailed results obtained for the FSMVRPTW instances of \cite[LS99]{LiuShen1999}, considering the objective of minimizing the sum of the route durations, compared with those found by the UHGS of \cite[VCGP14]{Vidaletal2014} and the HEA of \cite[KBJL15]{Kocetal2015} (Tables \ref{tbl:Resultados.HH.FSMTW_Al}-\ref{tbl:Resultados.HH.FSMTW_Cl}).





\begin{table}[htbp]
\caption{Results for the FSMTW (\small{minimize duration, fleet A})}
\setlength{\tabcolsep}{1.5mm}
\scalebox{0.78} {
\begin{tabular}{lrrlrrlrrlrrrrr}
\hline\noalign{\smallskip} 
 &  &  &  & \multicolumn{2}{l}{VCGP14} &  & \multicolumn{2}{l}{KBJL15} &  & \multicolumn{5}{l}{HILS-RVRP} \\
\cline{5-6} \cline{8-9} \cline{11-15}\noalign{\smallskip}
\multicolumn{1}{c}{Inst.} & \multicolumn{1}{c}{n} & \multicolumn{1}{c}{BKS} &  & \multicolumn{1}{c}{Best Sol.} & \multicolumn{1}{c}{{Avg. T (s)}} &  & \multicolumn{1}{c}{Best Sol.} & {{Avg. T (s)}} & \multicolumn{1}{c}{} & \multicolumn{1}{c}{Best Sol.} & \multicolumn{1}{c}{Gap (\%)} & \multicolumn{1}{c}{{Avg. Sol.}} & \multicolumn{1}{c}{{Avg. T(s)}} & \multicolumn{1}{c}{{Avg. Gap (\%)}} \\
\noalign{\smallskip}\hline\noalign{\smallskip}
R101a  & 100 & 4536.40 &  & 4608.62 & 361.8 &  & 4541.70 & 315.6 &  & 4608.62 & 1.59 & 4616.92 & 62.56 & 1.77 \\ 
R102a  & 100 & 4348.92 &  & 4369.74 & 382.8 & \multicolumn{1}{c}{} & 4355.10 & 352.2 &  & 4368.70 & 0.45 & 4380.77 & 76.74 & 0.73 \\ 
R103a  & 100 & 4119.04 &  & 4145.68 & 279.0 & \multicolumn{1}{c}{} & 4131.23 & 251.4 &  & 4144.96 & 0.63 & 4152.64 & 77.66 & 0.82 \\ 
R104a  & 100 & 3961.39 &  & \textbf{3961.39} & 313.8 & \multicolumn{1}{c}{} & 3992.10 & 301.2 &  & 3962.87 & 0.04 & 3971.32 & 71.20 & 0.25 \\ 
R105a  & 100 & 4209.84 &  & \textbf{4209.84} & 306.6 &  & 4232.54 & 283.8 &  & \textbf{\textit{4204.87}} & -0.12 & 4222.45 & 76.06 & 0.30 \\ 
R106a  & 100 & 4109.08 &  & \textbf{4109.08} & 379.2 &  & 4138.30 & 307.8 &  & \textbf{\textit{4103.65}} & -0.13 & 4125.17 & 84.09 & 0.39 \\ 
R107a  & 100 & 4007.87 &  & \textbf{4007.87} & 327.0 &  & 4034.32 & 324.0 &  & \textbf{\textit{4003.95}} & -0.10 & 4014.98 & 75.44 & 0.18 \\ 
R108a  & 100 & 3934.48 &  & \textbf{3934.48} & 307.2 &  & 3966.10 & 286.8 &  & \textbf{3934.48} & 0.00 & 3946.40 & 69.80 & 0.30 \\ 
R109a  & 100 & 4020.75 &  & \textbf{4020.75} & 303.6 &  & 4059.02 & 276.0 &  & 4022.90 & 0.05 & 4034.21 & 78.46 & 0.33 \\ 
R110a  & 100 & 3965.88 &  & \textbf{3965.88} & 318.0 &  & 3996.31 & 250.2 &  & \textbf{\textit{3963.80}} & -0.05 & 3973.92 & 75.62 & 0.20 \\ 
R111a  & 100 & 3985.68 &  & \textbf{3985.68} & 363.0 &  & 4020.10 & 298.8 &  & \textbf{\textit{3984.92}} & -0.02 & 3999.81 & 83.57 & 0.35 \\ 
R112a  & 100 & 3918.88 &  & \textbf{3918.88} & 406.2 &  & 3957.60 & 346.8 &  & \textbf{3918.88} & 0.00 & 3927.38 & 68.16 & 0.22 \\ 
\noalign{\smallskip}\hline\noalign{\smallskip}
C101a  & 100 & 7226.51 &  & \textbf{7226.51} & 191.4 &  & \textbf{7226.51} & 178.2 &  & \textbf{7226.51} & 0.00 & 7231.05 & 83.95 & 0.06 \\ 
C102a  & 100 & 7119.35 &  & \textbf{7119.35} & 169.2 &  & 7145.65 & 186.0 &  & \textbf{7119.35} & 0.00 & \textbf{7119.35} & 90.42 & 0.00 \\ 
C103a  & 100 & 7102.86 &  & \textbf{7102.86} & 147.6 &  & 7143.88 & 162.0 &  & \textbf{7102.86} & 0.00 & 7103.02 & 86.13 & 0.00 \\ 
C104a  & 100 & 7081.50 &  & 7081.51 & 133.2 &  & 7082.92 & 120.6 &  & 7081.51 & 0.00 & 7081.51 & 73.36 & 0.00 \\ 
C105a  & 100 & 7175.00 &  & 7196.06 & 204.0 &  & \textbf{7175.00} & 147.0 &  & 7196.06 & 0.29 & 7196.39 & 98.76 & 0.30 \\ 
C106a  & 100 & 7163.32 &  & 7176.68 & 220.2 &  & \textbf{7163.32} & 180.6 &  & 7176.68 & 0.19 & 7180.15 & 98.72 & 0.23 \\ 
C107a  & 100 & 7140.20 &  & 7144.49 & 191.4 &  & \textbf{7140.20} & 166.8 &  & 7144.49 & 0.06 & 7144.53 & 97.32 & 0.06 \\ 
C108a  & 100 & 7111.23 &  & \textbf{7111.23} & 166.8 &  & 7120.98 & 147.0 &  & \textbf{7111.23} & 0.00 & \textbf{7111.23} & 87.48 & 0.00 \\ 
C109a  & 100 & 7091.66 &  & \textbf{7091.66} & 143.4 &  & \textbf{7091.66} & 142.2 &  & \textbf{7091.66} & 0.00 & \textbf{7091.66} & 82.51 & 0.00 \\ 
\noalign{\smallskip}\hline\noalign{\smallskip}
RC101a  & 100 & 5217.90 &  & \textbf{5217.90} & 301.8 &  & 5235.42 & 298.2 &  & \textbf{\textit{5213.66}} & -0.08 & 5227.92 & 61.18 & 0.19 \\ 
RC102a  & 100 & 5018.47 &  & \textbf{5018.47} & 342.6 &  & 5029.69 & 338.4 &  & 5019.99 & 0.03 & 5042.64 & 66.00 & 0.48 \\ 
RC103a  & 100 & 4822.21 &  & \textbf{4822.21} & 361.8 &  & 4870.00 & 308.4 &  & \textbf{4822.21} & 0.00 & 4859.48 & 66.50 & 0.77 \\ 
RC104a  & 100 & 4737.00 &  & \textbf{4737.00} & 246.0 &  & 4769.30 & 298.2 &  & \textbf{\textit{4736.13}} & -0.02 & 4748.10 & 54.96 & 0.23 \\ 
RC105a  & 100 & 5097.35 &  & \textbf{5097.35} & 336.6 &  & 5118.10 & 319.2 &  & \textbf{\textit{5096.66}} & -0.01 & 5111.72 & 70.68 & 0.28 \\ 
RC106a  & 100 & 4935.91 &  & \textbf{4935.91} & 396.0 &  & 4958.62 & 360.6 &  & \textbf{\textit{4928.60}} & -0.15 & 4946.41 & 65.98 & 0.21 \\ 
RC107a  & 100 & 4783.08 &  & \textbf{4783.08} & 319.2 &  & 4825.21 & 322.2 &  & \textbf{4783.08} & 0.00 & 4803.81 & 61.07 & 0.43 \\ 
RC108a  & 100 & 4708.85 &  & \textbf{4708.85} & 310.2 &  & 4754.77 & 282.6 &  & 4708.98 & 0.00 & 4723.87 & 51.10 & 0.32 \\ 
\noalign{\smallskip}\hline\noalign{\smallskip}
R201a  & 100 & 3753.42 &  & 3782.88 & 459.6 &  & 3760.43 & 538.2 &  & 3787.55 & 0.91 & 3838.10 & 265.18 & 2.26 \\ 
R202a  & 100 & 3540.03 &  & \textbf{3540.03} & 802.2 &  & 3554.20 & 598.8 &  & 3540.59 & 0.02 & 3562.54 & 296.59 & 0.64 \\ 
R203a  & 100 & 3311.35 &  & \textbf{3311.35} & 544.2 &  & 3315.50 & 525.6 &  & 3315.16 & 0.12 & 3318.61 & 288.80 & 0.22 \\ 
R204a  & 100 & 3075.95 &  & \textbf{3075.95} & 532.2 &  & \textbf{3075.95} & 478.8 &  & \textbf{3075.95} & 0.00 & 3080.01 & 209.39 & 0.13 \\ 
R205a  & 100 & 3334.27 &  & \textbf{3334.27} & 555.0 &  & \textbf{3334.27} & 507.0 &  & \textbf{3334.27} & 0.00 & 3361.92 & 217.21 & 0.83 \\ 
R206a  & 100 & 3242.40 &  & \textbf{3242.40} & 540.6 &  & 3263.40 & 490.2 &  & \textbf{3242.40} & 0.00 & 3259.19 & 226.03 & 0.52 \\ 
R207a  & 100 & 3145.08 &  & \textbf{3145.08} & 564.0 &  & 3152.29 & 557.4 &  & 3145.83 & 0.02 & 3157.37 & 210.22 & 0.39 \\ 
R208a  & 100 & 3017.12 &  & \textbf{3017.12} & 484.2 &  & \textbf{3017.12} & 510.6 &  & \textbf{3017.12} & 0.00 & 3021.98 & 181.73 & 0.16 \\ 
R209a  & 100 & 3183.36 &  & \textbf{3183.36} & 569.4 &  & 3194.28 & 562.2 &  & 3184.41 & 0.03 & 3190.30 & 174.19 & 0.22 \\ 
R210a  & 100 & 3287.66 &  & \textbf{3287.66} & 612.6 &  & 3309.26 & 527.4 &  & 3288.46 & 0.02 & 3301.07 & 219.59 & 0.41 \\ 
R211a  & 100 & 3019.93 &  & \textbf{3019.93} & 544.8 &  & 3020.56 & 479.4 &  & \textbf{3019.93} & 0.00 & 3022.14 & 155.30 & 0.07 \\ 
\noalign{\smallskip}\hline\noalign{\smallskip}
C201a  & 100 & 5820.78 &  & 5878.54 & 310.2 &  & 5830.20 & 300.0 &  & 5853.90 & 0.57 & 5892.77 & 271.05 & 1.24 \\ 
C202a  & 100 & 5776.88 &  & \textbf{5776.88} & 309.0 &  & \textbf{5776.88} & 310.2 &  & \textbf{5776.88} & 0.00 & \textbf{5776.88} & 242.04 & 0.00 \\ 
C203a  & 100 & 5736.94 &  & 5741.12 & 343.2 &  & 5741.12 & 285.6 &  & 5741.12 & 0.07 & 5742.71 & 196.11 & 0.10 \\ 
C204a  & 100 & 5680.46 &  & \textbf{5680.46} & 258.6 &  & \textbf{5680.46} & 252.6 &  & \textbf{5680.46} & 0.00 & \textbf{5680.46} & 180.73 & 0.00 \\ 
C205a  & 100 & 5747.67 &  & 5781.15 & 393.6 &  & 5751.40 & 407.4 &  & 5781.15 & 0.58 & 5789.29 & 254.67 & 0.72 \\ 
C206a  & 100 & 5738.09 &  & 5767.70 & 284.4 &  & 5741.30 & 258.0 &  & 5767.70 & 0.52 & 5771.59 & 213.11 & 0.58 \\ 
C207a  & 100 & 5721.16 &  & 5731.44 & 308.4 &  & 5725.10 & 250.2 &  & 5731.44 & 0.18 & 5732.97 & 203.22 & 0.21 \\ 
C208a  & 100 & 5725.03 &  & \textbf{5725.03} & 271.2 &  & \textbf{5725.03} & 312.6 &  & \textbf{5725.03} & 0.00 & 5727.22 & 171.01 & 0.04 \\ 
\noalign{\smallskip}\hline\noalign{\smallskip}
RC201a  & 100 & 4701.88 &  & 4737.59 & 316.8 &  & 4707.80 & 270.0 &  & 4737.59 & 0.76 & 4747.24 & 107.41 & 0.96 \\ 
RC202a  & 100 & 4487.48 &  & \textbf{4487.48} & 268.8 &  & 4519.40 & 280.2 &  & \textbf{4487.48} & 0.00 & 4487.57 & 125.19 & 0.00 \\ 
RC203a  & 100 & 4305.49 &  & \textbf{4305.49} & 352.8 &  & 4319.10 & 316.2 &  & \textbf{4305.49} & 0.00 & 4314.27 & 137.30 & 0.20 \\ 
RC204a  & 100 & 4137.93 &  & \textbf{4137.93} & 400.8 &  & 4155.77 & 311.4 &  & \textbf{\textit{4137.84}} & 0.00 & 4142.65 & 124.05 & 0.11 \\ 
RC205a  & 100 & 4585.20 &  & 4615.04 & 384.0 &  & 4595.67 & 413.4 &  & 4615.40 & 0.66 & 4634.25 & 125.16 & 1.07 \\ 
RC206a  & 100 & 4405.16 &  & \textbf{4405.16} & 308.4 &  & 4434.30 & 301.8 &  & \textbf{4405.16} & 0.00 & 4420.53 & 116.34 & 0.35 \\ 
RC207a  & 100 & 4290.14 &  & \textbf{4290.14} & 391.2 &  & 4315.90 & 376.2 &  & \textbf{4290.14} & 0.00 & 4301.27 & 123.83 & 0.26 \\ 
RC208a  & 100 & 4075.04 &  & \textbf{4075.04} & 344.4 &  & 4081.37 & 310.2 &  & \textbf{4075.04} & 0.00 & 4075.96 & 112.59 & 0.02 \\ 
\noalign{\smallskip}\hline\noalign{\smallskip}
\multicolumn{2}{l}{Average}  &  &  &  0.14 & 351.50 &  & 0.36 & 326.51  &  &  & 0.13 &  & 131.13 & 0.38 \\
\hline
\end{tabular}
} 
\label{tbl:Resultados.HH.FSMTW_Al}
\end{table}

\begin{table}[htbp]
\caption{Results for the FSMTW (\small{minimize duration, fleet B})}
\setlength{\tabcolsep}{1.5mm}
\scalebox{0.78} {
\begin{tabular}{lrrlrrlrrlrrrrr}
\hline\noalign{\smallskip} 
 &  &  &  & \multicolumn{2}{l}{VCGP14} &  & \multicolumn{2}{l}{KBJL15} &  & \multicolumn{5}{l}{HILS-RVRP} \\
\cline{5-6} \cline{8-9} \cline{11-15}\noalign{\smallskip}
\multicolumn{1}{c}{Inst.} & \multicolumn{1}{c}{n} & \multicolumn{1}{c}{BKS} &  & \multicolumn{1}{c}{Best Sol.} & \multicolumn{1}{c}{{Avg. T (s)}} &  & \multicolumn{1}{c}{Best Sol.} & {{Avg. T (s)}} & \multicolumn{1}{c}{} & \multicolumn{1}{c}{Best Sol.} & \multicolumn{1}{c}{Gap (\%)} & \multicolumn{1}{c}{{Avg. Sol.}} & \multicolumn{1}{c}{{Avg. T(s)}} & \multicolumn{1}{c}{{Avg. Gap (\%)}} \\
\noalign{\smallskip}\hline\noalign{\smallskip}
R101b  & 100 & 2421.19 &  & 2486.77 & 233.40 &  & 2425.10 & 226.80 &  & 2486.77 & 2.71 & 2494.10 & 54.59 & 3.01 \\ 
R102b  & 100 & 2209.50 &  & 2222.15 & 258.60 & \multicolumn{1}{c}{} & 2212.37 & 238.20 &  & 2222.15 & 0.57 & 2226.88 & 59.11 & 0.79 \\ 
R103b  & 100 & 1930.21 &  & \textbf{1930.21} & 250.80 & \multicolumn{1}{c}{} & 1951.99 & 256.80 &  & 1930.69 & 0.02 & 1941.17 & 65.88 & 0.57 \\ 
R104b  & 100 & 1688.12 &  & \textbf{1688.12} & 260.40 & \multicolumn{1}{c}{} & 1714.86 & 240.60 &  & \textbf{1688.12} & 0.00 & 1710.55 & 65.13 & 1.33 \\ 
R105b  & 100 & 2017.56 &  & \textbf{2017.56} & 229.80 &  & 2024.91 & 220.80 &  & \textbf{2017.56} & 0.00 & 2024.55 & 62.85 & 0.35 \\ 
R106b  & 100 & 1913.84 &  & \textbf{1913.84} & 306.00 &  & 1922.10 & 251.40 &  & \textbf{\textit{1913.04}} & -0.04 & 1918.48 & 70.66 & 0.24 \\ 
R107b  & 100 & 1774.50 &  & \textbf{1774.50} & 256.20 &  & 1783.20 & 318.00 &  & \textbf{1774.50} & 0.00 & 1785.81 & 66.61 & 0.64 \\ 
R108b  & 100 & 1649.24 &  & 1654.68 & 349.80 &  & 1661.58 & 286.80 &  & 1654.78 & 0.34 & 1663.92 & 61.99 & 0.89 \\ 
R109b  & 100 & 1818.15 &  & \textbf{1818.15} & 305.40 &  & 1829.10 & 294.60 &  & \textbf{1818.15} & 0.00 & 1829.18 & 65.23 & 0.61 \\ 
R110b  & 100 & 1761.53 &  & \textbf{1761.53} & 346.20 &  & 1778.80 & 312.60 &  & 1764.64 & 0.18 & 1777.41 & 76.36 & 0.90 \\ 
R111b  & 100 & 1751.10 &  & \textbf{1751.10} & 334.20 &  & 1775.24 & 286.80 &  & 1754.89 & 0.22 & 1773.07 & 77.50 & 1.25 \\ 
R112b  & 100 & 1663.09 &  & \textbf{1663.09} & 379.80 &  & 1677.00 & 372.60 &  & \textbf{\textit{1661.97}} & -0.07 & 1672.61 & 58.08 & 0.57 \\ 
\noalign{\smallskip}\hline\noalign{\smallskip}
C101b  & 100 & 2417.52 &  & \textbf{2417.52} & 123.60 &  & \textbf{2417.52} & 119.40 &  & \textbf{2417.52} & 0.00 & \textbf{2417.52} & 50.16 & 0.00 \\ 
C102b  & 100 & 2350.54 &  & 2350.55 & 178.80 &  & \textbf{2350.54} & 147.00 &  & \textbf{2350.54} & 0.00 & \textbf{2350.54} & 70.56 & 0.00 \\ 
C103b  & 100 & 2345.31 &  & \textbf{2345.31} & 241.20 &  & \textbf{2345.31} & 208.20 &  & \textbf{2345.31} & 0.00 & \textbf{2345.31} & 72.54 & 0.00 \\ 
C104b  & 100 & 2325.78 &  & 2327.84 & 145.80 &  & 2330.59 & 185.40 &  & 2327.84 & 0.09 & 2327.98 & 77.14 & 0.09 \\ 
C105b  & 100 & 2373.53 &  & \textbf{2373.53} & 201.00 &  & 2376.45 & 183.60 &  & \textbf{2373.53} & 0.00 & 2377.74 & 67.68 & 0.18 \\ 
C106b  & 100 & 2381.14 &  & 2386.03 & 190.20 &  & 2386.43 & 177.00 &  & 2386.03 & 0.21 & 2390.56 & 66.59 & 0.40 \\ 
C107b  & 100 & 2357.52 &  & 2364.21 & 186.00 &  & 2359.00 & 147.00 &  & 2364.21 & 0.28 & 2366.44 & 72.71 & 0.38 \\ 
C108b  & 100 & 2346.38 &  & \textbf{2346.38} & 196.80 &  & 2348.15 & 167.40 &  & \textbf{2346.38} & 0.00 & 2347.66 & 68.35 & 0.05 \\ 
C109b  & 100 & 2336.29 &  & \textbf{2336.29} & 156.00 &  & 2337.60 & 153.60 &  & \textbf{2336.29} & 0.00 & \textbf{2336.29} & 90.44 & 0.00 \\ 
\noalign{\smallskip}\hline\noalign{\smallskip}
RC101b  & 100 & 2456.10 &  & \textbf{2456.10} & 281.40 &  & 2464.19 & 268.20 &  & \textbf{2456.10} & 0.00 & 2463.85 & 61.82 & 0.32 \\ 
RC102b  & 100 & 2259.25 &  & \textbf{2259.25} & 254.40 &  & 2270.43 & 247.20 &  & \textbf{2259.25} & 0.00 & 2267.13 & 61.31 & 0.35 \\ 
RC103b  & 100 & 2025.30 &  & \textbf{2025.30} & 284.40 &  & 2041.20 & 238.80 &  & \textbf{2025.30} & 0.00 & 2033.25 & 63.21 & 0.39 \\ 
RC104b  & 100 & 1901.04 &  & \textbf{1901.04} & 262.20 &  & 1922.27 & 252.60 &  & \textbf{1901.04} & 0.00 & 1923.40 & 61.79 & 1.18 \\ 
RC105b  & 100 & 2308.59 &  & 2329.14 & 285.60 &  & 2327.70 & 273.60 &  & 2329.14 & 0.89 & 2332.74 & 60.53 & 1.05 \\ 
RC106b  & 100 & 2146.00 &  & \textbf{2146.00} & 220.80 &  & 2147.14 & 252.60 &  & \textbf{\textit{2140.94}} & -0.24 & 2151.26 & 64.45 & 0.25 \\ 
RC107b  & 100 & 1989.34 &  & \textbf{1989.34} & 245.40 &  & 1996.09 & 251.40 &  & \textbf{1989.34} & 0.00 & 1998.37 & 60.63 & 0.45 \\ 
RC108b  & 100 & 1898.96 &  & \textbf{1898.96} & 195.60 &  & 1908.89 & 186.60 &  & \textbf{1898.96} & 0.00 & 1908.89 & 54.64 & 0.52 \\ 
\noalign{\smallskip}\hline\noalign{\smallskip}
R201b  & 100 & 1953.42 &  & 1973.43 & 383.40 &  & 1956.21 & 372.60 &  & 1980.60 & 1.39 & 1989.18 & 217.75 & 1.83 \\ 
R202b  & 100 & 1740.03 &  & \textbf{1740.03} & 483.00 &  & 1752.40 & 480.00 &  & 1743.82 & 0.22 & 1775.69 & 258.74 & 2.05 \\ 
R203b  & 100 & 1511.35 &  & \textbf{1511.35} & 386.40 &  & 1515.17 & 346.80 &  & 1515.08 & 0.25 & 1518.14 & 254.35 & 0.45 \\ 
R204b  & 100 & 1275.95 &  & \textbf{1275.95} & 454.80 &  & 1279.57 & 413.40 &  & \textbf{1275.95} & 0.00 & 1279.16 & 190.43 & 0.25 \\ 
R205b  & 100 & 1534.27 &  & \textbf{1534.27} & 387.00 &  & 1549.39 & 389.40 &  & 1542.47 & 0.53 & 1566.91 & 196.58 & 2.13 \\ 
R206b  & 100 & 1441.35 &  & \textbf{1441.35} & 355.20 &  & 1450.37 & 312.60 &  & 1448.03 & 0.46 & 1464.74 & 193.32 & 1.62 \\ 
R207b  & 100 & 1345.08 &  & \textbf{1345.08} & 418.20 &  & 1359.18 & 378.60 &  & 1346.18 & 0.08 & 1362.59 & 191.29 & 1.30 \\ 
R208b  & 100 & 1217.12 &  & \textbf{1217.12} & 360.00 &  & 1220.36 & 328.20 &  & \textbf{1217.12} & 0.00 & 1220.92 & 168.08 & 0.31 \\ 
R209b  & 100 & 1380.79 &  & \textbf{1380.79} & 465.00 &  & 1385.65 & 428.40 &  & 1384.42 & 0.26 & 1391.29 & 172.35 & 0.76 \\ 
R210b  & 100 & 1485.65 &  & \textbf{1485.65} & 463.20 &  & 1495.75 & 415.80 &  & 1492.12 & 0.44 & 1507.17 & 205.78 & 1.45 \\ 
R211b  & 100 & 1219.93 &  & \textbf{1219.93} & 441.60 &  & \textbf{1219.93} & 447.00 &  & \textbf{1219.93} & 0.00 & 1226.01 & 152.18 & 0.50 \\ 
\noalign{\smallskip}\hline\noalign{\smallskip}
C201b  & 100 & 1816.14 &  & 1820.64 & 174.00 &  & 1820.64 & 186.60 &  & 1820.64 & 0.25 & 1820.64 & 194.23 & 0.25 \\ 
C202b  & 100 & 1768.51 &  & \textbf{1768.51} & 313.20 &  & 1770.10 & 274.80 &  & \textbf{1768.51} & 0.00 & 1778.44 & 173.71 & 0.56 \\ 
C203b  & 100 & 1733.63 &  & \textbf{1733.63} & 197.40 &  & \textbf{1733.63} & 191.40 &  & \textbf{1733.63} & 0.00 & \textbf{1733.63} & 159.40 & 0.00 \\ 
C204b  & 100 & 1680.46 &  & \textbf{1680.46} & 198.00 &  & \textbf{1680.46} & 190.20 &  & \textbf{1680.46} & 0.00 & \textbf{1680.46} & 164.82 & 0.00 \\ 
C205b  & 100 & 1747.68 &  & 1778.30 & 328.80 &  & 1756.54 & 312.60 &  & 1778.30 & 1.75 & 1784.42 & 193.49 & 2.10 \\ 
C206b  & 100 & 1756.01 &  & 1767.70 & 230.40 &  & 1773.17 & 207.60 &  & 1767.70 & 0.67 & 1776.92 & 178.08 & 1.19 \\ 
C207b  & 100 & 1729.39 &  & 1729.49 & 208.80 &  & \textbf{1729.39} & 178.20 &  & 1729.49 & 0.01 & 1733.41 & 167.76 & 0.23 \\ 
C208b  & 100 & 1723.20 &  & 1724.20 & 204.00 &  & 1724.20 & 187.80 &  & 1724.20 & 0.06 & 1724.20 & 167.91 & 0.06 \\ 
\noalign{\smallskip}\hline\noalign{\smallskip}
RC201b  & 100 & 2230.54 &  & 2329.59 & 260.40 &  & 2235.90 & 250.20 &  & 2331.14 & 4.51 & 2345.62 & 124.50 & 5.16 \\ 
RC202b  & 100 & 2002.62 &  & 2057.66 & 401.40 &  & 2022.00 & 328.20 &  & 2057.66 & 2.75 & 2078.10 & 144.73 & 3.77 \\ 
RC203b  & 100 & 1824.54 &  & \textbf{1824.54} & 319.80 &  & 1840.40 & 307.20 &  & 1824.88 & 0.02 & 1847.01 & 165.37 & 1.23 \\ 
RC204b  & 100 & 1555.74 &  & 1555.75 & 330.00 &  & \textbf{1555.74} & 298.80 &  & 1555.75 & 0.00 & 1560.82 & 146.18 & 0.33 \\ 
RC205b  & 100 & 2166.62 &  & 2174.74 & 325.80 &  & 2169.00 & 388.20 &  & 2174.74 & 0.37 & 2196.93 & 144.16 & 1.40 \\ 
RC206b  & 100 & 1883.08 &  & \textbf{1883.08} & 259.80 &  & 1898.70 & 248.40 &  & \textbf{1883.08} & 0.00 & 1892.52 & 137.24 & 0.50 \\ 
RC207b  & 100 & 1714.14 &  & \textbf{1714.14} & 339.00 &  & 1730.00 & 308.40 &  & \textbf{1714.14} & 0.00 & 1733.17 & 145.07 & 1.11 \\ 
RC208b  & 100 & 1483.20 &  & \textbf{1483.20} & 288.00 &  & 1490.64 & 265.80 &  & \textbf{1483.20} & 0.00 & 1487.17 & 141.23 & 0.27 \\ 
\noalign{\smallskip}\hline\noalign{\smallskip}
 \multicolumn{2}{l}{Average}  &  &  & 0.29  & 288.15 &  & 0.46 & 271.48 &  &  & 0.34 &  & 118.83 & 0.85 \\
 \hline
\end{tabular}
}
\label{tbl:Resultados.HH.FSMTW_Bl}
\end{table}

\begin{table}[htbp]
\caption{Results for the FSMTW (\small{minimize duration, fleet C})}
\setlength{\tabcolsep}{1.5mm}
\scalebox{0.78} {
\begin{tabular}{lrrlrrlrrlrrrrr}
\hline\noalign{\smallskip} 
 &  &  &  & \multicolumn{2}{l}{VCGP14} &  & \multicolumn{2}{l}{KBJL15} &  & \multicolumn{5}{l}{HILS-RVRP} \\
\cline{5-6} \cline{8-9} \cline{11-15} \noalign{\smallskip}
\multicolumn{1}{c}{Inst.} & \multicolumn{1}{c}{n} & \multicolumn{1}{c}{BKS} &  & \multicolumn{1}{c}{Best Sol.} & \multicolumn{1}{c}{{Avg. T(s)}} &  & \multicolumn{1}{c}{Best Sol.} & {{Avg. T(s)}} & \multicolumn{1}{c}{} & \multicolumn{1}{c}{Best Sol.} & \multicolumn{1}{c}{Gap (\%)} & \multicolumn{1}{c}{{Avg. Sol.}} & \multicolumn{1}{c}{{Avg. T(s)}} & \multicolumn{1}{c}{{Avg. Gap (\%)}} \\
\noalign{\smallskip}\hline\noalign{\smallskip}
R101c  & 100 & 2134.90 &  & 2199.79 & 202.80 &  & 2137.20 & 188.40 &  & 2199.79 & 3.04 & 2200.87 & 54.19 & 3.09 \\ 
R102c  & 100 & 1913.37 &  & 1925.56 & 303.60 & \multicolumn{1}{c}{} & 1914.87 & 372.60 &  & 1925.56 & 0.64 & 1927.82 & 62.85 & 0.76 \\ 
R103c  & 100 & 1609.94 &  & 1615.38 & 217.20 & \multicolumn{1}{c}{} & 1621.20 & 194.40 &  & 1615.38 & 0.34 & 1619.49 & 67.40 & 0.59 \\ 
R104c  & 100 & 1363.26 &  & \textbf{1363.26} & 274.80 & \multicolumn{1}{c}{} & 1375.60 & 268.20 &  & \textbf{1363.26} & 0.00 & 1373.84 & 66.13 & 0.78 \\ 
R105c  & 100 & 1722.05 &  & \textbf{1722.05} & 216.00 &  & \textbf{1722.05} & 190.20 &  & \textbf{1722.05} & 0.00 & 1724.62 & 62.42 & 0.15 \\ 
R106c  & 100 & 1599.04 &  & \textbf{1599.04} & 286.20 &  & 1610.40 & 244.80 &  & \textbf{1599.04} & 0.00 & 1605.26 & 69.50 & 0.39 \\ 
R107c  & 100 & 1442.97 &  & \textbf{1442.97} & 223.20 &  & 1454.30 & 210.60 &  & \textbf{1442.97} & 0.00 & 1458.92 & 70.03 & 1.11 \\ 
R108c  & 100 & 1321.68 &  & \textbf{1321.68} & 326.40 &  & 1329.92 & 319.80 &  & \textbf{\textit{1317.43}} & -0.32 & 1333.87 & 62.38 & 0.92 \\ 
R109c  & 100 & 1506.59 &  & \textbf{1506.59} & 301.20 &  & 1507.10 & 283.80 &  & \textbf{1506.59} & 0.00 & 1511.26 & 68.29 & 0.31 \\ 
R110c  & 100 & 1443.92 &  & \textbf{1443.92} & 343.80 &  & 1451.06 & 327.60 &  & \textbf{\textit{1443.33}} & -0.04 & 1457.89 & 76.92 & 0.97 \\ 
R111c  & 100 & 1423.47 &  & \textbf{1423.47} & 419.40 &  & 1436.32 & 368.40 &  & 1425.18 & 0.12 & 1444.20 & 71.06 & 1.46 \\ 
R112c  & 100 & 1329.07 &  & \textbf{1329.07} & 286.20 &  & 1341.10 & 250.20 &  & \textbf{1329.07} & 0.00 & 1346.41 & 65.15 & 1.30 \\ 
\noalign{\smallskip}\hline\noalign{\smallskip}
C101c  & 100 & 1628.31 &  & 1628.94 & 109.80 &  & 1628.94 & 118.20 &  & 1628.94 & 0.04 & 1628.94 & 49.98 & 0.04 \\ 
C102c  & 100 & 1610.96 &  & \textbf{1610.96} & 145.80 &  & \textbf{1610.96} & 151.80 &  & \textbf{1610.96} & 0.00 & \textbf{1610.96} & 63.04 & 0.00 \\ 
C103c  & 100 & 1607.14 &  & \textbf{1607.14} & 166.20 &  & \textbf{1607.14} & 227.40 &  & \textbf{1607.14} & 0.00 & \textbf{1607.14} & 61.19 & 0.00 \\ 
C104c  & 100 & 1598.50 &  & 1599.90 & 162.00 &  & 1599.21 & 173.40 &  & 1599.90 & 0.09 & 1599.90 & 62.97 & 0.09 \\ 
C105c  & 100 & 1628.38 &  & 1628.94 & 112.80 &  & 1628.94 & 118.20 &  & 1628.94 & 0.03 & 1628.94 & 61.15 & 0.03 \\ 
C106c  & 100 & 1628.94 &  & \textbf{1628.94} & 112.80 &  & \textbf{1628.94} & 120.60 &  & \textbf{1628.94} & 0.00 & \textbf{1628.94} & 56.94 & 0.00 \\ 
C107c  & 100 & 1628.38 &  & 1628.94 & 117.60 &  & 1628.94 & 119.40 &  & 1628.94 & 0.03 & 1628.94 & 59.20 & 0.03 \\ 
C108c  & 100 & 1622.89 &  & \textbf{1622.89} & 180.00 &  & 1625.00 & 147.00 &  & \textbf{1622.89} & 0.00 & 1627.13 & 58.86 & 0.26 \\ 
C109c  & 100 & 1614.99 &  & 1615.93 & 217.20 &  & 1618.61 & 212.40 &  & \textbf{1614.99} & 0.00 & 1615.95 & 64.59 & 0.06 \\ 
\noalign{\smallskip}\hline\noalign{\smallskip}
RC101c  & 100 & 2082.95 &  & \textbf{2082.95} & 294.60 &  & 2092.10 & 272.40 &  & \textbf{2082.95} & 0.00 & 2086.57 & 60.68 & 0.17 \\ 
RC102c  & 100 & 1895.05 &  & \textbf{1895.05} & 256.80 &  & 1901.89 & 251.40 &  & 1895.56 & 0.03 & 1901.70 & 63.18 & 0.35 \\ 
RC103c  & 100 & 1650.30 &  & \textbf{1650.30} & 238.80 &  & 1660.70 & 213.60 &  & \textbf{1650.30} & 0.00 & 1666.02 & 62.18 & 0.95 \\ 
RC104c  & 100 & 1526.04 &  & \textbf{1526.04} & 214.20 &  & 1540.60 & 208.20 &  & \textbf{1526.04} & 0.00 & 1551.14 & 62.96 & 1.64 \\ 
RC105c  & 100 & 1953.99 &  & 1957.14 & 282.60 &  & 1956.09 & 249.60 &  & 1957.14 & 0.16 & 1958.54 & 59.19 & 0.23 \\ 
RC106c  & 100 & 1774.94 &  & \textbf{1774.94} & 228.00 &  & 1780.45 & 209.40 &  & \textbf{1774.94} & 0.00 & 1778.88 & 64.92 & 0.22 \\ 
RC107c  & 100 & 1607.11 &  & \textbf{1607.11} & 229.80 &  & 1620.30 & 184.20 &  & \textbf{1607.11} & 0.00 & 1620.06 & 59.59 & 0.81 \\ 
RC108c  & 100 & 1523.96 &  & \textbf{1523.96} & 202.80 &  & 1532.60 & 213.60 &  & \textbf{1523.96} & 0.00 & 1526.00 & 57.50 & 0.13 \\ 
\noalign{\smallskip}\hline\noalign{\smallskip}
R201c  & 100 & 1716.02 &  & \textbf{1716.02} & 272.40 &  & 1731.20 & 406.80 &  & \textbf{1716.02} & 0.00 & 1720.65 & 200.64 & 0.27 \\ 
R202c  & 100 & 1515.03 &  & \textbf{1515.03} & 530.40 &  & 1529.70 & 488.40 &  & 1519.41 & 0.29 & 1536.94 & 234.74 & 1.45 \\ 
R203c  & 100 & 1286.35 &  & \textbf{1286.35} & 374.40 &  & 1296.72 & 390.00 &  & 1290.16 & 0.30 & 1297.61 & 266.88 & 0.87 \\ 
R204c  & 100 & 1050.95 &  & \textbf{1050.95} & 457.20 &  & 1052.90 & 473.40 &  & \textbf{1050.95} & 0.00 & 1057.00 & 192.24 & 0.58 \\ 
R205c  & 100 & 1309.27 &  & \textbf{1309.27} & 386.40 &  & 1315.20 & 402.60 &  & \textbf{1309.27} & 0.00 & 1317.49 & 183.91 & 0.63 \\ 
R206c  & 100 & 1216.35 &  & \textbf{1216.35} & 320.40 &  & 1226.93 & 395.40 &  & 1223.07 & 0.55 & 1234.31 & 180.06 & 1.48 \\ 
R207c  & 100 & 1120.08 &  & \textbf{1120.08} & 433.80 &  & 1125.50 & 418.80 &  & 1121.80 & 0.15 & 1137.13 & 201.30 & 1.52 \\ 
R208c  & 100 & 992.12 &  & \textbf{992.12} & 360.60 &  & 997.97 & 352.20 &  & \textbf{992.12} & 0.00 & 1000.17 & 182.55 & 0.81 \\ 
R209c  & 100 & 1155.79 &  & \textbf{1155.79} & 450.00 &  & 1164.31 & 428.40 &  & 1159.67 & 0.34 & 1168.88 & 170.29 & 1.13 \\ 
R210c  & 100 & 1257.89 &  & \textbf{1257.89} & 392.40 &  & 1269.70 & 368.40 &  & 1264.38 & 0.52 & 1272.87 & 201.36 & 1.19 \\ 
R211c  & 100 & 994.93 &  & \textbf{994.93} & 395.40 &  & 995.58 & 370.20 &  & \textbf{994.93} & 0.00 & 1000.99 & 156.78 & 0.61 \\ 
\noalign{\smallskip}\hline\noalign{\smallskip}
C201c  & 100 & 1250.97 &  & 1269.41 & 171.60 &  & \textbf{1250.97} & 178.20 &  & 1269.41 & 1.47 & 1277.27 & 176.97 & 2.10 \\ 
C202c  & 100 & 1239.54 &  & \textbf{1239.54} & 231.00 &  & 1240.86 & 212.40 &  & \textbf{1239.54} & 0.00 & 1242.29 & 165.61 & 0.22 \\ 
C203c  & 100 & 1193.63 &  & \textbf{1193.63} & 181.80 &  & \textbf{1193.63} & 188.40 &  & \textbf{1193.63} & 0.00 & \textbf{1193.63} & 148.25 & 0.00 \\ 
C204c  & 100 & 1176.52 &  & \textbf{1176.52} & 234.00 &  & \textbf{1176.52} & 220.20 &  & \textbf{1176.52} & 0.00 & \textbf{1176.52} & 161.86 & 0.00 \\ 
C205c  & 100 & 1238.30 &  & \textbf{1238.30} & 261.60 &  & 1240.10 & 257.40 &  & \textbf{1238.30} & 0.00 & 1245.38 & 180.27 & 0.57 \\ 
C206c  & 100 & 1229.23 &  & 1238.30 & 292.20 &  & \textbf{1229.23} & 262.80 &  & 1238.30 & 0.74 & 1241.36 & 163.15 & 0.99 \\ 
C207c  & 100 & 1209.48 &  & 1209.49 & 180.00 &  & \textbf{1209.48} & 213.60 &  & 1209.49 & 0.00 & 1214.57 & 158.96 & 0.42 \\ 
C208c  & 100 & 1204.20 &  & \textbf{1204.20} & 181.80 &  & \textbf{1204.20} & 180.60 &  & \textbf{1204.20} & 0.00 & \textbf{1204.20} & 158.52 & 0.00 \\ 
\noalign{\smallskip}\hline\noalign{\smallskip}
RC201c  & 100 & 1915.42 &  & 1996.79 & 220.20 &  & 1917.90 & 279.00 &  & 1996.79 & 4.25 & 2010.99 & 127.73 & 4.99 \\ 
RC202c  & 100 & 1677.62 &  & 1732.66 & 391.80 &  & 1680.00 & 366.00 &  & 1732.66 & 3.28 & 1753.28 & 150.22 & 4.51 \\ 
RC203c  & 100 & 1496.11 &  & \textbf{1496.11} & 369.00 &  & 1500.20 & 376.20 &  & 1499.54 & 0.23 & 1515.94 & 168.38 & 1.33 \\ 
RC204c  & 100 & 1220.75 &  & \textbf{1220.75} & 327.00 &  & 1222.16 & 328.20 &  & \textbf{1220.75} & 0.00 & 1229.53 & 146.27 & 0.72 \\ 
RC205c  & 100 & 1822.07 &  & 1844.74 & 304.20 &  & 1823.00 & 317.40 &  & 1844.74 & 1.24 & 1866.61 & 161.43 & 2.44 \\ 
RC206c  & 100 & 1553.65 &  & \textbf{1553.65} & 267.00 &  & 1564.30 & 282.00 &  & 1558.08 & 0.29 & 1572.54 & 138.66 & 1.22 \\ 
RC207c  & 100 & 1377.52 &  & \textbf{1377.52} & 363.60 &  & 1381.71 & 340.20 &  & 1379.14 & 0.12 & 1393.87 & 153.51 & 1.19 \\ 
RC208c  & 100 & 1140.10 &  & \textbf{1140.10} & 379.20 &  & 1151.40 & 310.20 &  & \textbf{1140.10} & 0.00 & 1148.27 & 143.06 & 0.72 \\ 
\noalign{\smallskip}\hline\noalign{\smallskip}
\multicolumn{2}{l}{Average}  &  &  & 0.28 & 275.04 &  & 0.37 & 271.74 &  &  & 0.32 &  & 116.04 & 0.87 \\
\hline
\end{tabular}
}
\label{tbl:Resultados.HH.FSMTW_Cl}
\end{table}

\newpage
\subsection{FSMTW -- {{minimizing}} total distance} \label{sec:Result.HH.FSMTW-distance}



Detailed results obtained for the FSMVRPTW instances of \cite[LS99]{LiuShen1999}, considering the objective of minimizing the sum of the total distance, compared with those found by the UHGS of \cite[VCGP14]{Vidaletal2014} and the HEA of \cite[KBJL15]{Kocetal2015} (Tables \ref{tbl:Resultados.HH.FSMTW_An}-\ref{tbl:Resultados.HH.FSMTW_Cn}).


\begin{table}[htbp]
\caption{Results for the FSMTW (\small{minimize distance, fleet A})}
\setlength{\tabcolsep}{1.5mm}
\scalebox{0.78} {
\begin{tabular}{lrrlrrlrrlrrrrr}
\hline\noalign{\smallskip}
 &  &  &  & \multicolumn{2}{l}{VCGP14} &  & \multicolumn{2}{l}{KBJL15} &  & \multicolumn{5}{l}{HILS-RVRP} \\
\cline{5-6} \cline{8-9} \cline{11-15}\noalign{\smallskip}
\multicolumn{1}{c}{Inst.} & \multicolumn{1}{c}{n} & \multicolumn{1}{c}{BKS} &  & \multicolumn{1}{c}{Best Sol.} & \multicolumn{1}{c}{{Avg. T (s)}} &  & \multicolumn{1}{c}{Best Sol.} & {{Avg. T (s)}} & \multicolumn{1}{c}{} & \multicolumn{1}{c}{Best Sol.} & \multicolumn{1}{c}{Gap (\%)} & \multicolumn{1}{c}{{Avg. Sol.}} & \multicolumn{1}{c}{{Avg. T(s)}} & \multicolumn{1}{c}{{Avg. Gap (\%)}} \\
\noalign{\smallskip}\hline\noalign{\smallskip}
R101 & 100 & {4314.36} &  & \textbf{4314.36} & 276.60 &  & 4317.52 & 248.40 &  & \textbf{4314.36} & 0.00 & 4325.76 & 120.56 & 0.26 \\ 
R102 & 100 & 4166.28 &  & \textbf{4166.28} & 361.80 & \multicolumn{1}{c}{} & 4173.84 & 358.80 &  & \textbf{4166.28} & 0.00 & 4181.25 & 124.45 & 0.36 \\ 
R103 & 100 & 4027.36 &  & \textbf{4027.36} & 321.00 & \multicolumn{1}{c}{} & 4031.40 & 312.60 &  & \textbf{\textit{4024.14}} & -0.08 & 4038.85 & 134.79 & 0.29 \\ 
R104 & 100 & 3936.40 &  & \textbf{3936.40} & 288.60 & \multicolumn{1}{c}{} & 3946.44 & 247.20 &  & 3936.55 & 0.00 & 3948.00 & 113.83 & 0.29 \\ 
R105 & 100 & 4122.50 &  & \textbf{4122.50} & 389.40 &  & 4134.06 & 360.60 &  & \textbf{4122.50} & 0.00 & 4133.06 & 131.61 & 0.26 \\ 
R106 & 100 & 4048.59 &  & \textbf{4048.59} & 334.20 &  & 4060.05 & 307.20 &  & 4050.17 & 0.04 & 4059.07 & 127.47 & 0.26 \\ 
R107 & 100 & 3970.51 &  & \textbf{3970.51} & 333.60 &  & 3985.12 & 286.80 &  & 3976.40 & 0.15 & 3988.73 & 126.16 & 0.46 \\ 
R108 & 100 & 3928.12 &  & \textbf{3928.12} & 280.80 &  & 3932.60 & 392.40 &  & \textbf{3928.12} & 0.00 & 3935.45 & 108.28 & 0.19 \\ 
R109 & 100 & 4015.71 &  & \textbf{4015.71} & 288.00 &  & 4024.83 & 367.20 &  & \textbf{4015.71} & 0.00 & 4023.34 & 128.21 & 0.19 \\ 
R110 & 100 & 3961.68 &  & \textbf{3961.68} & 389.40 &  & 3973.51 & 312.60 &  & \textbf{3961.68} & 0.00 & 3965.89 & 114.65 & 0.11 \\ 
R111 & 100 & 3964.99 &  & \textbf{3964.99} & 316.80 &  & 3988.00 & 307.20 &  & 3971.90 & 0.17 & 3989.84 & 127.58 & 0.63 \\ 
R112 & 100 & 3918.88 &  & \textbf{3918.88} & 295.20 &  & 3930.19 & 282.60 &  & \textbf{\textit{3917.88}} & -0.03 & 3927.20 & 105.96 & 0.21 \\ 
\noalign{\smallskip}\hline\noalign{\smallskip}
C101 & 100 & 7093.45 &  & \textbf{7093.45} & 177.60 &  & \textbf{7093.45} & 148.20 &  & \textbf{7093.45} & 0.00 & 7093.59 & 167.56 & 0.00 \\ 
C102 & 100 & 7080.17 &  & \textbf{7080.17} & 128.40 &  & \textbf{7080.17} & 159.00 &  & \textbf{7080.17} & 0.00 & 7080.17 & 134.68 & 0.00 \\ 
C103 & 100 & 7079.21 &  & \textbf{7079.21} & 125.40 &  & \textbf{7079.21} & 120.60 &  & \textbf{7079.21} & 0.00 & 7079.21 & 119.00 & 0.00 \\ 
C104 & 100 & 7075.06 &  & \textbf{7075.06} & 131.40 &  & \textbf{7075.06} & 118.20 &  & \textbf{7075.06} & 0.00 & 7075.06 & 101.36 & 0.00 \\ 
C105 & 100 & 7093.45 &  & \textbf{7093.45} & 199.80 &  & \textbf{7093.45} & 159.00 &  & \textbf{7093.45} & 0.00 & 7093.60 & 154.47 & 0.00 \\ 
C106 & 100 & 7083.87 &  & \textbf{7083.87} & 136.80 &  & \textbf{7083.87} & 130.20 &  & \textbf{7083.87} & 0.00 & 7083.87 & 140.25 & 0.00 \\ 
C107 & 100 & 7084.61 &  & \textbf{7084.61} & 133.80 &  & \textbf{7084.61} & 143.40 &  & \textbf{7084.61} & 0.00 & 7084.61 & 142.94 & 0.00 \\ 
C108 & 100 & 7079.66 &  & \textbf{7079.66} & 131.40 &  & \textbf{7079.66} & 118.20 &  & \textbf{7079.66} & 0.00 & 7079.66 & 124.35 & 0.00 \\ 
C109 & 100 & {7077.30} &  & \textbf{7077.30} & 122.40 &  & \textbf{7077.30} & 131.40 &  & \textbf{7077.30} & 0.00 & 7077.30 & 107.98 & 0.00 \\ 
\noalign{\smallskip}\hline\noalign{\smallskip}
RC101 & 100 & 5150.86 &  & \textbf{5150.86} & 312.60 &  & 5173.47 & 308.40 &  & \textbf{5150.86} & 0.00 & 5160.03 & 121.21 & 0.18 \\ 
RC102 & 100 & 4987.24 &  & \textbf{4987.24} & 288.60 &  & 5018.83 & 255.60 &  & \textbf{\textit{4974.82}} & -0.25 & 4999.64 & 122.10 & 0.25 \\ 
RC103 & 100 & 4804.61 &  & \textbf{4804.61} & 424.80 &  & 4850.20 & 388.20 &  & \textbf{4804.61} & 0.00 & 4837.40 & 122.75 & 0.68 \\ 
RC104 & 100 & 4717.63 &  & \textbf{4717.63} & 318.00 &  & 4725.40 & 317.40 &  & 4721.44 & 0.08 & 4734.77 & 94.74 & 0.36 \\ 
RC105 & 100 & 5035.35 &  & \textbf{5035.35} & 334.20 &  & 5048.86 & 286.80 &  & 5036.50 & 0.02 & 5047.72 & 119.73 & 0.25 \\ 
RC106 & 100 & 4936.74 &  & \textbf{4936.74} & 337.80 &  & 4964.13 & 317.40 &  & \textbf{\textit{4921.13}} & -0.32 & 4941.27 & 118.48 & 0.09 \\ 
RC107 & 100 & 4788.69 &  & \textbf{4788.69} & 304.80 &  & 4825.60 & 250.20 &  & \textbf{\textit{4787.59}} & -0.02 & 4807.65 & 107.49 & 0.40 \\ 
RC108 & 100 & 4708.85 &  & \textbf{4708.85} & 286.80 &  & 4724.79 & 277.80 &  & 4711.31 & 0.05 & 4726.02 & 86.50 & 0.36 \\ 
\noalign{\smallskip}\hline\noalign{\smallskip}
R201 & 100 & 3446.78 &  & \textbf{3446.78} & 390.60 &  & \textbf{3446.78} & 367.80 &  & \textbf{3446.78} & 0.00 & 3452.08 & 320.01 & 0.15 \\ 
R202 & 100 & 3297.42 &  & 3308.16 & 460.80 &  & \textbf{3297.42} & 447.60 &  & 3308.16 & 0.33 & 3313.28 & 342.72 & 0.48 \\ 
R203 & 100 & 3141.09 &  & \textbf{3141.09} & 339.00 &  & \textbf{3141.09} & 368.40 &  & \textbf{3141.09} & 0.00 & 3143.21 & 321.32 & 0.07 \\ 
R204 & 100 & 3018.14 &  & \textbf{3018.14} & 417.60 &  & \textbf{3018.14} & 376.80 &  & \textbf{3018.14} & 0.00 & 3019.95 & 296.69 & 0.06 \\ 
R205 & 100 & 3218.97 &  & \textbf{3218.97} & 384.00 &  & \textbf{3218.97} & 382.80 &  & \textbf{3218.97} & 0.00 & 3228.75 & 305.08 & 0.30 \\ 
R206 & 100 & 3146.34 &  & \textbf{3146.34} & 618.00 &  & \textbf{3146.34} & 488.40 &  & 3147.41 & 0.03 & 3155.39 & 325.90 & 0.29 \\ 
R207 & 100 & 3077.36 &  & 3077.58 & 522.00 &  & \textbf{3077.36} & 388.20 &  & 3077.58 & 0.01 & 3082.82 & 313.42 & 0.18 \\ 
R208 & 100 & 2997.24 &  & \textbf{2997.24} & 322.20 &  & 2997.25 & 380.40 &  & \textbf{2997.24} & 0.00 & 2999.97 & 282.11 & 0.09 \\ 
R209 & 100 & 3119.56 &  & 3122.42 & 382.20 &  & \textbf{3119.56} & 299.40 &  & 3122.42 & 0.09 & 3129.93 & 301.95 & 0.33 \\ 
R210 & 100 & 3170.41 &  & 3174.85 & 415.80 &  & \textbf{3170.41} & 328.20 &  & 3174.31 & 0.12 & 3181.35 & 331.93 & 0.35 \\ 
R211 & 100 & 3019.93 &  & \textbf{3019.93} & 546.00 &  & \textbf{3019.93} & 475.80 &  & \textbf{3019.93} & 0.00 & 3022.78 & 285.35 & 0.09 \\ 
\noalign{\smallskip}\hline\noalign{\smallskip}
C201 & 100 & 5695.02 &  & \textbf{5695.02} & 222.60 &  & \textbf{5695.02} & 207.60 &  & \textbf{5695.02} & 0.00 & 5695.02 & 345.40 & 0.00 \\ 
C202 & 100 & 5685.24 &  & \textbf{5685.24} & 226.80 &  & \textbf{5685.24} & 190.20 &  & \textbf{5685.24} & 0.00 & 5685.24 & 287.90 & 0.00 \\ 
C203 & 100 & 5681.55 &  & \textbf{5681.55} & 252.60 &  & \textbf{5681.55} & 257.40 &  & \textbf{5681.55} & 0.00 & 5681.79 & 271.20 & 0.00 \\ 
C204 & 100 & 5677.66 &  & \textbf{5677.66} & 256.20 &  & \textbf{5677.66} & 238.20 &  & \textbf{5677.66} & 0.00 & 5677.88 & 281.35 & 0.00 \\ 
C205 & 100 & 5691.36 &  & \textbf{5691.36} & 238.80 &  & \textbf{5691.36} & 207.60 &  & \textbf{5691.36} & 0.00 & 5691.36 & 300.58 & 0.00 \\ 
C206 & 100 & 5689.32 &  & \textbf{5689.32} & 229.20 &  & \textbf{5689.32} & 178.20 &  & \textbf{5689.32} & 0.00 & 5689.32 & 270.71 & 0.00 \\ 
C207 & 100 & 5687.35 &  & \textbf{5687.35} & 254.40 &  & \textbf{5687.35} & 246.00 &  & \textbf{5687.35} & 0.00 & 5687.35 & 294.10 & 0.00 \\ 
C208 & 100 & {5686.50} &  & \textbf{5686.50} & 231.60 &  & \textbf{5686.50} & 213.60 &  & \textbf{5686.50} & 0.00 & 5686.50 & 265.14 & 0.00 \\ 
\noalign{\smallskip}\hline\noalign{\smallskip}
RC201 & 100 & 4374.09 &  & \textbf{4374.09} & 355.20 &  & 4376.82 & 308.40 &  & \textbf{4374.09} & 0.00 & 4380.07 & 193.35 & 0.14 \\ 
RC202 & 100 & 4244.63 &  & \textbf{4244.63} & 277.80 &  & \textbf{4244.63} & 255.60 &  & \textbf{4244.63} & 0.00 & 4246.90 & 203.69 & 0.05 \\ 
RC203 & 100 & 4170.17 &  & \textbf{4170.17} & 463.80 &  & \textbf{4170.17} & 368.40 &  & \textbf{4170.17} & 0.00 & 4177.62 & 196.96 & 0.18 \\ 
RC204 & 100 & 4087.11 &  & \textbf{4087.11} & 347.40 &  & \textbf{4087.11} & 328.20 &  & \textbf{4087.11} & 0.00 & 4094.04 & 171.97 & 0.17 \\ 
RC205 & 100 & 4291.93 &  & \textbf{4291.93} & 327.60 &  & 4293.73 & 251.40 &  & \textbf{4291.93} & 0.00 & 4294.56 & 197.21 & 0.06 \\ 
RC206 & 100 & 4251.88 &  & \textbf{4251.88} & 307.20 &  & \textbf{4251.88} & 256.20 &  & \textbf{4251.88} & 0.00 & 4257.78 & 207.92 & 0.14 \\ 
RC207 & 100 & 4182.44 &  & 4185.98 & 289.20 &  & \textbf{4182.44} & 338.40 &  & 4185.98 & 0.08 & 4188.14 & 192.60 & 0.14 \\ 
RC208 & 100 & 4075.04 &  & \textbf{4075.04} & 244.80 &  & \textbf{4075.04} & 318.60 &  & \textbf{4075.04} & 0.00 & 4077.57 & 166.98 & 0.06 \\ 
\noalign{\smallskip}\hline\noalign{\smallskip}
\multicolumn{2}{l}{Average}  &  &  & 0.01 & 305.24 &  & 0.13 & 283.60 &  &  & 0.01 &  & 193.26 & 0.17 \\
\hline
\end{tabular}
}
\label{tbl:Resultados.HH.FSMTW_An}
\end{table}

\begin{table}[htbp]
\caption{Results for the FSMTW (\small{minimize distance, fleet B})}
\setlength{\tabcolsep}{1.5mm}
\scalebox{0.78} {
\begin{tabular}{lrrlrrlrrlrrrrr}
\hline\noalign{\smallskip}
 &  &  &  & \multicolumn{2}{l}{VCGP14} &  & \multicolumn{2}{l}{KBJL15} &  & \multicolumn{5}{l}{HILS-RVRP} \\
\cline{5-6} \cline{8-9} \cline{11-15}\noalign{\smallskip}
\multicolumn{1}{c}{Inst.} & \multicolumn{1}{c}{n} & \multicolumn{1}{c}{BKS} &  & \multicolumn{1}{c}{Best Sol.} & \multicolumn{1}{c}{{Avg. T (s)}} &  & \multicolumn{1}{c}{Best Sol.} & {{Avg. T (s)}} & \multicolumn{1}{c}{} & \multicolumn{1}{c}{Best Sol.} & \multicolumn{1}{c}{Gap (\%)} & \multicolumn{1}{c}{{Avg. Sol.}} & \multicolumn{1}{c}{{Avg. T(s)}} & \multicolumn{1}{c}{{Avg. Gap (\%)}} \\
\noalign{\smallskip}\hline\noalign{\smallskip}
R101 & 100 & 2222.56 &  & 2228.67 & 303.00 &  & \textbf{2222.56} & 256.20 &  & 2228.67 & 0.27 & 2229.99 & 67.06 & 0.33 \\ 
R102 & 100 & 2048.12 &  & 2073.63 & 215.40 & \multicolumn{1}{c}{} & \textbf{2048.12} & 196.80 &  & 2071.90 & 1.16 & 2073.92 & 71.71 & 1.26 \\ 
R103 & 100 & 1853.66 &  & \textbf{1853.66} & 274.20 & \multicolumn{1}{c}{} & 1855.74 & 316.20 &  & \textbf{1853.66} & 0.00 & 1857.83 & 73.56 & 0.23 \\ 
R104 & 100 & 1683.33 &  & \textbf{1683.33} & 322.20 & \multicolumn{1}{c}{} & 1686.42 & 305.40 &  & 1685.49 & 0.13 & 1691.53 & 65.13 & 0.49 \\ 
R105 & 100 & 1980.96 &  & 1988.86 & 198.00 &  & \textbf{1980.96} & 202.20 &  & 1988.86 & 0.40 & 1992.22 & 70.02 & 0.57 \\ 
R106 & 100 & 1888.31 &  & \textbf{1888.31} & 278.40 &  & 1890.28 & 251.40 &  & \textbf{1888.31} & 0.00 & 1896.83 & 74.06 & 0.45 \\ 
R107 & 100 & 1752.02 &  & 1753.35 & 249.00 &  & \textbf{1752.02} & 315.60 &  & 1753.35 & 0.08 & 1765.13 & 73.84 & 0.75 \\ 
R108 & 100 & 1647.88 &  & \textbf{1647.88} & 272.40 &  & 1649.37 & 238.20 &  & \textbf{1647.88} & 0.00 & 1664.56 & 67.97 & 1.01 \\ 
R109 & 100 & 1818.15 &  & \textbf{1818.15} & 219.60 &  & 1819.10 & 239.40 &  & \textbf{1818.15} & 0.00 & 1827.03 & 73.66 & 0.49 \\ 
R110 & 100 & 1758.64 &  & \textbf{1758.64} & 306.60 &  & 1761.96 & 328.20 &  & 1762.39 & 0.21 & 1773.40 & 78.90 & 0.84 \\ 
R111 & 100 & 1740.86 &  & \textbf{1740.86} & 319.20 &  & 1743.16 & 341.40 &  & 1743.16 & 0.13 & 1761.93 & 81.04 & 1.21 \\ 
R112 & 100 & 1661.85 &  & \textbf{1661.85} & 321.60 &  & 1663.09 & 300.60 &  & 1663.09 & 0.07 & 1672.36 & 68.90 & 0.63 \\ 
\noalign{\smallskip}\hline\noalign{\smallskip}
C101 & 100 & 2340.15 &  & \textbf{2340.15} & 187.20 &  & \textbf{2340.15} & 178.80 &  & \textbf{2340.15} & 0.00 & 2342.43 & 70.09 & 0.10 \\ 
C102 & 100 & 2325.70 &  & \textbf{2325.70} & 156.60 &  & \textbf{2325.70} & 163.80 &  & \textbf{2325.70} & 0.00 & 2325.70 & 73.58 & 0.00 \\ 
C103 & 100 & 2324.60 &  & \textbf{2324.60} & 181.80 &  & \textbf{2324.60} & 218.40 &  & \textbf{2324.60} & 0.00 & 2324.60 & 80.46 & 0.00 \\ 
C104 & 100 & 2318.04 &  & \textbf{2318.04} & 146.40 &  & \textbf{2318.04} & 178.80 &  & \textbf{2318.04} & 0.00 & 2318.04 & 70.17 & 0.00 \\ 
C105 & 100 & 2340.15 &  & \textbf{2340.15} & 180.00 &  & \textbf{2340.15} & 162.60 &  & \textbf{2340.15} & 0.00 & 2340.28 & 72.24 & 0.01 \\ 
C106 & 100 & 2340.15 &  & \textbf{2340.15} & 207.60 &  & \textbf{2340.15} & 191.40 &  & \textbf{2340.15} & 0.00 & 2340.84 & 68.22 & 0.03 \\ 
C107 & 100 & 2340.15 &  & \textbf{2340.15} & 192.00 &  & \textbf{2340.15} & 176.40 &  & \textbf{2340.15} & 0.00 & 2340.28 & 70.80 & 0.01 \\ 
C108 & 100 & 2338.58 &  & \textbf{2338.58} & 190.80 &  & \textbf{2338.58} & 232.80 &  & \textbf{2338.58} & 0.00 & 2338.58 & 84.86 & 0.00 \\ 
C109 & 100 & 2328.55 &  & \textbf{2328.55} & 163.20 &  & \textbf{2328.55} & 187.20 &  & \textbf{2328.55} & 0.00 & 2328.55 & 80.69 & 0.00 \\ 
\noalign{\smallskip}\hline\noalign{\smallskip}
RC101 & 100 & 2407.43 &  & 2412.71 & 249.60 &  & \textbf{2407.43} & 207.60 &  & 2412.71 & 0.22 & 2416.42 & 67.53 & 0.37 \\ 
RC102 & 100 & 2213.92 &  & \textbf{2213.92} & 351.60 &  & 2219.23 & 308.40 &  & \textbf{2213.92} & 0.00 & 2221.66 & 70.19 & 0.35 \\ 
RC103 & 100 & 2015.55 &  & 2016.28 & 210.00 &  & \textbf{2015.55} & 221.40 &  & 2016.28 & 0.04 & 2022.06 & 75.49 & 0.32 \\ 
RC104 & 100 & 1896.40 &  & 1897.04 & 244.20 &  & \textbf{1896.40} & 274.20 &  & 1907.56 & 0.59 & 1920.20 & 68.47 & 1.26 \\ 
RC105 & 100 & 2274.28 &  & 2287.51 & 272.40 &  & \textbf{2274.28} & 341.40 &  & 2287.51 & 0.58 & 2299.23 & 65.91 & 1.10 \\ 
RC106 & 100 & 2132.13 &  & 2140.86 & 237.00 &  & \textbf{2132.13} & 187.20 &  & 2139.36 & 0.34 & 2143.93 & 68.42 & 0.55 \\ 
RC107 & 100 & 1984.67 &  & 1989.34 & 179.40 &  & \textbf{1984.67} & 147.00 &  & 1989.34 & 0.24 & 1997.33 & 63.66 & 0.64 \\ 
RC108 & 100 & 1895.97 &  & 1898.96 & 234.60 &  & \textbf{1895.97} & 160.20 &  & 1898.96 & 0.16 & 1905.85 & 61.43 & 0.52 \\ 
\noalign{\smallskip}\hline\noalign{\smallskip}
R201 & 100 & 1646.78 &  & \textbf{1646.78} & 357.60 &  & \textbf{1646.78} & 407.40 &  & \textbf{1646.78} & 0.00 & 1658.07 & 178.63 & 0.69 \\ 
R202 & 100 & 1501.81 &  & 1508.16 & 549.00 &  & \textbf{1501.81} & 433.80 &  & 1510.25 & 0.56 & 1516.94 & 196.30 & 1.01 \\ 
R203 & 100 & 1341.09 &  & \textbf{1341.09} & 225.00 &  & \textbf{1341.09} & 273.60 &  & \textbf{1341.09} & 0.00 & 1346.55 & 183.62 & 0.41 \\ 
R204 & 100 & 1218.14 &  & \textbf{1218.14} & 306.60 &  & \textbf{1218.14} & 246.60 &  & \textbf{1218.14} & 0.00 & 1222.38 & 169.45 & 0.35 \\ 
R205 & 100 & 1418.97 &  & \textbf{1418.97} & 436.80 &  & 1420.81 & 388.20 &  & \textbf{1418.97} & 0.00 & 1429.55 & 176.47 & 0.75 \\ 
R206 & 100 & 1346.34 &  & \textbf{1346.34} & 433.80 &  & 1347.41 & 419.40 &  & 1349.12 & 0.21 & 1360.02 & 190.02 & 1.02 \\ 
R207 & 100 & 1277.58 &  & \textbf{1277.58} & 376.80 &  & 1278.57 & 406.80 &  & \textbf{1277.58} & 0.00 & 1282.99 & 179.80 & 0.42 \\ 
R208 & 100 & 1197.24 &  & \textbf{1197.24} & 259.80 &  & 1198.70 & 328.20 &  & \textbf{1197.24} & 0.00 & 1200.17 & 167.17 & 0.24 \\ 
R209 & 100 & 1322.42 &  & \textbf{1322.42} & 369.00 &  & \textbf{1322.42} & 328.20 &  & \textbf{1322.42} & 0.00 & 1330.11 & 182.65 & 0.58 \\ 
R210 & 100 & 1370.41 &  & 1374.31 & 418.20 &  & \textbf{1370.41} & 355.80 &  & 1376.40 & 0.44 & 1382.57 & 186.79 & 0.89 \\ 
R211 & 100 & 1219.93 &  & \textbf{1219.93} & 407.40 &  & 1220.57 & 468.60 &  & \textbf{1219.93} & 0.00 & 1222.20 & 169.81 & 0.19 \\ 
\noalign{\smallskip}\hline\noalign{\smallskip}
C201 & 100 & 1695.02 &  & \textbf{1695.02} & 178.20 &  & \textbf{1695.02} & 126.60 &  & \textbf{1695.02} & 0.00 & 1695.02 & 201.24 & 0.00 \\ 
C202 & 100 & 1685.24 &  & \textbf{1685.24} & 169.80 &  & \textbf{1685.24} & 139.80 &  & \textbf{1685.24} & 0.00 & 1685.24 & 165.84 & 0.00 \\ 
C203 & 100 & 1681.55 &  & \textbf{1681.55} & 201.60 &  & \textbf{1681.55} & 154.20 &  & \textbf{1681.55} & 0.00 & 1682.06 & 160.47 & 0.03 \\ 
C204 & 100 & 1677.66 &  & \textbf{1677.66} & 208.20 &  & \textbf{1677.66} & 221.40 &  & \textbf{1677.66} & 0.00 & 1678.07 & 177.00 & 0.02 \\ 
C205 & 100 & 1691.36 &  & \textbf{1691.36} & 192.60 &  & \textbf{1691.36} & 184.20 &  & \textbf{1691.36} & 0.00 & 1691.36 & 170.75 & 0.00 \\ 
C206 & 100 & 1689.32 &  & \textbf{1689.32} & 179.40 &  & \textbf{1689.32} & 191.40 &  & \textbf{1689.32} & 0.00 & 1689.32 & 165.75 & 0.00 \\ 
C207 & 100 & 1687.35 &  & \textbf{1687.35} & 206.40 &  & \textbf{1687.35} & 225.60 &  & \textbf{1687.35} & 0.00 & 1687.35 & 165.04 & 0.00 \\ 
C208 & 100 & 1686.50 &  & \textbf{1686.50} & 194.40 &  & \textbf{1686.50} & 144.60 &  & \textbf{1686.50} & 0.00 & 1686.50 & 150.51 & 0.00 \\ 
\noalign{\smallskip}\hline\noalign{\smallskip}
RC201 & 100 & 1938.36 &  & \textbf{1938.36} & 356.40 &  & 1941.16 & 418.80 &  & \textbf{1938.36} & 0.00 & 1943.46 & 131.74 & 0.26 \\ 
RC202 & 100 & 1768.04 &  & 1772.81 & 355.20 &  & \textbf{1768.04} & 388.20 &  & 1772.49 & 0.25 & 1773.94 & 138.32 & 0.33 \\ 
RC203 & 100 & 1603.55 &  & 1604.04 & 359.40 &  & \textbf{1603.55} & 369.00 &  & 1604.03 & 0.03 & 1622.55 & 153.21 & 1.18 \\ 
RC204 & 100 & 1489.27 &  & 1490.25 & 256.80 &  & \textbf{1489.27} & 208.20 &  & 1490.25 & 0.07 & 1492.28 & 142.89 & 0.20 \\ 
RC205 & 100 & 1832.53 &  & \textbf{1832.53} & 306.60 &  & 1833.34 & 238.80 &  & \textbf{1832.53} & 0.00 & 1836.37 & 139.21 & 0.21 \\ 
RC206 & 100 & 1724.41 &  & 1725.44 & 347.40 &  & \textbf{1724.41} & 272.40 &  & 1725.44 & 0.06 & 1736.04 & 159.45 & 0.67 \\ 
RC207 & 100 & 1646.37 &  & \textbf{1646.37} & 339.00 &  & 1650.23 & 300.60 &  & \textbf{\textit{1644.44}} & -0.12 & 1653.60 & 157.99 & 0.44 \\ 
RC208 & 100 & 1481.74 &  & 1483.20 & 263.40 &  & \textbf{1481.74} & 244.80 &  & 1483.20 & 0.10 & 1486.01 & 156.95 & 0.29 \\ 
\noalign{\smallskip}\hline\noalign{\smallskip}
\multicolumn{2}{l}{Average} & \multicolumn{1}{l}{} &  & 0.09 & 269.98 &  & 0.04 & 262.76 &  & \multicolumn{1}{l}{} & 0.11 & \multicolumn{1}{l}{} & 117.77 & 0.42 \\ 
\hline

\end{tabular}
}
\label{tbl:Resultados.HH.FSMTW_Bn}
\end{table}

\begin{table}[htbp]
\caption{Results for the FSMTW (\small{minimize distance, fleet C})}
\setlength{\tabcolsep}{1.5mm}
\scalebox{0.78} {
\begin{tabular}{lrrlrrlrrlrrrrr}
\hline\noalign{\smallskip}
 &  &  &  & \multicolumn{2}{l}{VCGP14} &  & \multicolumn{2}{l}{KBJL15} &  & \multicolumn{5}{l}{HILS-RVRP} \\
\cline{5-6} \cline{8-9} \cline{11-15}\noalign{\smallskip}
\multicolumn{1}{c}{Inst.} & \multicolumn{1}{c}{n} & \multicolumn{1}{c}{BKS} &  & \multicolumn{1}{c}{Best Sol.} & \multicolumn{1}{c}{{Avg. T (s)}} &  & \multicolumn{1}{c}{Best Sol.} & {{Avg. T (s)}} & \multicolumn{1}{c}{} & \multicolumn{1}{c}{Best Sol.} & \multicolumn{1}{c}{Gap (\%)} & \multicolumn{1}{c}{{Avg. Sol.}} & \multicolumn{1}{c}{{Avg. T(s)}} & \multicolumn{1}{c}{{Avg. Gap (\%)}} \\
\noalign{\smallskip}\hline\noalign{\smallskip}
R101 & 100 & 1937.38 &  & 1951.20 & 276.00 &  & \textbf{1937.38} & 250.20 &  & 1951.20 & 0.71 & 1951.64 & 65.13 & 0.74 \\ 
R102 & 100 & 1762.22 &  & 1785.35 & 175.20 & \multicolumn{1}{c}{} & \textbf{1762.22} & 193.80 &  & 1778.29 & 0.91 & 1782.92 & 70.67 & 1.17 \\ 
R103 & 100 & 1546.98 &  & 1552.34 & 213.00 & \multicolumn{1}{c}{} & \textbf{1546.98} & 221.40 &  & 1550.73 & 0.24 & 1552.81 & 82.77 & 0.38 \\ 
R104 & 100 & 1352.37 &  & 1355.15 & 322.20 & \multicolumn{1}{c}{} & \textbf{1352.37} & 310.20 &  & 1355.15 & 0.21 & 1366.54 & 72.12 & 1.05 \\ 
R105 & 100 & 1681.44 &  & 1694.56 & 195.00 &  & \textbf{1681.44} & 247.80 &  & 1694.56 & 0.78 & 1696.65 & 69.18 & 0.90 \\ 
R106 & 100 & 1583.17 &  & \textbf{1583.17} & 247.20 &  & 1585.65 & 220.20 &  & \textbf{1583.17} & 0.00 & 1590.23 & 78.15 & 0.45 \\ 
R107 & 100 & 1424.37 &  & 1428.08 & 321.60 &  & \textbf{1424.37} & 358.80 &  & 1428.08 & 0.26 & 1444.99 & 77.80 & 1.45 \\ 
R108 & 100 & 1314.88 &  & \textbf{1314.88} & 319.20 &  & 1318.44 & 286.80 &  & \textbf{1314.88} & 0.00 & 1330.05 & 71.87 & 1.15 \\ 
R109 & 100 & 1506.59 &  & \textbf{1506.59} & 280.80 &  & 1507.10 & 246.60 &  & \textbf{1506.59} & 0.00 & 1511.16 & 71.68 & 0.30 \\ 
R110 & 100 & 1443.37 &  & 1443.92 & 303.60 &  & \textbf{1443.37} & 286.80 &  & \textbf{\textit{1439.42}} & -0.27 & 1452.53 & 86.86 & 0.63 \\ 
R111 & 100 & 1419.43 &  & 1420.15 & 334.80 &  & \textbf{1419.43} & 308.40 &  & 1423.41 & 0.28 & 1439.63 & 79.32 & 1.42 \\ 
R112 & 100 & 1327.58 &  & \textbf{1327.58} & 298.20 &  & 1328.01 & 280.20 &  & 1328.47 & 0.07 & 1346.27 & 71.60 & 1.41 \\ 
\noalign{\smallskip}\hline\noalign{\smallskip}
C101 & 100 & 1628.94 &  & \textbf{1628.94} & 127.20 &  & \textbf{1628.94} & 119.40 &  & \textbf{1628.94} & 0.00 & 1628.94 & 56.82 & 0.00 \\ 
C102 & 100 & 1597.66 &  & \textbf{1597.66} & 141.00 &  & \textbf{1597.66} & 128.40 &  & \textbf{1597.66} & 0.00 & 1597.84 & 60.67 & 0.01 \\ 
C103 & 100 & 1596.56 &  & \textbf{1596.56} & 172.80 &  & \textbf{1596.56} & 159.00 &  & \textbf{1596.56} & 0.00 & 1596.56 & 65.06 & 0.00 \\ 
C104 & 100 & 1590.76 &  & \textbf{1590.76} & 139.20 &  & \textbf{1590.76} & 126.60 &  & \textbf{1590.76} & 0.00 & 1590.76 & 62.50 & 0.00 \\ 
C105 & 100 & 1628.94 &  & \textbf{1628.94} & 117.60 &  & \textbf{1628.94} & 144.60 &  & \textbf{1628.94} & 0.00 & 1628.94 & 57.33 & 0.00 \\ 
C106 & 100 & 1628.94 &  & \textbf{1628.94} & 123.00 &  & \textbf{1628.94} & 104.40 &  & \textbf{1628.94} & 0.00 & 1628.94 & 55.11 & 0.00 \\ 
C107 & 100 & 1628.94 &  & \textbf{1628.94} & 130.20 &  & \textbf{1628.94} & 121.80 &  & \textbf{1628.94} & 0.00 & 1628.94 & 58.33 & 0.00 \\ 
C108 & 100 & 1622.75 &  & \textbf{1622.75} & 192.00 &  & \textbf{1622.75} & 153.60 &  & \textbf{1622.75} & 0.00 & 1628.32 & 62.69 & 0.34 \\ 
C109 & 100 & 1614.99 &  & 1615.93 & 232.80 &  & \textbf{1614.99} & 178.20 &  & \textbf{1614.99} & 0.00 & 1616.05 & 70.56 & 0.07 \\ 
\noalign{\smallskip}\hline\noalign{\smallskip}
RC101 & 100 & 2033.89 &  & 2043.48 & 263.40 &  & \textbf{2033.89} & 249.60 &  & 2040.61 & 0.33 & 2046.28 & 63.45 & 0.61 \\ 
RC102 & 100 & 1847.92 &  & \textbf{1847.92} & 246.00 &  & \textbf{1847.92} & 241.80 &  & \textbf{1847.92} & 0.00 & 1859.92 & 70.91 & 0.65 \\ 
RC103 & 100 & 1646.35 &  & \textbf{1646.35} & 251.40 &  & \textbf{1646.35} & 250.20 &  & \textbf{1646.35} & 0.00 & 1656.22 & 77.59 & 0.60 \\ 
RC104 & 100 & 1518.96 &  & 1522.04 & 338.40 &  & \textbf{1518.96} & 308.40 &  & 1522.04 & 0.20 & 1544.84 & 67.29 & 1.70 \\ 
RC105 & 100 & 1884.92 &  & 1913.06 & 240.60 &  & \textbf{1884.92} & 274.20 &  & 1913.06 & 1.49 & 1926.98 & 64.25 & 2.23 \\ 
RC106 & 100 & 1753.99 &  & 1770.95 & 226.20 &  & \textbf{1753.99} & 206.40 &  & 1770.95 & 0.97 & 1772.43 & 65.58 & 1.05 \\ 
RC107 & 100 & 1601.12 &  & 1607.11 & 244.80 &  & \textbf{1601.12} & 208.20 &  & 1607.11 & 0.37 & 1614.00 & 63.87 & 0.80 \\ 
RC108 & 100 & 1516.36 &  & 1523.96 & 201.00 &  & \textbf{1516.36} & 218.40 &  & 1523.96 & 0.50 & 1525.47 & 60.51 & 0.60 \\ 
\noalign{\smallskip}\hline\noalign{\smallskip}
R201 & 100 & 1429.50 &  & 1443.41 & 294.00 &  & \textbf{1429.50} & 272.40 &  & \textbf{\textit{1421.78}} & -0.54 & 1431.20 & 176.19 & 0.12 \\ 
R202 & 100 & 1273.11 &  & 1283.16 & 474.60 &  & \textbf{1273.11} & 427.20 &  & 1283.86 & 0.84 & 1291.28 & 189.45 & 1.43 \\ 
R203 & 100 & 1116.09 &  & \textbf{1116.09} & 237.00 &  & \textbf{1116.09} & 274.80 &  & \textbf{1116.09} & 0.00 & 1119.58 & 178.95 & 0.31 \\ 
R204 & 100 & 993.14 &  & \textbf{993.14} & 397.80 &  & \textbf{993.14} & 408.60 &  & \textbf{993.14} & 0.00 & 995.54 & 173.31 & 0.24 \\ 
R205 & 100 & 1193.97 &  & \textbf{1193.97} & 439.80 &  & 1195.81 & 372.60 &  & \textbf{1193.97} & 0.00 & 1207.74 & 174.27 & 1.15 \\ 
R206 & 100 & 1121.34 &  & \textbf{1121.34} & 360.00 &  & \textbf{1121.34} & 308.40 &  & 1123.43 & 0.19 & 1129.96 & 180.35 & 0.77 \\ 
R207 & 100 & 1052.58 &  & \textbf{1052.58} & 418.20 &  & \textbf{1052.58} & 313.80 &  & 1055.45 & 0.27 & 1061.29 & 171.75 & 0.83 \\ 
R208 & 100 & 969.90 &  & \textbf{969.90} & 346.80 &  & 973.70 & 328.20 &  & \textbf{969.90} & 0.00 & 974.31 & 170.43 & 0.45 \\ 
R209 & 100 & 1094.97 &  & 1097.42 & 356.40 &  & \textbf{1094.97} & 338.40 &  & 1097.42 & 0.22 & 1105.09 & 175.98 & 0.92 \\ 
R210 & 100 & 1145.48 &  & 1149.85 & 408.00 &  & \textbf{1145.48} & 370.20 &  & 1153.08 & 0.66 & 1157.78 & 187.50 & 1.07 \\ 
R211 & 100 & 994.93 &  & \textbf{994.93} & 390.60 &  & \textbf{994.93} & 370.20 &  & \textbf{994.93} & 0.00 & 997.74 & 167.85 & 0.28 \\ 
\noalign{\smallskip}\hline\noalign{\smallskip}
C201 & 100 & 1194.33 &  & \textbf{1194.33} & 289.20 &  & \textbf{1194.33} & 270.00 &  & \textbf{1194.33} & 0.00 & 1194.33 & 201.61 & 0.00 \\ 
C202 & 100 & 1185.24 &  & \textbf{1185.24} & 154.20 &  & \textbf{1185.24} & 141.60 &  & \textbf{1185.24} & 0.00 & 1185.24 & 168.50 & 0.00 \\ 
C203 & 100 & 1176.25 &  & \textbf{1176.25} & 216.00 &  & \textbf{1176.25} & 184.20 &  & \textbf{1176.25} & 0.00 & 1176.41 & 161.78 & 0.01 \\ 
C204 & 100 & 1175.37 &  & \textbf{1175.37} & 259.20 &  & \textbf{1175.37} & 185.40 &  & \textbf{1175.37} & 0.00 & 1175.37 & 182.77 & 0.00 \\ 
C205 & 100 & 1190.36 &  & \textbf{1190.36} & 267.60 &  & \textbf{1190.36} & 270.00 &  & \textbf{1190.36} & 0.00 & 1190.36 & 176.91 & 0.00 \\ 
C206 & 100 & 1188.62 &  & \textbf{1188.62} & 240.60 &  & \textbf{1188.62} & 239.40 &  & \textbf{1188.62} & 0.00 & 1188.62 & 162.28 & 0.00 \\ 
C207 & 100 & 1184.88 &  & \textbf{1184.88} & 219.60 &  & \textbf{1184.88} & 190.20 &  & \textbf{1184.88} & 0.00 & 1185.37 & 159.37 & 0.04 \\ 
C208 & 100 & 1186.50 &  & \textbf{1186.50} & 179.40 &  & \textbf{1186.50} & 172.20 &  & \textbf{1186.50} & 0.00 & 1186.50 & 156.34 & 0.00 \\ 
\noalign{\smallskip}\hline\noalign{\smallskip}
RC201 & 100 & 1623.36 &  & \textbf{1623.36} & 408.60 &  & 1625.71 & 360.60 &  & \textbf{1623.36} & 0.00 & 1627.62 & 125.38 & 0.26 \\ 
RC202 & 100 & 1445.12 &  & 1447.27 & 317.40 &  & \textbf{1445.12} & 247.20 &  & 1447.27 & 0.15 & 1452.60 & 140.32 & 0.52 \\ 
RC203 & 100 & 1273.55 &  & 1274.04 & 269.40 &  & \textbf{1273.55} & 220.20 &  & 1274.03 & 0.04 & 1291.12 & 152.44 & 1.38 \\ 
RC204 & 100 & 1157.94 &  & 1159.00 & 360.60 &  & \textbf{1157.94} & 308.40 &  & 1159.00 & 0.09 & 1161.80 & 152.38 & 0.33 \\ 
RC205 & 100 & 1512.53 &  & \textbf{1512.53} & 314.40 &  & 1515.34 & 300.60 &  & \textbf{1512.53} & 0.00 & 1517.15 & 136.56 & 0.31 \\ 
RC206 & 100 & 1395.18 &  & \textbf{1395.18} & 235.20 &  & 1399.41 & 196.20 &  & 1400.44 & 0.38 & 1412.04 & 161.53 & 1.21 \\ 
RC207 & 100 & 1314.44 &  & \textbf{1314.44} & 374.40 &  & 1317.50 & 328.20 &  & \textbf{1314.44} & 0.00 & 1324.55 & 154.91 & 0.77 \\ 
RC208 & 100 & 1140.10 & & 1140.10 & 408.60 &  & \textbf{1140.10} & 359.40 &  & \textbf{1140.10} & 0.00 & 1148.45 & 151.26 & 0.73 \\ 
\noalign{\smallskip}\hline\noalign{\smallskip}
\multicolumn{2}{l}{Average} & \multicolumn{1}{l}{} &  & 0.19 & 273.43 &  & 0.03 & 252.91 &  & \multicolumn{1}{l}{} & 0.17 & \multicolumn{1}{l}{} & 115.54 & 0.59 \\ 
\hline
\end{tabular}
}
\label{tbl:Resultados.HH.FSMTW_Cn}
\end{table}

\newpage
\subsection{HFFVRPMBTW} \label{sec:Result.HH.HFFVRPMBTW.c}

Detailed results obtained for the HFFVRPMBTW instances of \citep[BPYA13]{Belmecherietal2013}, compared with the results obtained by the PSO of the same authors and the EBBO of \cite[BB15]{BerghidaBoukra2015} (Tables~\ref{tbl:Resultados.HH.HFFVRPMBTW.c}-\ref{tbl:Resultados.HH.HFFVRPMBTW.rc}).

\begin{table}[htbp]
\caption{Results for the HFFVRPMBTW (Type C)}
\setlength{\tabcolsep}{1.5mm}
\scalebox{0.82} {
\begin{tabular}{lrrlrclrclrrrrr}
\hline\noalign{\smallskip} 
 &  &  &  & \multicolumn{2}{l}{PSO} &  \multicolumn{2}{l}{EBBO} & & & \multicolumn{5}{l}{HILS-RVRP} \\
  &  &  &  & \multicolumn{2}{l}{BPYA13} &  \multicolumn{2}{l}{BB15} &  & & & & & \\
\cline{5-6} \cline{8-9} \cline{11-15}\hline\noalign{\smallskip}
\multicolumn{1}{c}{Inst.} & \multicolumn{1}{c}{n} & \multicolumn{1}{c}{BKS} &  & \multicolumn{1}{c}{{Sol.}} & \multicolumn{1}{c}{T (s)} &  & \multicolumn{1}{c}{{Sol.}} & \multicolumn{1}{c}{T (s)} & \multicolumn{1}{c}{} & \multicolumn{1}{c}{Best Sol.} & \multicolumn{1}{c}{Gap (\%)} & \multicolumn{1}{c}{{Avg. Sol.}} & \multicolumn{1}{c}{{Avg. T(s)}} & \multicolumn{1}{c}{{Avg. Gap (\%)}} \\
\noalign{\smallskip}\hline\noalign{\smallskip}
C101 & 100 & 2331.54 &  & 2560.02 & -- &  & \textbf{2331.54} & -- &  & \textbf{\textit{1348.34}} & -42.17 & \textbf{\textit{1348.34}} & 79.32 & -42.17 \\ 
C102 & 100 & 2410.18 &  & 2615.32 & -- &  & \textbf{2410.18} & -- &  & \textbf{\textit{1344.72}} & -44.21 & \textbf{\textit{1344.72}} & 77.69 & -44.21 \\ 
C103 & 100 & 2102.41 &  & 2405.30 & -- &  & \textbf{2102.41} & -- &  & \textbf{\textit{1343.33}} & -36.11 & \textbf{\textit{1343.33}} & 72.41 & -36.11 \\ 
C104 & 100 & 2021.55 &  & 2333.95 & -- &  & \textbf{2021.55} & -- &  & \textbf{\textit{1332.78}} & -34.07 & \textbf{\textit{1333.68}} & 78.20 & -34.03 \\ 
C105 & 100 & 1998.83 &  & 2055.90 & -- &  & \textbf{1998.83} & -- &  & \textbf{\textit{1348.10}} & -32.56 & \textbf{\textit{1348.10}} & 80.27 & -32.56 \\ 
C106 & 100 & 2193.54 &  & 2366.05 & -- &  & \textbf{2193.54} & -- &  & \textbf{\textit{1348.34}} & -38.53 & \textbf{\textit{1348.34}} & 91.37 & -38.53 \\ 
C107 & 100 & 1992.83 &  & \textbf{1992.83} & -- &  & 2188.36 & -- &  & \textbf{\textit{1348.10}} & -32.35 & \textbf{\textit{1348.10}} & 82.29 & -32.35 \\ 
C108 & 100 & 1938.54 &  & \textbf{1938.54} & -- &  & 1950.36 & -- &  & \textbf{\textit{1348.10}} & -30.46 & \textbf{\textit{1348.10}} & 76.58 & -30.46 \\ 
C109 & 100 & 1786.66 &  & 2234.79 & -- &  & \textbf{1786.66} & -- &  & \textbf{\textit{1342.30}} & -24.87 & \textbf{\textit{1342.30}} & 60.48 & -24.87 \\ 
\noalign{\smallskip}\hline\noalign{\smallskip}
C201 & 100 & 1298.17 &  & 1420.62 & -- &  & 1326.76 & -- &  & \textbf{\textit{1172.71}} & -9.66 & \textbf{\textit{1174.37}} & 99.42 & -9.54 \\ 
C202 & 100 & 1327.33 &  & 1590.20 & -- &  & \textbf{1327.33} & -- &  & \textbf{\textit{1172.71}} & -11.65 & \textbf{\textit{1174.06}} & 99.10 & -11.55 \\ 
C203 & 100 & 1432.36 &  & 1823.84 & -- &  & \textbf{1432.36} & -- &  & \textbf{\textit{1160.29}} & -18.99 & \textbf{\textit{1167.44}} & 105.78 & -18.50 \\ 
C204 & 100 & 1727.46 &  & 1856.26 & -- &  & 1884.86 & -- &  & \textbf{\textit{1162.24}} & -32.72 & \textbf{\textit{1162.24}} & 119.88 & -32.72 \\ 
C205 & 100 & 1296.01 &  & 1504.98 & -- &  & \textbf{1296.01} & -- &  & \textbf{\textit{1167.88}} & -9.89 & \textbf{\textit{1168.89}} & 103.38 & -9.81 \\ 
C206 & 100 & 1444.86 &  & 1528.31 & -- &  & 1687.04 & -- &  & \textbf{\textit{1163.11}} & -19.50 & \textbf{\textit{1168.35}} & 97.39 & -19.14 \\ 
C207 & 100 & 1376.64 &  & 1391.91 & -- &  & 1416.27 & -- &  & \textbf{\textit{1166.22}} & -15.29 & \textbf{\textit{1166.32}} & 91.84 & -15.28 \\ 
C208 & 100 & 1453.05 &  & 1626.96 & -- &  & \textbf{1453.05} & -- &  & \textbf{\textit{1166.78}} & -19.70 & \textbf{\textit{1169.48}} & 96.51 & -19.52 \\ 
\noalign{\smallskip}\hline\noalign{\smallskip}
 \multicolumn{2}{l}{Average} &  &  &  10.76 &  &  &  2.43 & &  &  & -26.63 &  & 88.94 & -26.55 \\
 \hline
\end{tabular}
}
\label{tbl:Resultados.HH.HFFVRPMBTW.c}
\end{table}

\begin{table}[htbp]
\caption{Results for the HFFVRPMBTW (Type R)}
\setlength{\tabcolsep}{1.5mm}
\scalebox{0.82} {
\begin{tabular}{lrrlrclrclrrrrr}
\hline\noalign{\smallskip} 
 &  &  &  & \multicolumn{2}{l}{PSO} &  & \multicolumn{2}{l}{EBBO} &  & \multicolumn{5}{l}{HILS-RVRP}  \\
   &  &  &  & \multicolumn{2}{l}{BPYA13} &  \multicolumn{2}{l}{BB15} &  & & & & & \\
\cline{5-6} \cline{8-9} \cline{11-15}\noalign{\smallskip}
\multicolumn{1}{c}{Inst.} & \multicolumn{1}{c}{n} & \multicolumn{1}{c}{BKS} &  & \multicolumn{1}{c}{{Sol.}} & \multicolumn{1}{c}{T (s)} &  & \multicolumn{1}{c}{{Sol.}} & \multicolumn{1}{c}{T (s)} & \multicolumn{1}{c}{} & \multicolumn{1}{c}{Best Sol.} & \multicolumn{1}{c}{Gap (\%)} & \multicolumn{1}{c}{{Avg. Sol.}} & \multicolumn{1}{c}{{Avg. T(s)}} & \multicolumn{1}{c}{{Avg. Gap (\%)}} \\
\noalign{\smallskip}\hline\noalign{\smallskip}
R101 & 100 & 2567.14 &  & 2632.13 & -- &  & \textbf{2567.14} & -- &  & \textbf{\textit{2239.53}} & -12.76 & \textbf{\textit{2247.28}} & 99.69 & -12.46 \\ 
R102 & 100 & 2369.30 &  & 2375.41 & -- &  & \textbf{2369.30} & -- &  & \textbf{\textit{1988.11}} & -16.09 & \textbf{\textit{2000.44}} & 107.87 & -15.57 \\ 
R103 & 100 & 2006.80 &  & \textbf{2006.80} & -- &  & 2080.34 & -- &  & \textbf{\textit{1658.71}} & -17.35 & \textbf{\textit{1680.60}} & 86.95 & -16.25 \\ 
R104 & 100 & 1853.21 &  & \textbf{1853.21} & -- &  & 1971.83 & -- &  & \textbf{\textit{1444.76}} & -22.04 & \textbf{\textit{1461.90}} & 82.12 & -21.12 \\ 
R105 & 100 & 2253.42 &  & \textbf{2253.42} & -- &  & 2334.56 & -- &  & \textbf{\textit{1873.51}} & -16.86 & \textbf{\textit{1881.50}} & 98.75 & -16.50 \\ 
R106 & 100 & 2031.19 &  & \textbf{2031.19} & -- &  & 2121.66 & -- &  & \textbf{\textit{1699.09}} & -16.35 & \textbf{\textit{1714.60}} & 110.12 & -15.59 \\ 
R107 & 100 & 1905.43 &  & 1928.90 & -- &  & 1943.65 & -- &  & \textbf{\textit{1546.98}} & -18.81 & \textbf{\textit{1567.01}} & 102.49 & -17.76 \\ 
R108 & 100 & 1877.52 &  & \textbf{1877.52} & -- &  & 2002.36 & -- &  & \textbf{\textit{1392.10}} & -25.85 & \textbf{\textit{1416.40}} & 80.96 & -24.56 \\ 
R109 & 100 & 2001.56 &  & \textbf{2001.56} & -- &  & 2069.38 & -- &  & \textbf{\textit{1621.53}} & -18.99 & \textbf{\textit{1641.62}} & 101.90 & -17.98 \\ 
R110 & 100 & 1979.53 &  & 1983.98 & -- &  & 2065.76 & -- &  & \textbf{\textit{1551.83}} & -21.61 & \textbf{\textit{1576.86}} & 99.14 & -20.34 \\ 
R111 & 100 & 1881.21 &  & 1896.70 & -- &  & \textbf{1881.21} & -- &  & \textbf{\textit{1524.55}} & -18.96 & \textbf{\textit{1551.63}} & 93.69 & -17.52 \\ 
R112 & 100 & 1689.12 &  & 1895.77 & -- &  & \textbf{1689.12} & -- &  & \textbf{\textit{1404.26}} & -16.86 & \textbf{\textit{1423.15}} & 86.35 & -15.75 \\ 
\noalign{\smallskip}\hline\noalign{\smallskip}
R201 & 100 & 1344.47 &  & 1990.47 & -- &  & \textbf{1344.47} & -- &  & 1637.34 & 21.78 & 1661.42 & 132.41 & 23.57 \\ 
R202 & 100 & 1922.72 &  & 1932.74 & -- &  & 1941.04 & -- &  & \textbf{\textit{1573.33}} & -18.17 & \textbf{\textit{1596.65}} & 157.32 & -16.96 \\ 
R203 & 100 & 1736.20 &  & 1745.37 & -- &  & 1910.74 & -- &  & \textbf{\textit{1419.67}} & -18.23 & \textbf{\textit{1430.72}} & 141.93 & -17.59 \\ 
R204 & 100 & 1522.50 &  & \textbf{1522.50} & -- &  & 1876.82 & -- &  & \textbf{\textit{1260.50}} & -17.21 & \textbf{\textit{1275.06}} & 154.38 & -16.25 \\ 
R205 & 100 & 1753.49 &  & 1885.75 & -- &  & \textbf{1753.49} & -- &  & \textbf{\textit{1474.12}} & -15.93 & \textbf{\textit{1505.58}} & 138.98 & -14.14 \\ 
R206 & 100 & 1758.78 &  & 1813.48 & -- &  & 1792.46 & -- &  & \textbf{\textit{1425.15}} & -18.97 & \textbf{\textit{1444.85}} & 153.72 & -17.85 \\ 
R207 & 100 & 1650.12 &  & 1654.84 & -- &  & 1806.25 & -- &  & \textbf{\textit{1333.60}} & -19.18 & \textbf{\textit{1353.71}} & 141.27 & -17.96 \\ 
R208 & 100 & 1536.68 &  & 1589.42 & -- &  & 1737.67 & -- &  & \textbf{\textit{1244.82}} & -18.99 & \textbf{\textit{1256.92}} & 155.73 & -18.21 \\ 
R209 & 100 & 1729.58 &  & \textbf{1729.58} & -- &  & 1798.76 & -- &  & \textbf{\textit{1375.12}} & -20.49 & \textbf{\textit{1415.13}} & 136.16 & -18.18 \\ 
R210 & 100 & 1754.44 &  & \textbf{1754.44} & -- &  & 1868.55 & -- &  & \textbf{\textit{1396.01}} & -20.43 & \textbf{\textit{1430.41}} & 125.22 & -18.47 \\ 
R211 & 100 & 1615.85 &  & 1699.39 & -- &  & 1641.38 & -- &  & \textbf{\textit{1267.12}} & -21.58 & \textbf{\textit{1282.46}} & 117.58 & -20.63 \\ 
\noalign{\smallskip}\hline\noalign{\smallskip}
\multicolumn{2}{l}{Average} &  &  &  3.74  &  &  & 4.58 & &  &  & -16.95 &  & 117.60 & -15.83 \\
\hline
\end{tabular}
}
\label{tbl:Resultados.HH.HFFVRPMBTW.r}
\end{table}

\begin{table}[htbp]
\caption{Results for the HFFVRPMBTW (Type RC)}
\setlength{\tabcolsep}{1.5mm}
\scalebox{0.82} {
\begin{tabular}{lrrlrclrclrrrrr}
\hline\noalign{\smallskip}
 &  &  &  & \multicolumn{2}{l}{PSO} &  & \multicolumn{2}{l}{EBBO} &  & \multicolumn{5}{l}{HILS-RVRP}  \\
   &  &  &  & \multicolumn{2}{l}{BPYA13} &  \multicolumn{2}{l}{BB15} &  & & & & & \\
\cline{5-6} \cline{8-9} \cline{11-15}\noalign{\smallskip}
\multicolumn{1}{c}{Inst.} & \multicolumn{1}{c}{n} & \multicolumn{1}{c}{BKS} &  & \multicolumn{1}{c}{{Sol.}} & \multicolumn{1}{c}{T (s)} &  & \multicolumn{1}{c}{{Sol.}} & \multicolumn{1}{c}{T (s)} & \multicolumn{1}{c}{} & \multicolumn{1}{c}{Best Sol.} & \multicolumn{1}{c}{Gap (\%)} & \multicolumn{1}{c}{{Avg. Sol.}} & \multicolumn{1}{c}{{Avg. T(s)}} & \multicolumn{1}{c}{{Avg. Gap (\%)}} \\
\noalign{\smallskip}\hline\noalign{\smallskip}
RC101 & 100 & 2387.96 &  & 2957.49 & -- &  & \textbf{2387.96} & -- &  & \textbf{\textit{2314.12}} & -3.09 & \textbf{\textit{2329.56}} & 89.90 & -2.45 \\ 
RC102 & 100 & 2464.51 &  & \textbf{2464.51} & -- &  & 2664.55 & -- &  & \textbf{\textit{2069.34}} & -16.03 & \textbf{\textit{2076.25}} & 84.58 & -15.75 \\ 
RC103 & 100 & 2426.88 &  & \textbf{2426.88} & -- &  & 2553.62 & -- &  & \textbf{\textit{1904.09}} & -21.54 & \textbf{\textit{1927.90}} & 91.91 & -20.56 \\ 
RC104 & 100 & 2244.58 &  & \textbf{2244.58} & -- &  & 2253.76 & -- &  & \textbf{\textit{1691.27}} & -24.65 & \textbf{\textit{1708.22}} & 81.14 & -23.90 \\ 
RC105 & 100 & 2385.27 &  & 2711.05 & -- &  & \textbf{2385.27} & -- &  & \textbf{\textit{2155.60}} & -9.63 & \textbf{\textit{2171.73}} & 86.08 & -8.95 \\ 
RC106 & 100 & 2254.16 &  & 2495.57 & -- &  & \textbf{2254.16} & -- &  & \textbf{\textit{1968.21}} & -12.69 & \textbf{\textit{1987.42}} & 90.58 & -11.83 \\ 
RC107 & 100 & 2414.86 &  & 2420.42 & -- &  & 2420.42 & -- &  & \textbf{\textit{1809.99}} & -25.05 & \textbf{\textit{1825.15}} & 89.60 & -24.42 \\ 
RC108 & 100 & 2166.96 &  & 2381.45 & -- &  & \textbf{2166.96} & -- &  & \textbf{\textit{1670.53}} & -22.91 & \textbf{\textit{1694.15}} & 70.62 & -21.82 \\ 
\noalign{\smallskip}\hline\noalign{\smallskip}
RC201 & 100 & 2401.11 &  & 2401.11 & -- &  & 2571.17 & -- &  & \textbf{\textit{1965.54}} & -18.14 & \textbf{\textit{1998.80}} & 111.73 & -16.76 \\ 
RC202 & 100 & 2100.22 &  & 2251.39 & -- &  & 2100.22 & -- &  & \textbf{\textit{1782.27}} & -15.14 & \textbf{\textit{1800.03}} & 122.01 & -14.29 \\ 
RC203 & 100 & 1931.69 &  & 2022.90 & -- &  & 1941.74 & -- &  & \textbf{\textit{1600.20}} & -17.16 & \textbf{\textit{1620.29}} & 121.84 & -16.12 \\ 
RC204 & 100 & 1673.10 &  & 1827.48 & -- &  & 1673.10 & -- &  & \textbf{\textit{1422.68}} & -14.97 & \textbf{\textit{1428.31}} & 136.25 & -14.63 \\ 
RC205 & 100 & 2226.16 &  & 2274.91 & -- &  & 2304.21 & -- &  & \textbf{\textit{1846.52}} & -17.05 & \textbf{\textit{1869.31}} & 118.44 & -16.03 \\ 
RC206 & 100 & 1953.99 &  & 2123.08 & -- &  & \textbf{1953.99} & -- &  & \textbf{\textit{1779.05}} & -8.95 & \textbf{\textit{1793.24}} & 131.24 & -8.23 \\ 
RC207 & 100 & 1867.97 &  & 2084.50 & -- &  & 1867.97 & -- &  & \textbf{\textit{1631.16}} & -12.68 & \textbf{\textit{1636.28}} & 114.72 & -12.40 \\ 
RC208 & 100 & 1836.63 &  & \textbf{1836.63} & -- &  & \textbf{1836.63} & -- &  & \textbf{\textit{1396.52}} & -23.96 & \textbf{\textit{1406.85}} & 120.99 & -23.40 \\ 
\noalign{\smallskip}\hline\noalign{\smallskip}
\multicolumn{2}{l}{Average} &  &  &  6.37  &  &  & 1.57  & &  &  & -16.48 &  & 103.85 & -15.72 \\
\hline
\end{tabular}
} 
\label{tbl:Resultados.HH.HFFVRPMBTW.rc}
\end{table}


\newpage
\subsection{SDepVRPTW} \label{sec:Result.H.SDepVRPTW}

Detailed results obtained for the SDepVRPTW instances of \cite[CL01]{CordeauLaporte2001}, compared with those found by the ITS1 of \cite[CM12]{CordeauMaischberger2012} and by the HGSADC of \cite[VCGP13]{Vidaletal2013} (Table~\ref{tbl:Resultados.HH.SDepVRPTW}).

\begin{table}[htbp]
\caption{Results for the SDepVRPTW}
\setlength{\tabcolsep}{1.5mm}
\scalebox{0.8} {
\begin{tabular}{lrrlrclrrlrrrrr}
\hline\noalign{\smallskip} 
 &  &  &  & \multicolumn{2}{l}{ITS1} & & \multicolumn{2}{l}{HGSADC} & & \multicolumn{5}{l}{HILS-RVRP} \\
 &  &  &  & \multicolumn{2}{l}{CM12} & & \multicolumn{2}{l}{VCGP13} & &  \\
\cline{5-6} \cline{8-9} \cline{11-15} \noalign{\smallskip}
\multicolumn{1}{c}{Inst.} & \multicolumn{1}{c}{n} & \multicolumn{1}{c}{BKS} &  & \multicolumn{1}{c}{Best Sol.} & \multicolumn{1}{c}{T (s)} &  & \multicolumn{1}{c}{Best Sol.} & \multicolumn{1}{c}{{Avg. T (s)}} & \multicolumn{1}{c}{} & \multicolumn{1}{c}{Best Sol.} & \multicolumn{1}{c}{Gap (\%)} & \multicolumn{1}{c}{{Avg. Sol.}} & \multicolumn{1}{c}{{Avg. T(s)}} & \multicolumn{1}{c}{{Avg. Gap (\%)}} \\
\hline 
p01a & 48 & 1655.42 &  & \textbf{1655.42} & -- &  & \textbf{1655.42} & 13.80 &  & \textbf{1655.42} & 0.00 & \textbf{1655.42} & 18.34 & 0.00 \\ 
p02a & 96 & 2904.13 &  & \textbf{2904.13} & -- &  & \textbf{2904.13} & 42.00 &  & \textbf{2904.13} & 0.00 & 2906.01 & 134.18 & 0.06 \\ 
p03a & 144 & 3304.13 &  & 3317.33 & -- &  & \textbf{3304.13} & 96.00 &  & 3304.91 & 0.02 & 3317.04 & 522.57 & 0.39 \\ 
p04a & 192 & 4427.25 &  & 4461.13 & -- &  & \textbf{4427.25} & 351.00 &  & 4437.95 & 0.24 & 4488.39 & 1516.05 & 1.38 \\ 
p05a & 240 & 5626.42 &  & 5663.32 & -- &  & 5647.76 & 698.40 &  & 5642.37 & 0.28 & 5729.89 & 2505.33 & 1.84 \\ 
p06a & 288 & 5627.82 &  & 5698.93 & -- &  & 5637.48 & 760.80 &  & 5658.17 & 0.54 & 5710.12 & 4601.31 & 1.46 \\ 
p07a & 72 & 2166.88 &  & \textbf{2166.88} & -- &  & \textbf{2166.88} & 25.20 &  & \textbf{2166.88} & 0.00 & \textbf{2166.88} & 64.52 & 0.00 \\ 
p08a & 144 & 3873.40 &  & 3880.58 & -- &  & \textbf{3873.40} & 141.00 &  & \textbf{3873.40} & 0.00 & 3887.65 & 502.51 & 0.37 \\ 
p09a & 216 & 4772.55 &  & 4818.32 & -- &  & 4777.61 & 336.00 &  & 4807.58 & 0.73 & 4834.70 & 1196.14 & 1.30 \\ 
p10a & 288 & 5817.28 &  & 5908.53 & -- &  & 5858.82 & 694.80 &  & 5882.78 & 1.13 & 5926.17 & 3211.65 & 1.87 \\ 
\hline
p01b & 48 & 1429.35 &  & \textbf{1429.35} & -- &  & \textbf{1429.35} & 13.20 &  & \textbf{1429.35} & 0.00 & \textbf{1429.35} & 16.58 & 0.00 \\ 
p02b & 96 & 2479.56 &  & \textbf{2479.56} & -- &  & \textbf{2479.56} & 59.40 &  & \textbf{2479.56} & 0.00 & 2481.20 & 151.82 & 0.07 \\ 
p03b & 144 & 2774.30 &  & 2781.22 & -- &  & 2775.61 & 136.80 &  & 2776.45 & 0.08 & 2786.42 & 499.96 & 0.44 \\ 
p04b & 192 & 3649.72 &  & 3674.53 & -- &  & \textbf{3649.72} & 394.20 &  & 3658.71 & 0.25 & 3716.93 & 1553.47 & 1.84 \\ 
p05b & 240 & 4609.20 &  & 4613.58 & -- &  & 4611.16 & 483.60 &  & 4613.45 & 0.09 & 4671.35 & 3276.20 & 1.35 \\ 
p06b & 288 & 4716.36 &  & 4788.39 & -- &  & 4729.96 & 917.40 &  & 4780.56 & 1.36 & 4828.04 & 4148.33 & 2.37 \\ 
p07b & 72 & 1837.94 &  & \textbf{1837.94} & -- &  & \textbf{1837.94} & 30.60 &  & \textbf{1837.94} & 0.00 & \textbf{1837.94} & 69.52 & 0.00 \\ 
p08b & 144 & 3144.91 &  & 3149.77 & -- &  & 3149.77 & 129.00 &  & \textbf{3144.91} & 0.00 & 3153.64 & 446.55 & 0.28 \\ 
p09b & 216 & 3883.94 &  & 3937.53 & -- &  & \textbf{3883.94} & 534.00 &  & 3900.17 & 0.42 & 3944.60 & 1191.17 & 1.56 \\ 
p10b & 288 & 4927.95 &  & 4996.72 & -- &  & 4932.40 & 721.80 &  & 5007.57 & 1.62 & 5070.12 & 2865.22 & 2.88 \\ 
\noalign{\smallskip}\hline\noalign{\smallskip}
\multicolumn{2}{l}{Average} &  &  & 0.56 & \multicolumn{1}{l}{} & \multicolumn{1}{l}{} & 0.1 & 328.95 &  & \multicolumn{1}{l}{} & 0.34 & \multicolumn{1}{l}{} & 1424.57 & 0.97 \\ 
\hline
\end{tabular}
} 
\label{tbl:Resultados.HH.SDepVRPTW}
\end{table}